\documentclass[11pt,leqno]{article}

\usepackage{amsmath,amsthm}
\usepackage {latexsym}
\usepackage{amssymb}
\usepackage[dvips]{epsfig}
\input{epsf.tex}

\newcommand{\les}{\lesssim}
\newcommand{\sppt}{\mbox{\rm sppt}}

\newcommand{\const}{\mbox{\rm const}}
\newcommand{\bse}{\begin{equation}}

\newcommand{\bea}{\begin{eqnarray}}
\newcommand{\eea}{\end{eqnarray}}
\newcommand{\be}{\begin{equation}}
\newcommand{\ee}{\end{equation}}

\newcommand{\spec}{{\rm spec}}

\newcommand{\half}{\frac{1}{2}}

\newcommand{\eps}{{\varepsilon}}
\newcommand{\R}{{\mathbb R}}
\newcommand{\C}{{\mathbb C}}

\newcommand{\Compl}{{\mathbb C}}

\newcommand{\calg}{\,{\mathfrak g}}
\newcommand{\calG}{{\mathcal G}}

\newcommand{\calL}{{\mathcal L}}
\newcommand{\calS}{{\mathcal S}}

\newcommand{\calM}{{\mathcal M}}

\newcommand{\sign}{\mbox{sign}}

\newcommand{\Arg}{{\rm Arg}}

\newcommand{\Laplace}{\frac{\triangle}{2}}
\newcommand{\Lapl}{\frac12{\triangle}}

\newcommand{\la}{\langle}
\newcommand{\ra}{\rangle}

\def\pr{\partial}

\def\nn{\nonumber}
\def\calge1{\calg_{\vec{e_1}}}

\def\bm{\left( \begin{array}{cc}}
\def\endm{\end{array}\right)}
\def\ker{{\rm ker}}
\def\Ran{{\rm Ran}}
\def\Hil{{\mathcal H}}
\def\Dom{{\rm Dom}}

\def\Reg{{\mathcal R}}
\def\Res{{\rm Res}}
\def\T{{\mathcal T}}
\def\vpsi{{\vec \psi}}
\def\vphi{{\vec \phi}}
\def\vf{\vec f}
\def\vg{\vec g}
\def\pa{\partial}

\def\nn{\nonumber}
\def\Xl{{\cal{X}}}
\def\Yl{{\cal{Y}}}
\def\veta{\vec\eta}

\newtheorem{theorem}{Theorem}[section]
\newtheorem{lemma}[theorem]{Lemma}
\newtheorem{defi}[theorem]{Definition}
\newtheorem{cor}[theorem]{Corollary}
\newtheorem{prop}[theorem]{Proposition}

\theoremstyle{remark}
\newtheorem{remark}[theorem]{Remark}

\def\ga{\gamma}

\numberwithin{equation}{section}

\begin{document}

\title{ Dispersive Analysis of Charge Transfer Models}
\author{ I. Rodnianski, W. Schlag and A. Soffer
\thanks{The research of I.R. was partially conducted during the period
he served as a Clay Mathematics Institute Long-Term  
Prize Fellow.  He was also supported in part by the NSF grant DMS-0107791.
W.S. was partially supported by the NSF grant DMS-0070538
and a Sloan Fellowship. A.S.  was partially supported 
by the NSF grant DMS-0100490.}}
\maketitle

\begin{abstract}
We prove $L^p$ estimates for charge transfer Hamiltonians,
including matrix and inhomogeneous generalizations; such
equations appear naturally in the study of multi-soliton
systems.
\end{abstract}

\section{Introduction}
This paper is devoted to the study of dispersive properties of the 
model corresponding to the time-dependent charge transfer Hamiltonian
$$
H(t) = -\Lapl + \sum_{j=1}^{m} V_j (x - \vec v_j t)
$$
with rapidly decaying smooth potentials $V_{j}(x)$ and a set 
of mutually non-parallel constant velocities $\vec v_{j}$.
Our main focus is on the $L^{p}$ decay estimates for the
solutions of the time-dependent problem 
$$
\frac 1i \partial_{t }\psi + H(t) \psi = 0
$$
associated with a charge transfer Hamiltonian $H(t)$.
The well-known $L^p$ estimates for the free 
Schr\"odinger equation ($H(t)=-\Lapl$) on $\R^n$ are 
$$
\| e^{it\Laplace} f \|_{L^p}\leq C_p\, |t|^{-n(\frac{1}{2} -
\frac{1}{p})}\|f\|_{L^{p'}},\quad p\ge 2,\,\,\,
\frac{1}{p}+\frac{1}{p'}=1.
$$
They imply the Strichartz estimates
$$
\| e^{it\Laplace} f \|_{L^q_t L^r_x}\leq C_q\|f\|_{L^2},\quad  
2 \leq r, q
\leq \infty,\,\,\,\, \frac{n}{r} + \frac{2}{q} =\frac{n}{2},\quad n\geq 
3\qquad[ GV, KT]
$$
Such estimates play a fundamental role, among other things in the
theory of nonlinear dispersive equations.  The extension of such
theories to inhomogeneous problems (either due to curvature, local
potentials or coherent structure such as solitons, vortices etc. ...) 
then
motivated the efforts to establish the $L^p$ decay estimates for the general
{\it time independent} Schr\"odinger operators of the type 
$H= -\Lapl + V(x)$. 
In this case the estimates should take the form 
\begin{align}
&\|e^{-it H} P_{c}(H) \psi_{0}\|_{L^p}\leq C_p\, |t|^{-n(\frac{1}{2} -
\frac{1}{p})}\|\psi_{0}\|_{L^{p'}}\quad p\ge 2,\,\,\,
\frac{1}{p}+\frac{1}{p'}=1,\label{eq:inddec}\\
&\| e^{-itH} P_{c}(H)f \|_{L^q_t L^r_x}\leq C_q\|f\|_{L^2} \text{ for } 
2 \leq q
\leq \infty,\,\, \frac{n}{r} + \frac{2}{q} =\frac{n}{2},\quad n\geq 
3\label{eq:indStr}
\end{align}
where $P_{c} (H)$ is the projection onto the continuous part of the
spectrum of the self-adjoint operator $H$. Its purpose is to prevent
the emergence of any bound state, i.e., an $L^{2}$ eigenfunctions of $H$.
Under the evolution $e^{-itH}$ such bound states are merely  multiplied by 
 oscillating factors and thus do not disperse. 

The first approach to the proof of estimates \eqref{eq:inddec} was 
developed by
Journ\'e, Soffer, Sogge~\cite{JSS}. They  used a time dependent method
which  combined spectral and scattering theory with harmonic analysis. 
Their method involved splitting solutions into  high and low energy
parts and using  Kato's smoothing and the local energy decay on the 
corresponding pieces.
Later, a stationary method was
used by Yajima~\cite{Ya1} who proved that the wave operators are 
$L^p$-bounded under
weaker assumptions on the potential than in~\cite{JSS}. 
His theorem implies the dispersive bounds.
See also Weder~\cite{We} for  results in one 
dimension $n=1$ and Yajima~\cite{Ya3} for $n=2$.

Most recently, Rodnianski 
and Schlag~\cite{RS} have addressed the issue of determining the
optimal class of potentials $V(x)$ for which the 
estimates~\eqref{eq:inddec} and~\eqref{eq:indStr} 
hold true. In particular, they  obtained $L^1\to L^\infty$  decay 
estimates for the scaling invariant ``$L^{p}$-like" classes of 
small potentials
in dimension $n=3$.
 In \cite{NS}, Nier and Soffer were able to establish the
decay and Strichartz estimates 
for finite rank perturbations of $H_{0}=-\Laplace$. 
The Strichartz estimates for the inverse-square potential 
were obtained by  Burq, Planchon, Stalker, and Tahvildar-Zadeh 
in \cite{BPST} .

The situation in the case of  time-dependent potentials $V(t,x)$
becomes more complicated. On a superficial level, this can be 
attributed to the disappearance of a connection between the decay
estimates and spectral theory of a corresponding Hamiltonian.
In particular, solutions that remain trapped in a compact 
region, and thus do not disperse, might still exist. However, they 
can no longer be easily  characterized as corresponding to the 
bound states of some Hamiltonian as in the case of time independent 
Hamiltonians.

\noindent   Let $U(t)$ be the solution operator corresponding 
to the time-dependent Schr\"odinger equation
\bea
&&\frac 1i \partial_{t}\psi -\Laplace \psi + V(t,x) \psi =0, 
\label{eq:timeSc}\\
&&\psi|_{t=0} =\psi_{0},\nn
\eea
i.e., $\psi(t,\cdot) = U(t) \psi_{0}$. Then the desired decay 
estimates take the form 
\begin{equation}
\|U(t) \psi_{0}\|_{L^p}\leq C_p\, |t|^{-n(\frac{1}{2} -
\frac{1}{p})}\|\psi_{0}\|_{L^{p'}}\quad p\ge 2,\,\,\,
\frac{1}{p}+\frac{1}{p'}=1\label{eq:tddec}
\end{equation}
holding for an appropriate set of initial data.
The first results in this direction were obtained in \cite{RS},
where the $L^1\to L^\infty$ decay estimates for {\it all} 
initial data were established in dimension $n=3$ 
for a large class  of small time-dependent potentials. 
The smallness assumption rules out the emergence of 
trapped solutions.
 
\vskip 1pc

\noindent In this paper we prove $L^p$ decay estimates for 
solutions of a time-dependent Schr\"odinger equation corresponding
to a charge transfer Hamiltonian
$$
H(t)= -\Lapl + \sum_j V_j (x - \vec v_j t), \qquad 
\vec v_j \neq \vec v_i, \,\, i\neq j,
$$
see Theorem~\ref{thm:main} below. 
Our main motivation comes from the problem of asymptotic stability
of noninteracting multi-soliton states. This question refers to 
solving a NLS
\begin{equation}
\label{eq:NLS}
 i\partial_t \psi + \Laplace\psi + \beta(|\psi|^2)\psi = 0
\end{equation}
in $\R^{n},\, n\ge 3$ 
with initial data $\psi_0 = \sum_{j=1}^k w_j(0,\cdot) + R_0$ where 
$w_j$ are special standing wave
solutions called solitons and~$R_0$ is a small perturbation.
In our forthcoming paper~\cite{RSS}  we show that
if the solitons are sufficiently separated at time~$t=0$, and 
if~$R_0$ is sufficiently small
in a suitable norm, then  the
solution~$\psi$ evolves like 
a sum of solitons with time-dependent parameters approaching a limit, 
plus a radiation term that goes to zero in~$L^\infty(\R^3)$. 
This argument requires 
linearizing the NLS around the bulk-term,
which is given by the sum of the solitons. This leads to the problem 
of establishing dispersive
estimates, which typically means $L^1\to L^\infty$ decay estimates, 
for the linear problem.
It is well-known that the linearized equation needs to be written as 
a system in $\delta\psi$
 and~$\overline{\delta\psi}$ where $\delta{\psi}$ is the variation of 
the bulk-term.
A first step in the understanding of this problem is to consider the 
{\em scalar}  model that is closest
to the system at hand. It is easy to see that these scalar equations 
are Schr\"odinger equations
of the charge-transfer type
\begin{equation}
\label{eq:scal_mod}
 \frac 1i\partial\psi - \Laplace\psi +\sum_{j=1}^m V_j(\cdot-t\vec v_j)\psi 
=0,
\end{equation}
with real potentials $V_j$ and distinct vectors $\vec v_j\in\R^n$. 

This problem has been extensively studied in the literature in 
connection with the question of asymptotic completeness.  Asymptotic 
completeness for the equation \eqref{eq:scal_mod} 
had been first established by Yajima~\cite{Ya2}. Later,
Graf~\cite{Gr} has obtained a new proof containing a crucial 
argument showing that the energy\footnote{Energy here simply means 
$\|\psi(t)\|_{H^{1}}$ for a solution~$\psi(t)$. While the conservation of 
the $L^{2}$ norm is a trivial consequence of the potentials $V_{j}$ being
real, the energy bound is by no means a simple matter since the problem
contains a time-dependent potential.} associated with a solution of 
\eqref{eq:scal_mod} remains bounded in time.  For similar results see
also the work of W\"uller~\cite{Wu}, and for more general charge transfer models
with PDO's, the work of Zielinski~\cite{Z}. Their results imply that for a dense 
set of initial data, solutions converge to the states represented by 
radiation and the sum of the bound states, associated with each of the 
Hamiltonians $H_{j} = -\Laplace + V_{j}(x)$, traveling with velocity $\vec v_{j}$.

In this paper we obtain
$L^p$ decay estimates for~\eqref{eq:scal_mod}, which in particular lead to a very
short and self-contained proof of the asymptotic completeness theorem. 
The idea is to decompose the solution into $k+1$ channels that move 
(and expand) along
with the potentials according to a straight motion. On the $j^{\rm th}$~channel 
one has well-known dispersive estimates by Journ\'e, Soffer, 
Sogge~\cite{JSS} and 
Yajima~\cite{Ya1} for the equation with the single potential 
\footnote{The results of \cite{JSS} and \cite{Ya1} apply to 
problems with time-independent potentials. However, the problem
with a {\it single} potential $V_{j}(x-\vec v_{j}t)$ can be reduced 
to a problem with the potential $V_{j}(x)$ by means of a Galilei
transform.}
$V_j(\cdot-\vec v_j t)$, which dominates on that channel. 
One then proceeds to show that the interaction between the channels 
is weak enough not
to destroy this property. This weak interaction is accomplished by 
distinguishing low from high
momenta. In the former case one uses the slow propagation rate, 
in 
the latter a variant of
Kato's smoothing bound, to derive the desired bounds. As always, we 
cannot hope to prove
dispersive bounds for all initial conditions. In the scalar equation 
with a single stationary
potential this simply means being perpendicular to the bound states. 
In the charge
transfer case one expects this to be replaced with perpendicularity 
to the (traveling)
bound states, or at least that the projection of the solution onto 
the traveling bound states should go
to zero. This is precisely the "asymptotic orthogonality condition" 
which we use below. 

For this type of initial data we establish the decay estimates 
\begin{equation}
\label{eq:chdec}
\|\psi(t)\|_{L^{\infty}}\le C(n)\, |t|^{-\frac n2} \|\psi_{0}\|_{L^{1}\cap 
L^{2}},
\end{equation}
which, by interpolation with a trivial $L^{2}$ bound, also imply the 
the full range of the $L^{p}$ estimates
\begin{equation}
\label{eq:chdecp}
\|\psi(t)\|_{L^{p}}\le C(n,p)\, |t|^{-n(\frac 12 -\frac 1p)}
\|\psi_{0}\|_{L^{p'}\cap L^{2}},\quad 2\le p\le \infty
\end{equation}
The presence of the $L^{2}$ norm on the right-hand side 
of \eqref{eq:chdec} and \eqref{eq:chdecp} distinguishes them from
the standard $L^{p'}\to L^{p}$ estimates. This, however, is no more than
a technical issue. In fact, our main result is stated as an 
$L^{1}\cap L^{2}\to L^{2}+L^{\infty}$ estimate, which, for example, is
sufficient for a proof of asymptotic completeness.
We then establish a procedure converting this estimate into 
an $L^{1}\cap L^{2}\to L^{\infty}$ bounds, using an argument 
along the lines of \cite{JSS}.

\vskip 1pc

To motivate the results of the second part of this paper, we recall
that stability and asymptotic stability of solitons and other coherent
structures require analyzing the spectral properties, scattering,
and local decay of the linearized matrix operator. A comprehensive
description of results in this direction can be found in the book by 
C.~Sulem and P.\ L.\ Sulem~\cite{SuSu}.
 In particular, the spectral properties of the linearized operators 
for NLS solitons were developed in Weinstein~\cite{W1}, Shatah, Strauss~\cite{ShSt},
and Grillakis, Shatah, Strauss~\cite{GrShSt}. The linearized theory of other
structures was also considered, e.g., vortices in Weinstein, Xin~\cite{WeinXi}, 
and Hartree-type solitons in Fr\"ohlich, Tsai, Yau~\cite{FTY}. 
Many of these results are examples to which our $L^p$ theory applies.

On a more technical level, 
the second part of this paper extends our approach to the full system 
that one obtains by linearizing~\eqref{eq:NLS} around a sum of solitons. 
Several complications arise,
even for a system with a single stationary matrix potential:
Firstly, these systems are no longer self-adjoint and there is the 
well-known issue
of existence of unstable modes, see Weinstein~\cite{W1}. Secondly, 
and somewhat related to the 
aforementioned instabilities, the $L^1\to L^\infty$ decay estimates 
are more involved for systems (\cite{JSS} and~\cite{Ya1} do not 
apply). 

As far as the latter is concerned, Cuccagna~\cite{Cuc} recently 
established the necessary bounds by 
transferring Yajima's $L^p$-wave operator approach to the system 
case. While the wave operators
are of course of independent interest, we feel that it is desirable 
to have a more direct 
approach to the dispersive estimates for systems. It turns out that 
one can easily adapt
a simple but powerful method of Rauch~\cite{R} to systems, which is 
what we do in Section~\ref{sec:compl}
below, see Theorem~\ref{thm:vec_dec}.
 Rauch's method, which is related to earlier work of Dolph, 
McLeod,  Thoe~\cite{DMT} and 
Vainberg~\cite{V}, hinges on analytic continuation of the resolvent 
across the spectrum as
a bounded operator in an exponentially weighted~$L^2$-space. This can 
only be done for exponentially
decreasing potentials, but the potentials arising from linearizing 
around solitons are of this type
(since solitons decay exponentially). It is important to realize, 
however, that  Rauch's method
only leads to local~$L^2$-decay, which is inadequate for the 
nonlinear applications. To obtain
$L^\infty$-bounds, we then use an observation of Ginibre~\cite{Gin} 
that allows one to pass from local $L^2$-decay 
to $L^\infty$-decay by means of a double application of Duhamel's formula.

After deriving these bounds for systems with a single stationary 
matrix potentials, we 
then introduce charge transfer models in the system case, see Definition~\ref{def:chargetransm}.
 These are 
of course
analogous to the scalar models~\eqref{eq:scal_mod} but differ in that 
one not only
needs to translate the potentials according to linear motion, but one 
also needs to modulate
the matrix potentials in a suitable manner. In particular, these matrix 
potentials have complex off-diagonal
entries. We then proceed to apply the aforementioned method of 
"channel-decomposition" to these matrix charge transfer equations,
which yields bounds of the form~\eqref{eq:chdec}. 
Finally, we extend the $L^p$ decay estimates to the inhomogeneous
equation with charge transfer Hamiltonians. We give a general $L^p$-decay
estimate in terms of a properly chosen norm of the inhomogeneous term, 
see Theorem~\ref{thm:inhom}.

\section{Charge transfer model}

\begin{defi}
\label{def:chargetrans} 
By a {\em charge transfer model} we mean a Schr\"odinger equation
\bea
&& \frac{1}{i} \partial_t \psi - \Lapl \psi + \sum^m_{\kappa =1}
V_\kappa(x - \vec v_\kappa t) \psi = 0 \label{eq:transfer} \\
&& \psi |_{t=0} = \psi_0,  x\in \R^n, \nn
\eea
where $\vec v_\kappa$ are distinct vectors in $\R^n$, $n\ge3$, and 
the real potentials $V_\kappa$ are such that
for every $1\le\kappa\le m$, 
\begin{enumerate}
\item $V_\kappa$ has compact support (or fast decay), 
$V_{\kappa}, \nabla V_\kappa \in L^\infty$ 
\item  $0$ is neither a zero eigenvalue nor a zero resonance of
the operators
$$
H_\kappa = - \Lapl + V_\kappa(x).$$
\end{enumerate}
\end{defi}

The regularity assumption on the potential can be relaxed
considerably without too much effort, but we use the $C^1$-property
for convenience. 
The assumption $n\ge3$ will be made throughout. This definition is 
standard, see~\cite{Gr}, \cite{Ya2}. 
The Schr\"odinger group
$e^{-itH_\kappa}$ is known to satisfy the decay estimates (see 
Journ\'e, Soffer, Sogge~\cite{JSS} and Yajima~\cite{Ya1}) 
\begin{equation}
\| e^{-itH_\kappa} P_c (H_\kappa) \psi_0 \|_{L^\infty} \lesssim
|t|^{-n/2} \|\psi_0\|_{L^1}. \label{eq:nodecay}
\end{equation}
Here $P_c(H_\kappa)$ is the spectral projection onto the continuous
spectrum of $H_\kappa$. The estimate \eqref{eq:nodecay} holds under 
the "no 0 eigenvalue/resonance" assumption 2. above for potentials 
obeying conditions roughly equivalent to the assumption 1.
\begin{remark}
It follows from the results of Yajima \cite{Ya1} that in dimension
$n=3$ the estimate \eqref{eq:nodecay} holds for any bounded potential
with sufficiently fast decay at infinity. In higher dimensions one 
needs to add assumptions on the derivatives of the potential 
$V_{\kappa}$:
\begin{align}
&|\nabla^{\alpha} V_{\kappa} (x)|\le C_{\alpha} \la x\ra^{-\delta},
\qquad \forall |\alpha|\le \frac {n+1}2-3,\quad 
\delta > \max\{n+2, \frac {3n}2 -2\},\label{eq:ya1}\\
&{\mathcal{F}} (\la x\ra^{\sigma} V ) \in L^{\frac {n-1}{n-2}},\qquad 
\sigma > \frac {2(n-2)}{n-1}\label{eq:ya2}
\end{align}
We shall assume that potentials $V_{\kappa}(x)$ always obey 
conditions guaranteeing estimate \eqref{eq:nodecay}, e.g.
conditions \eqref{eq:ya1}, \eqref{eq:ya2}.
\end{remark}

\begin{defi}
\label{def:2+infty} 
Let $L^2 + L^\infty : = \{f : \exists
h; g : f = h + g, \| h \|_{L^2} + \| g \|_{L^\infty} < \infty\}$
with norm
$$
\| f \|_{L^2 + L^\infty} := {\inf}_{f=h+g} (\| h \|_{L^2} + 
\|g\|_{L^\infty})
$$
\end{defi}
We shall use a weaker version of the decay estimate
which follows from~\eqref{eq:nodecay} combined with the unitarity of 
$e^{-itH_\kappa}$, namely
\begin{equation}
\| e^{-itH_\kappa} P_c (H_\kappa) \psi_0 \|_{L^\infty + L^2} \lesssim
\langle t\rangle^{-n/2} \| \psi_0 \|_{L^1 \cap L^2} \label{eq:decay}
\end{equation}
where $\langle t\rangle = (1+t^2)^{\half}$. 
In the following we shall assume that the
number of potentials is $m = 2$ and that the velocities are 
$\vec v_1 = 0,
\vec v_2 = (1,0, \dots 0) = \vec e_1$. 
This can be done without loss of generality. An indispensable tool in 
the study 
of the charge transfer model are the Galilei transforms
\begin{equation}
\label{eq:gal} 
\calg_{\vec{v},y}(t)= e^{-i\frac{|\vec v|^2}{2} t} e^{-i{x\cdot\vec 
v}} e^{i(y+t\vec v)\cdot\vec p},
\end{equation}
cf.~\cite{Gr}, where $\vec p=-i\vec \nabla$. These are the quantum 
analogues of the classical
Galilei transforms 
\[ x\mapsto x-t\vec v -y,\quad \vec p\mapsto \vec p-\vec v\]
in the following sense: if $f$ is a Schwartz function, say,  such 
that $f$ and $\hat{f}$ are supported around $0$,
 then $\calg_{\vec{v},y}(t) f$ is supported around $t\vec v+y$, and 
$\widehat{\calg_{\vec{v},y}(t)f}$ is supported around~${\vec v}$.
Under  $\calg_{\vec{v},y}(t)$ the Schr\"odinger equation transforms 
as follows:
\[
 \calg_{\vec v,y}(t) e^{it\Laplace} = e^{it\Laplace} \calg_{\vec 
v,y}(0)
\]
and moreover, with $H=-\Lapl + V$, 
\begin{equation}
\label{eq:invg}  \psi(t) := \calg_{\vec v,y}(t)^{-1} e^{-itH} 
\calg_{\vec v,y}(0) \phi_0, ,\qquad \calg_{\vec v,y}(t)^{-1} = 
e^{-i{y\cdot \vec v}}\calg_{-\vec v,-y}(t)
\end{equation}
solves
\bea
&& \frac{1}{i} \partial_t \psi - \Lapl \psi + V(\cdot - t\vec v -y) 
\psi = 0 \label{eq:shifted} \\
&& \psi |_{t=0} = \phi_0. \nn
\eea
These properties are not only easy to check, but are of course to be 
expected in view of
the classical interpretation. We will make frequent use of these 
transformation laws without
further notice. Another property that we will use often without 
further mention is that 
the transformations $\calg_{\vec v,y}(t)$ are isometries on all $L^p$ 
spaces.
Finally, since in our case always $y=0$, we set $\calg_{\vec v}(t):= 
\calg_{\vec v,0}(t)$. 
By~\eqref{eq:invg}, $\calge1(t)^{-1}=\calg_{-\vec{e_1}}(t)$ in that 
case. 

\noindent We now return to the problem
\bea
&& \frac{1}{i} \partial_t \psi - \Lapl \psi + V_1 \psi +V_2
(\cdot -t\vec{e_1}) \psi = 0 \label{eq:2pot} \\
&& \psi|_{t=0} = \psi_0 \nn 
\eea
with $V_1, V_2$ compactly supported potentials. 
Let $u_1, \ldots, u_m$ and $w_1, \ldots, w_\ell$ be the normalized 
bound states of~$H_1$
and $H_2$ corresponding to the negative eigenvalues $\lambda_1, \dots,
\lambda_m$ and $\mu_1, \dots, \mu_\ell$, respectively 
(recall that we are assuming that $0$ is not an eigenvalue). 
We denote by $P_b(H_1)$ and $P_b(H_2)$ the 
corresponding projections onto 
the bound states of $H_1$ and $H_2$, respectively, and let 
$P_c(H_\kappa)=Id-P_b(H_\kappa)$, $\kappa=1,2$.
The projections $P_b(H_{1,2})$ have the form 
$$
P_b(H_1) =\sum^m_{i=1} \langle \cdot, u_i\rangle u_i, \qquad 
P_b(H_2) = \sum^\ell_{j=1} \langle \cdot, w_j\rangle w_j.
$$ 
In order to state our main theorem, we need to impose an orthogonality
condition in the context of the charge transfer 
Hamiltonian~\eqref{eq:2pot}.

\begin{defi} 
\label{def:asymp}  Let $U(t) \psi_0 = \psi(t, x)$ be the
solutions of~\eqref{eq:2pot}. We say that $\psi_0$ (or also $\psi(t,\cdot)$) 
is {\em asymptotically orthogonal} to the bounds states of $H_1$ and $H_2$ if 
\begin{equation}
\label{eq:scat}
\|P_b(H_1)U(t) \psi_0 \|_{L^2} + \| P_b(H_2,t) U(t)\psi_0 \|_{L^2}\to 
0\text{ as }t\to +\infty.
\end{equation}
Here
\begin{equation}
\label{eq:Proj2}
P_b(H_2,t) := \calg_{\vec{-e_1}}(t) P_b(H_2)\, \calge1(t)
\end{equation}
for all times $t$. 
\end{defi}

\begin{remark} Clearly, $P_b(H_2,t)$ is again an orthogonal 
projection for every~$t$.
It gives the projection onto the bound states of $H_2$ that have been 
translated to the
position of the potential $V_2(\cdot-t\vec{e_1})$. Equivalently, one 
can think of it as
translating the solution of~\eqref{eq:2pot} from that position to the 
origin, projection 
onto the bound states of~$H_2$, and then translating back. 
The explicit form of $ P_b(H_2,t)$ is given by 
$$\big(P_b(H_2,t)f\big)(x) =\sum^\ell_{j=1} e^{i{x_1}} w_j(x - 
t\vec{e_1}) \int_{\R^n} f(y)
e^{-i{y_1}} \overline{w_j(y- t\vec{e_1})} \, dy.
$$
We will make little use of this formula, though.
\end{remark}

\begin{remark}
\label{rem:graf}
It is clear that all $\psi_0$ that satisfy~\eqref{eq:scat}
form a closed subspace. This subspace coincides with the
space of {\em scattering states} for the charge transfer problem.
The latter is well-defined by Graf's asymptotic 
completeness result~\cite{Gr}. This is discussed in more 
detail in Section~\ref{sec:asymp}.
\end{remark}

\noindent We now formulate our main result.

\begin{theorem}
\label{thm:main}
Consider the charge transfer model as in 
Definition~\ref{def:chargetrans} with two potentials, 
cf.~\eqref{eq:2pot}. Let $U(t)$ denote the propagator of the 
equation~\eqref{eq:2pot}. Then 
for any initial data $\psi_0 \in L^1\cap L^2$, which is 
asymptotically orthogonal to the bound
states of $H_1$ and $H_2$ in the sense of Definition~\ref{def:asymp}, 
one has the 
decay estimates
\begin{equation}
\| U(t) \psi_0 \|_{L^2 + L^\infty} \lesssim
\langle t\rangle^{-n/2}\|\psi_0\|_{L^1\cap L^2}. \label{eq:main}
\end{equation}
An analogous statement holds for any number of potentials, i.e., with 
arbitrary $m$ in~\eqref{eq:transfer}.
\end{theorem}

We prove this theorem for the case of two potentials, but this is for 
simplicity only.
The same argument also applies to the general case of $m>2$ (one then 
needs to split into $m+1$ channels,
see below). A more substantial comment has to do with  the assumption 
of compact support
of the potentials. Inspection of the argument in the following 
section reveals that it equally
well applies to exponentially decaying potentials, say. But also 
power decay is allowed, provided
it is sufficiently fast. 
We shall prove~\eqref{eq:main} by means of a bootstrap argument. More 
precisely, 
we prove that the {\em bootstrap assumption}
\begin{equation}
\| U(t) \psi_0\|_{L^2+ L^\infty }\leq C_0 \langle t\rangle^{-n/2}
\|\psi_0\|_{L^1\cap L^2} \text{\ \ for all\ \ }0\le t\le 
T\label{eq:boot}
\end{equation}
implies that
\begin{equation}
\|U(t) \psi_0\|_{L^2+ L^\infty} \leq \frac{C_0}{2} \langle 
t\rangle^{-n/2}
\|\psi_0\|_{L^1\cap L^2} \text{\ \ for all\ \ }0\le t\le T. 
\label{eq:boot/2}
\end{equation}
Here $C_0$ is some sufficiently large positive constant which is 
assumed to be
much bigger than any universal constant appearing in our calculations
(e.g.~$\|V_{1,2} \|_{L^1\cap L^\infty}$, implicit constants in
Proposition~\ref{prop:bdstates} below, and decay 
bounds~\eqref{eq:decay}).
The logic here is that for arbitrary but fixed $T$, the 
assumption~\eqref{eq:boot}
can be made to hold for some $C_0$, depending on~$T$. Iterating the 
implication
\eqref{eq:boot}~$\Longrightarrow$~\eqref{eq:boot/2} then yields a 
constant that
does not depend on~$T$, thus proving the theorem.
We shall continue to use the notation $\lesssim$ to denote bounds 
involving
multiplicative constants independent of the constant~$C_0$.

\begin{remark}\label{rem:larget}
Since 
$$\|U(t) \psi_0\|_{L^2+ L^\infty} \leq
\|U(t)\psi_0\|_{L^2} = \|\psi_0 \|_{L^2} \leq \|\psi_0\|_{L^1\cap L^2}
$$
it suffices to prove \eqref{eq:boot/2} for $t\geq t_0 := 
\left(\frac{C_0}{2}\right)^{2/n}$.  
\end{remark} 

\begin{remark}
\label{rem:perturb}
For the nonlinear applications it might be also useful to have
estimates as in Theorem~\ref{thm:main} for 
\underline{perturbed} charge transfer Hamiltonians. It is easy to see
that the method of proof is stable under such small perturbations,
see Section~\ref{subsec:perturb} below.
\end{remark}

\section{Proof of the decay estimates}
\label{sec:dec}

\subsection{Bound states}

Our first result concerns the rate of convergence of the projections
onto the bound states of solutions which are asymptotically orthogonal
to the bound states.

\begin{prop}
\label{prop:bdstates}
Let $\psi (t, x) = (U(t) \psi_0)(x)$ be a solution of~\eqref{eq:2pot} 
which is
asymptotically orthogonal to the bound states of $H_1$ and $H_2$ in 
the
sense of Definition~\ref{def:asymp}. Then 
$$
\| P_b (H_1) U(t)\psi_0 \|_{L^2} + \| P_b (H_2,t) U(t) \psi_0\|_{L^2}
\lesssim e^{-\alpha t} \| \psi_0\|_{L^2}
$$
for some $\alpha>0$.
\end{prop}
\begin{proof}
By symmetry it suffices to prove the bound on the first part. More 
precisely,
let $\tilde U(t):= \calge1(t) U(t)$, and $\phi(t) = \tilde 
U(t)\psi_0$. Then $\phi(t)$
solves
\bea
&& \frac{1}{i} \partial_t \phi - \Lapl \phi + V_1(\cdot+t\vec{e_1}) 
\phi +V_2 \phi = 0, \label{eq:dual} \\
&& \phi|_{t=0}(x) = (\calge1(0) \psi_0)(x), \nn 
\eea
and
\[ \|P_b(H_2,t)U(t)\psi\|_2 = \|P_b(H_2)\calge1(t) U(t)\psi_0\|_2 = 
\|P_b(H_2) \tilde U(t) \psi_0\|_2.\]
This clearly allows one to reduce the treatment of $H_2$ to that 
of~$H_1$.
Decompose 
\begin{equation}
U(t) \psi_0 = \sum^m_{i=1} a_i (t) u_i + \psi_1 (t, 
x)\label{eq:decomp}
\end{equation}
relative to $H_1$ so that $\psi_1(t,\cdot)$ lies in the continuous
subspace of $H_1$, i.e.,  $P_c(H_1) \psi_1 = \psi_1$ and 
$P_b(H_1)\psi_1 = 0.$
By assumption, 
$$
\sum^m_{i=1} |a_i(t)|^2 \to 0 \text{ as } t \to \infty.
$$
Substituting \eqref{eq:decomp} into \eqref{eq:2pot} yields
\begin{align}
&\frac{1}{i}\partial_t \psi_1 - \Lapl \psi_1 + V_1 \psi_1 + V_2(\cdot 
-t\vec{e_1}) \psi_1 +\nn\\
&+  \sum^m_{j=1} \left[\frac{1}{i} \dot a_j(t) u_j - \Lapl u_j
a_j (t) + V_1\, u_j a_j(t) +  V_2(\cdot - t\vec{e_1}) u_j 
a_j(t)\right]=0.
\label{eq:ai}
\end{align}
Since $P_c(H_1) \psi_1 = \psi_1$ one has 
$$
(-\Lapl + V_1)\psi_1 = H_1 \psi_1 = P_c(H_1) H_1 \psi_1,\qquad 
\partial_t \psi_1 = P_c(H_1) \partial_t
\psi_1.
$$
In particular 
$$P_b(H_1) \left(\frac{1}{i}\partial_t \psi_1 -
\Lapl\psi_1 + V_1 \psi_1\right) =0.
$$
Thus taking an inner product of the equation~\eqref{eq:ai} with 
$u_\kappa$ and
using the fact that $\langle u_\kappa , u_j\rangle = 
\delta_{j\kappa}$ as
well as the identity
$$
-\Lapl u_j + V_1 u_j = \lambda_j u_j
$$
we obtain the ODE
$$
\frac{1}{i} \dot a_\kappa(t) + \lambda_\kappa a_\kappa(t) + \langle
V_2(\cdot - t \vec{e_1}) \psi_1, u_\kappa \rangle + \sum_{j=1}^m 
a_j(t) \la V_2(\cdot-t\vec{e_1})u_j,u_\kappa \ra = 0 \label{eq:ODE}
$$
for each $a_\kappa$ with the condition that 
$$
a_\kappa(t) \to 0 \text{ as } t \to + \infty.
$$
Recall that $u_\kappa$ is an eigenfunction of $H_1=-\Lapl + V_1$ with 
eigenvalue $\lambda_\kappa < 0$.
It is well-known (see e.g.~Agmon~\cite{Ag}) that such eigenfunctions 
are
exponentially localized, i.e.,
\begin{equation}
\label{eq:local}
\int_{\R^n} e^{2\alpha |x|}\, |u_\kappa(x) |^2\, dx  \le C= C(V_1,n) 
<\infty \text{ for some positive }
\alpha.
\end{equation}
Therefore, 
\begin{equation}
\|V_2(\cdot - t \vec{e_1}) u_\kappa \|_2 \lesssim e^{-\alpha t} 
\text{\ \ for all\ \ }t\ge 0
\label{eq:V2uk}
\end{equation}
which follows from the assumption that $V_2$ has compact support.
The implicit constant in~\eqref{eq:V2uk} depends on the size of the 
support of
$V_2$ and $\| V_2\|_{L^\infty}$. Therefore, 
\[ f_\kappa(t) : =  \la V_2(\cdot - t\vec{e_1})\psi_1, u_\kappa \ra \]
satisfies 
\begin{equation}
|f_\kappa(t)| \lesssim  e^{-\alpha t} \|\psi_1\|_{L^2} 
\les e^{-\alpha t} \|P_c(H_1)U(t)\psi_0\|_{L^2} \lesssim e^{-\alpha 
t} \|\psi_0\|_2.
\label{eq:fkappa}
\end{equation}
In view of~\eqref{eq:ODE},  $a_\kappa$ solves the equation
\bea
&& \frac{1}{i} \dot a_\kappa (t) + \lambda_\kappa a_\kappa(t) + 
\sum_{j=1}^m a_j(t) C_{j\kappa}(t) + f_\kappa(t) =0 \label{eq:ODE2} \\
&& a_\kappa(\infty) = 0,\nn 
\eea
where $C_{j\kappa}(t)=C_{\kappa j}(t)=\la 
V_2(\cdot-t\vec{e_1})u_j,u_\kappa \ra$. By~\eqref{eq:V2uk},
$\max_{j,\kappa} |C_{j\kappa}(t)| \lesssim e^{-\alpha t}.$
Solving \eqref{eq:ODE2} explicitly we obtain
$$
\vec a (t) = ie^{-i \int_0^t B(s)\,ds}\int^\infty_t e^{i\int_0^s 
B(\tau)\,d\tau} \vec f(s) \,ds,
$$
where $B_{j\kappa}(t)=\lambda_j\delta_{j\kappa}+C_{j\kappa}(t)$. 
Hence, by unitarity of $e^{i\int_0^s B(\tau)\, d\tau}$ 
and~\eqref{eq:fkappa}, we conclude that
$$
|\vec a(t) | \le \int_t^\infty |\vec f(s) |\, ds \lesssim \alpha^{-1} 
e^{-\alpha t} \|\psi_0 \|_{L^2}.
$$
This proves  the proposition.
\end{proof}

\subsection{The three channels}

\noindent Fix some small $\delta>0$. We introduce a partition of 
unity associated with the sets 
$$
B_{\delta t}(0)= \{ x: |x| \leq \delta t\},\quad B_{\delta t}(t 
\vec{e_1}) = \{ x:
|x - t\vec{e_1} | \leq \delta t\}
$$
and
$$
\R^n \setminus (B_{\delta t} (0) \cup B_{\delta t}(t \vec{e_1})).
$$
Namely, let $\chi_1(t,x)$ be a cut-off function s.t.
$$
\chi_1 (t,x) = 1,\quad  x \in B_{\delta t} (0) \text{\ \  and\ \ }
\chi_1 (t,x) = 0,\quad  x \in \R^n\setminus B_{2\delta t}(0).
$$
Here $t\geq t_0$ and $\delta>0$ is a fixed small constant.
Define
$$
\chi_2 (t, x) = \chi_1 ( t, x-t \vec e_1),\quad
\chi_3(t,x) = 1- \chi_1 (t, x) - \chi_2 (t, x).
$$
Observe that since $t_0 = \left(\frac{C_0}{2}\right)^{2/n}$ is large, 
see Remark~\ref{rem:larget}, the
support of $\chi_1 (t, \cdot)$ contains the support of $V_1$ for all 
$t\geq t_0$ 
and the support of $\chi_2(t, \cdot)$ contains the support of 
$V_2(\cdot - t\vec{e_1})$.  

\noindent There is the following natural decomposition of 
the solution $U(t) \psi_0$:
\bea
U(t)\psi_0 &=&   \chi_1(t,\cdot) P_b (H_1) U(t) \psi_0 + 
\chi_1(t,\cdot) P_c(H_1) U(t) \psi_0 +
\chi_2(t,\cdot) P_b (H_2,t) U(t) \psi_0 \nn \\
&& +
\chi_2(t,\cdot) P_c (H_2,t) U(t) \psi_0 + \chi_3(t,\cdot) U (t) 
\psi_0. \nn
\eea
The terms $\chi_1 P_b( H_1) U(t) \psi_0$ and $\chi_2 P_b(H_2,t) U(t) 
\psi_0$ 
can be estimated immediately using Proposition~\ref{prop:bdstates}.
Since  $U(t) \psi_0$ is asymptotically orthogonal to the
bound states of both $H_1$ and $H_2$ by assumption, one concludes 
that 
$$
\|\chi_1(t, \cdot) P_b (H_1) U(t) \psi_0 \|_{L^2+L^\infty} +
\|\chi_2(t,\cdot)  P_b (H_2,t) U(t) \psi_0 \|_{L^2 + L^\infty}
\lesssim e^{-\alpha t} \|\psi_0 \|_{L^1 \cap L^2}.\label{eq:easy}
$$
To estimate the remaining three terms 
\[ \chi_1(t,\cdot) P_c (H_1) U(t) \psi_0,\quad
\chi_2 (t, \cdot) P_c (H_2,t) U(t) \psi_0, \quad \chi_3 (t,\cdot) 
U(t) \psi_0\] 
we use Duhamel's formula.  However, on each of the three
``channels''  we compare $U(t) \psi_0$ to a different group.  We
claim that on the support of $\chi_1 (t, \cdot)$ the solution  $ U(t) 
\psi_0$ is best
approximated by $e^{-itH_1} \psi_0$, on the support of $\chi_2(t, 
\cdot)$ by
$e^{-itH_2} \psi_0$, and on $\chi_3 (t, \cdot)$ by $e^{it\Laplace} 
\psi_0$.
More precisely, with $\chi_1=\chi_1(t,x)$, 
\begin{equation}
\label{eq:DU1}
\chi_1 P_c (H_1) U(t) \psi_0 = \chi_1 e^{-it H_1} P_c(H_1)
\psi_0 -i\chi_1\int^t_0 e^{-i(t-s) H_1} P_c (H_1) V_2(\cdot-s 
\vec{e_1})U(s) \psi_0 \, ds.
\end{equation}
To analyze the second channel, define $\tilde 
U(t):=\calge1(t)U(t)\calg_{-\vec{e_1}}(0)$. Then 
$\phi(t):=\tilde U(t)\phi_0$ solves equation~\eqref{eq:dual}, where 
$\phi_0:=\calge1(0)\psi_0$. 
The Duhamel formula relative to $e^{-itH_2}$ can now be written in 
two different ways, namely with
respect to a stationary frame (this involves $\tilde U(t)$) or with 
respect to the
moving frame (this involves $U(t)$). Indeed, 
\bea
\label{eq:DU2stat}
\chi_1 P_c(H_2) \tilde U(t) \phi_0 &=& \chi_1 e^{-it H_2} P_c(H_2)
\phi_0 - i\chi_1\int^t_0 e^{-i(t-s) H_2} P_c (H_2) V_1(\cdot+s 
\vec{e_1})\tilde U(s) \phi_0 \, ds \\
\chi_2 P_c (H_2,t) U(t) \psi_0 &=& \chi_2 \calg_{-\vec{e_1}}(t) 
e^{-it H_2} P_c(H_2)\calge1(0)
\psi_0 \nn \\ 
&& \mbox{\hspace{.5in}} -i\chi_2\calg_{-\vec{e_1}}(t)\int^t_0 
e^{-i(t-s) H_2} P_c (H_2) V_1(\cdot+s \vec{e_1})\calge1(s)U(s) \psi_0 
\, ds.
\label{eq:DU2mov} 
\eea
The stationary formulation~\eqref{eq:DU2stat} is basically the same 
as~\eqref{eq:DU1}, 
whereas~\eqref{eq:DU2mov} follows from~\eqref{eq:DU2stat}  by 
means of~\eqref{eq:Proj2} from Definition~\ref{def:asymp}.
It is evident from~\eqref{eq:DU1} and~\eqref{eq:DU2stat} that it 
suffices to prove the
requisite amount of decay for~$\chi_1 P_c (H_1) U(t) \psi_0$. 
However, we will need~\eqref{eq:DU2mov}
in the treatment of the first channel, i.e.,~\eqref{eq:DU1}. 
Finally, for the third (the ``free'') channel, one has
\begin{equation}
\chi_3 U(t) \psi_0 = \chi_3\, e^{it\Laplace} \psi_0 - i\chi_3\int^t_0
e^{i(t-s)\Laplace} (V_1 + V_2(\cdot-s\vec{e_1}))U(s) \psi_0\, ds.
\label{eq:DU3}
\end{equation}

\subsection{Analysis of $\chi_1P_c(H_1) U(t) \psi_0$}

We start with the following simple lemma.

\begin{lemma} 
\label{lem:dual}  
Let $W\in L^p \cap L^{\frac{2p}{2-p}}$ for some
$p\in [1,2]$.  Then 
$$
\| W f \|_{L^p}\lesssim \| W \|_{L^p\cap L^{\frac{2p}{2-p}}} \| f
\|_{L^2 + L^\infty}.
$$
In particular, the dual space of $L^2+L^\infty$ is $L^1\cap L^2$. 
\end{lemma}

Fix sufficiently large constants $A, B> 0, B>>A$,
 which are independent of the constant $C_0$ and write the following
expansion
\bea
\chi_1 P_c(H_1) U(t) \psi_0 &=& \chi_1 e^{-itH_1} P_c (H_1)\psi_0 - i
\chi_1 \int^{t-A}_0 e^{-i(t-s) H_1} P_c (H_1) V_2(\cdot-s 
\vec{e_1})U(s) \psi_0\, ds \nn \\
&& \qquad -i\chi_1\int^t_{t-A} e^{-i(t-s)H_1} P_c(H_1) 
V_2(\cdot-s\vec{e_1}) U(s) \psi_0\, ds.
\label{eq:chan1exp}
\eea
According to~\eqref{eq:decay} 
\begin{equation}
\label{eq:20}
\| e^{-itH_1} P_c (H_1) \psi_0 \|_{L^2+L^\infty} \lesssim \langle
t\rangle^{-n/2} \|\psi_0 \|_{L^1 \cap L^2}.
\end{equation}
Using \eqref{eq:decay}, Lemma~\ref{lem:dual}, and the bootstrap 
assumption~\eqref{eq:boot}
\bea
&& \Bigl\| \int^{t-A}_0 e^{-i(t-s) H_1} P_c(H_1) V_2(\cdot-s 
\vec{e_1}) U(s)
\psi_0\, ds \Bigr\|_{L^2 + L^\infty} \nn \\
&& \lesssim
\int^{t-A}_0 \frac{1}{\langle t-s\rangle^{n/2}}\|V_2 
(\cdot-s\vec{e_1}) U(s)\psi_0 \|_{L^1\cap L^2}\, ds \nn \\
&& \lesssim \int^{t-A}_0 \frac{1}{\langle 
t-s\rangle^{n/2}}\Bigl[\|V_2\|_{L^1\cap L^2}+\|V_2\|_{L^2\cap 
L^\infty}\Bigr]  \| U(s) \psi_0\|_{L^2 + L^\infty}\, ds \nn \\
&&  \lesssim
C_0\int^{t-A}_A\frac{1}{\langle t-s\rangle^{n/2}}  \frac{1}{\langle
s\rangle^{n/2}}\, ds\; \|V_2 \|_{L^1\cap L^\infty} \| \psi_0 
\|_{L^1\cap L^2} \nn \\
&& \qquad + \int_0^A \frac{1}{\la t-s \ra^{n/2}} \, ds  \|V_2 
\|_{L^1\cap L^\infty} \| U(s) \psi_0\|_{2}
\nn \\
&& \lesssim C_0\,A^{-(n/2-1)} \langle t\rangle^{-n/2} \| 
V_2\|_{L^1\cap L^\infty} 
\|\psi_0 \|_{L^1\cap L^2} + A \la t\ra^{-\frac{n}{2}} \|\psi_0\|_2. 
\label{eq:At-A}
\eea
Thus, since $n\geq 3$, there exists a choice of a sufficiently large
constant $A$ that does not depend on $C_0$ such that
\begin{equation}
\Bigl\| \int^{t-A}_0 e^{-i(t-s)H_1} P_c(H_1) V_2(\cdot-s \vec{e_1}) 
U(s)
\psi_0\, ds \Bigr\|_{L^2+ L^\infty}\leq
 10^{-2} C_0 \langle t\rangle^{-n/2}\|\psi_0 \|_{L^1\cap L^2},
\label{eq:21}
\end{equation}
provided $C_0$ was chosen sufficiently large. 
Recall that  $t\geq t_0 =(C_0/2)^{2/n}$. We can therefore assume that 
$t>> A$.
The third term on the right hand-side of~\eqref{eq:chan1exp} is 
treated as follows.
With $\chi_1=\chi_1(t,x)$,  
\bea
&& \chi_1\int^t_{t-A} e^{-i(t-s)H_1} P_c (H_1) V_2(\cdot - s 
\vec{e_1})
U(s) \psi_0\, ds  \nn \\
&& = \chi_1 \int^t_{t-A} e^{-i(t-s)H_1} P_c (H_1) V_2(\cdot - 
s\vec{e_1})
 P_b(H_2,s) U(s) \psi_0\, ds \nn \\
&& \quad +\chi_1 \int^t_{t-A} e^{-i(t-s)H_1} P_c (H_1) V_2 
(\cdot-s\vec{e_1}) P_c (H_2,s) U(s) \psi_0\, ds = : J_b+ J_c. 
\label{eq:Jbcdef}
\eea
Since  $U(s)\psi_0$ is asymptotically orthogonal to the bound states 
of both $H_1$ and $H_2$ by assumption,  
Proposition~\ref{prop:bdstates} implies that
$$
\| P_b(H_2,s) U(s) \psi_0 \|_{L^2}\lesssim e^{-\alpha s} \| \psi_0
\|_{L^2}.
$$
Therefore, using unitarity of $e^{-i(t-s)H_1}$ and $L^2$ boundedness 
of the spectral projection $P_c(H_1)$ we obtain
\bea
\|J_b\|_{L^2 + L^\infty} &\leq& \| J_b\|_{L^2}\lesssim \int^t_{t-A}
\|V_2(\cdot-s\vec{e_1})  P_b(H_2,s) U(s) \psi_0 \|_{L^2}\, ds \nn \\
&\lesssim& \| V_2 \|_{L^\infty}\int^t_{t-A} \| P_b(H_2,s)U(s)
\psi_0 \|_{L^2}\,ds \nn \\
&\lesssim&
\| V_2 \|_{L^\infty} \int^t_{t-A} e^{-\alpha s} ds \|\psi_0\|_{L^2} 
\lesssim
\| V_2 \|_{L^\infty} A e^{-\alpha t/2}\|\psi_0\|_{L^1\cap L^2},
\eea
assuming as we may that $t-A\ge t/2$. 
This implies the bound 
\begin{equation}
\| J_b\|_{L^2 + L^\infty} \lesssim \langle 
t\rangle^{-n/2}\|\psi_0\|_{L^1\cap L^2}. 
\label{eq:23}
\end{equation}
To deal with the term $J_c$ in \eqref{eq:Jbcdef} we expand $\chi_2 
(s, \cdot)P_c(H_2,s) U(s) \psi_0$ 
via the Duhamel formula~\eqref{eq:DU2mov}. Then
\bea
J_c &=& \chi_1 \int^t_{t-A} e^{-i(t-s) H_1} P_c (H_1) 
V_2(\cdot-s\vec{e_1})\calg_{-\vec{e_1}}(s) e^{-isH_2} 
P_c(H_2)\calge1(0)\psi_0\, ds \nn \\
&& - i\chi_1\int^t_{t-A}\int^s_0 e^{-i(t-s)H_1} P_c(H_1) 
V_2(\cdot-s\vec{e_1})\calg_{-\vec{e_1}}(s) e^{-i(s-\tau) H_2} 
P_c(H_2) \nn \\
&& \mbox{\hspace{2.5in}} V_1(\cdot + \tau\vec{e_1}) \calge1(\tau) 
U(\tau) \psi_0\, d\tau ds \nn \\
&&  = : J_c^1 + J_c^2. \label{eq:Jc12def} 
\eea
Using \eqref{eq:decay} for $\kappa = 2$ and Lemma~\ref{lem:dual} we 
estimate $J_c^1$ as follows:
\bea
\| J_c^1 \|_{L^2 + L^\infty} &\leq& \| J_c^1 \|_{L^2} \lesssim 
\int^t_{t-A}
\Bigl\| V_2(\cdot - s \vec{e_1}) \calg_{-\vec{e_1}}(s) 
e^{-isH_2}P_c(H_2) \calge1(0) \psi_0 \Bigr\|_{L^2}\, ds \nn \\
&\lesssim& \| V_2\|_{L^2\cap L^\infty} \int^t_{t-A} \| 
\calg_{-\vec{e_1}}(s) 
e^{-isH_2}P_c(H_2) \calge1(0)\psi_0 \|_{L^2 + L^\infty} \, ds \nn \\
&\lesssim& \| V_2\|_{L^2\cap L^\infty}\int^t_{t-A} \frac{1}{\langle 
s\rangle^{n/2}} 
 \|\psi_0 \|_{L^1 \cap L^2} \, ds 
\lesssim \| V_2\|_{L^2 \cap L^\infty} A \langle t\rangle^{-n/2} 
\|\psi_0\|_{L^1\cap L^2}. \label{eq:zwisch}
\eea
To pass to \eqref{eq:zwisch} one uses that 
\[ \|P_c(H_2)f\|_{L^1\cap L^2} \le \|f\|_{L^1\cap L^2} + \|P_b(H_2) 
f\|_{L^1\cap L^2} \les \|f\|_{L^1\cap L^2}.\] 
Since $A$ is a fixed constant independent of the constant $C_0$ we
conclude that 
\begin{equation}
\label{eq:24}
\| J_c^1 \|_{L^2 + L^\infty} \lesssim \langle t\rangle^{-n/2} \| 
\psi_0\|_{L^1\cap L^2}.
\end{equation}
We decompose $J_c^2$ as follows:
\bea
J^2_c &=& \chi_1(t,\cdot) \int^t_{t-A} \int^{s-B}_0 e^{-i(t-s)H_1} 
P_c(H_1) V_2(\cdot-s\vec{e_1}) \calg_{-\vec{e_1}}(s) e^{-i(s-\tau) H_2} P_c(H_2) 
\nn \\
&& \mbox{\hspace{2in}}  V_1(\cdot + \tau\vec{e_1}) \calge1(\tau) 
U(\tau)\psi_0 \;d\tau ds \nn \\
&& +\chi_1(t,\cdot) \int^t_{t-A} \int^s_{s-B} e^{-i(t-s)H_1}P_c(H_1) 
F (|\vec p|\leq M) V_2(\cdot-s\vec{e_1}) \nn \\
&& \mbox{\hspace{2in}} \calg_{-\vec{e_1}}(s) e^{-i(s-\tau)H_2} 
P_c(H_2) V_1(\cdot+ \tau \vec{e_1}) \calge1(\tau) U(\tau) \psi_0 \; 
d\tau ds \nn \\
&& + \chi_1 (t,\cdot) \int^t_{t-A} \int^s_{s-B} e^{-i(t-s)H_1} 
P_c(H_1)
F(|\vec p|\geq M) V_2(\cdot-s\vec{e_1}) \nn \\
&& \mbox{\hspace{2in}} \calg_{-\vec{e_1}}(s) e^{-i(s-\tau) H_2} 
P_c(H_2) V_1(\cdot+\tau \vec{e_1}) \calge1(\tau)
U(\tau) \psi_0 \; d\tau d s  \nn \\
&& =: J_c^{2,\tau}+ J_c^{2,\text{ low }} + J_c^{2, \text{ high}}. \nn 
\eea
Here $B$ is a constant independent of $C_0$ so that $B>> A$, $M>> B$, 
and
 $F( |\vec p| \leq M), F(|\vec p|\geq M)$ denote smooth projections
 onto frequencies $|\vec p| \leq M, |\vec p| \geq M$, respectively.
The term $J_c^{2, \tau}$ is estimated similarly to $J_1^c$ above. 
Indeed, using that
$V_1(\cdot+\tau \vec{e_1})\calge1(\tau)=\calge1(\tau) V_1$, one has 
\bea
\| J_c^{2, \tau} \|_{L^2+ L^\infty} &\leq& \| J_c^{2,\tau}\|_{L^2} 
\lesssim
\int^t_{t-A} \int^{s-B}_0 \Big\| V_2(\cdot-s \vec{e_1}) 
\calg_{-\vec{e_1}}(s)e^{-i(s-\tau)H_2}P_c(H_2) \nn \\
&& \mbox{\hspace{2in}} V_1 (\cdot+\tau\vec{e_1}) \calge1(\tau) 
U(\tau) \psi_0 \Big\|_{L^2}\, d\tau ds \nn \\
&\lesssim& \int^t_{t-A} \int^{s-B}_0 \Big\| \calg_{-\vec{e_1}}(s) 
e^{-i(s-\tau)H_2} P_c(H_2)
\calge1(\tau) V_1 U(\tau) \psi_0 \Big\|_{L^2+ L^\infty}\,d\tau ds \nn 
\\
&\lesssim& \int^t_{t-A} \int^{s-B}_0\frac{1}{\langle 
s-\tau\rangle^{n/2}}
\|  V_1 U(\tau) \psi_0\|_{L^1\cap L^2} \, d\tau ds \nn \\
&\lesssim& \int^t_{t-A} \int^{s-B}_0\frac{1}{\langle 
s-\tau\rangle^{n/2}}\|U(\tau) \psi_0\|_{L^2+ L^\infty}\, d\tau d s 
\nn \\
&\lesssim& C_0 \int^t_{t-A} \int^{s-B}_B\frac{d\tau}{\langle
s-\tau\rangle^{n/2}}\frac{d s}{\langle \tau\rangle^{n/2}}\| \psi_0 
\|_{L^1 \cap L^2} 
+ \int_{t-A}^t \int_0^B \frac{d\tau}{\langle s-\tau\rangle^{n/2}} 
\|\psi_0\|_2\, d\tau\,ds
\nn \\
&\lesssim& C_0 A B^{-(n/2-1)}\langle t\rangle^{-n/2} \| \psi_0 
\|_{L^1 \cap L^2} + AB \la t\ra^{-\frac{n}{2}} \|\psi_0\|_{L^1\cap 
L^2}. \nn 
\eea
Since $B >> A$ and $n/2 -1 > 0$ we can conclude that 
\begin{equation}
\|J_c^{2,\tau}\|_{L^2 + L^\infty} \leq 10^{-2} C_0 \langle
t\rangle^{-n/2} \|\psi_0 \|_{L^1\cap L^2}, \label{eq:24'}
\end{equation}
provided $C_0$ is large.

\subsection{Low velocity estimates}

The idea behind the estimate for $J_c^{2, \text{ low}}$ is that the
norm of the operator 
\[ \chi_1(t,\cdot) e^{-i(t-s) H_1} P_c (H_1) F (|\vec p|\leq M) V_2 
(\cdot- s\vec{e_1})\]  
for all $s$ in the interval $[t-A, t]$ is small.
This can be explained as follows:

\noindent The support of $\chi_1 (t,x)$ with respect to $x$ belongs 
to the ball
$B_{2\delta t}(0)$.  On the other hand, since $V_2$ is compactly 
supported
$V_2(\cdot - s \vec{e_1})$ is different from $0$ on the set
$B_R(s\vec{e_1})$.  Here $R$ is the size of the support of $V_2$,
$R << t_0 \leq t$.  Since $s\in [t-A_1, t]$, we have that
$$
\bigcup_{s\in [t-A_1, t]}  B_R \subset B_{R +A} (\ t\vec{e_1})
$$
and thus the support of $V_2(\cdot-s\vec{e_1})$ is approximately
$(1-\delta)t$ units away from the support of $\chi_1(t,\cdot)$.
The operator $e^{-i(t-s)H_1} P_c (H_1) F( |\vec p| \leq M)$ can
``propagate'' the set $K$ into the set $L$ only if $(t-s) M
\ge \text{ dist} (K, L)$ according to the semiclassical picture.  
Setting
$K = \sppt\Big( V_2(\cdot - s\vec{e_1})\Big)$ and $L= \sppt\Big( 
\chi_1(t,\cdot)\Big)$, we see that this would require that
$$
(t-s)M \ge (1-\delta) t.
$$
However, $|t-s| \leq A$ and $M$ is fixed so that $A, M << t_0 \leq t$.

\noindent We make this argument rigorous in the following lemma.

\begin{lemma}
\label{lem:lowvel}
Let $A, M$ be large positive constants and $A, M <<  t$.
Then
$$
\sup_{|t-s|\leq A} \| \chi_1(t,\cdot) e^{-i(t-s)H_1} P_c(H_1) F(|\vec 
p|\leq M)
V_2 (\cdot-s\vec{e_1}) \|_{L^2\to L^2} \leq \frac{AM}{\delta t},
$$
where $\delta$ is the constant from the definition of $\chi_1, 
\chi_2$. 
\end{lemma}
\begin{proof}
The proof is a commutator argument. Let $\chi_1=\chi(t,x)$. 
Firstly, we claim that
\begin{equation}
\label{eq:1comm}
\|[\chi_1, P_c (H_1)]\|_{L^2\to L^2}\lesssim e^{-\delta\alpha t}.
\end{equation}
Clearly,
$$
[\chi_1, P_c (H_1)] = [\chi_1, I-P_b(H_1)] =-[\chi_1, P_b(H_1)].
$$ 
Recall that $u_1, \ldots , u_m$ are the exponentially decaying 
eigenvalues
of $H_1$.  Therefore,
\bea
[\chi_1,  P_b(H_1)]f &=& \sum^m_{i=1} (\chi_1 u_i\langle f,
u_i\rangle-u_i \langle f\chi_1, u_i\rangle) \nn \\
&=& 
\sum^m_{i=1} \left((-1 + \chi_1) u_i\langle f, u_i\rangle - u_i\langle
f,(\chi_1-1) u_i\rangle\right). \nn
\eea
In of the support of $\chi_1 - 1$ we have that 
$$
\|(1-\chi_1) u_i\|_2  \lesssim e^{-\alpha\delta t}
$$
and thus
$$
\|[\chi_1, P_b(H_1)]f\|_{L^2} \leq e^{-\alpha \delta t} \| f \|_{L^2},
$$
as desired. Secondly, we claim that
\begin{equation}
\label{eq:2comm}
\|[\chi_1, e^{-i(t-s)H_1}]F(|\vec p|\leq M) \|_{L^2\to L^2} \les
\frac{AM}{\delta t}.
\end{equation}
We write
$$[\chi_1, e^{-i(t-s)H_1}] = e^{-i(t-s)H_1} (e^{i(t-s)H_1}\chi_1
e^{-i(t-s)H_1}-\chi_1)
$$
and 
\bea 
e^{i(t-s) H_1} \chi_1\, e^{-i(t-s)H_1}-\chi_1 &=& i\int^{t-s}_0 
e^{i\tau
H_1} [H_1, \chi_1]e^{-i\tau H_1} \, d\tau \nn \\
&=& i\int^{t-s}_0 e^{i\tau H_1}(-\nabla \chi_1 \nabla - \Lapl \chi_1) 
e^{-i\tau H_1} \, d\tau. \nn
\eea
Observe now that 
$$
|\nabla \chi_1(t, x) | \lesssim \frac{1}{\delta t},\quad |\triangle 
\chi_1 (t,x) | \lesssim \frac{1}{(\delta t)^2}.
$$
Therefore,
\bea
&& \|[\chi_1, e^{-(t-s)H_1}] F(|\vec p|\leq M)\|_{L^2\to L^2} \nn\\
& &\lesssim |t-s|\Big(\|\nabla \chi_1 \, \nabla e^{-i\tau H_1} 
F(|\vec p|\leq M)\|_{L^2\to L^2}  +
\|\Delta \chi_1 e^{-i\tau H_1} F(|\vec p|\leq M) \|_{L^2\to L^2}\Big) 
\nn \\
&& \les \frac{A}{\delta t} \|\nabla e^{-i\tau H_1} F(|\vec p|\leq 
M)\|_{L^2\to
L^2} + \frac{A}{(\delta t)^2} \|e^{-i\tau H_1} F(|\vec p|\leq 
M)\|_{L^2\to L^2}. \nn
\eea
Since the potential $V_1$ is bounded it is standard that 
\begin{equation}
\label{eq:H1}
\sup_{\tau} \|\nabla e^{-i\tau H_1} f\|_{L^2}\lesssim \|\nabla 
f\|_{L^2} +
\|f\|_{L^2}.
\end{equation}
Indeed, 
\bea 
\sup_{\tau} \|\triangle e^{-i\tau H_1} f\|_{L^2} &\le& \sup_{\tau} \| 
e^{-i\tau H_1} H_1 f\|_{L^2} +
\sup_{\tau} \|V_1 e^{-i\tau H_1} f\|_{L^2} 
\les \|\nabla^2 f\|_{L^2} + \|f\|_{L^2}, \nn
\eea
and \eqref{eq:H1} follows by interpolation with $L^2$. 
Therefore,
$$
\|\nabla e^{-i\tau H_1} F(|\vec p|\leq M) \|_{L^2\to L^2} \leq M.
$$
Combining terms we obtain the bound
\begin{equation}
\|[\chi_1, e^{-i(t-s)H_1}] F(|\vec p|\leq M)\|_{L^2\to L^2} \lesssim
\frac{AM}{\delta t} + \frac{A}{(\delta t)^2} \leq 2\frac{AM}{\delta 
t} \nn
\end{equation}
since $M<<t$, which is \eqref{eq:2comm}. Finally, we invoke one more 
standard fact, namely
\begin{equation}
\|[\chi_1, F(|\vec p|\leq M)]\|_{L^2 \to L^2}\les 
M^{-1}\,\|\nabla\chi_1\|_\infty \les \frac{1}{\delta 
Mt}.\label{eq:3comm}
\end{equation}
To see this, write $F(|\vec p|\leq M) f = [\hat{\eta}(\xi/M) 
\hat{f}(\xi)]^{\vee}$ with some
smooth bump function~$\eta$. Hence the kernel $K$ of $[\chi_1, 
F(|\vec p|\leq M)]$ is
\[ K(x,y) = M^n\eta(M(x-y))(\chi_1(x)-\chi_1(y)),\]
and \eqref{eq:3comm} follows from Schur's test.
One concludes from estimates \eqref{eq:1comm}, \eqref{eq:2comm}, 
\eqref{eq:3comm} that
\bea
&& \Big\|\chi_1(t,\cdot) e^{-i(t-s)H_1} P_c(H_1) F(|\vec p|\leq M) 
V_2(\cdot - s\vec{e_1}) \nn \\
&&  \mbox{\hspace{1in}} 
-e^{-i(t-s)H_1} P_c(H_1)F(|\vec p|\leq M) \chi_1(t,\cdot) V_2(\cdot-s 
\vec{e_1})\Big\|_{L^\to L^2} 
 \lesssim \frac{AM}{\delta t}. \nn 
\eea 
It remains to observe that $\chi_1(t,\cdot) V_2 (\cdot-s\vec{e_1}) =0 
$ since the supports of $\chi_1(t,\cdot)$ and $V_2(\cdot-s\vec{e_1})$ 
are disjoint.
\end{proof}

\noindent We now estimate the term
\bea
J^{2, \text{ low}}_c &=& \chi_1\int^t_{t-A} \int^s_{s-B}
e^{-i(t-s)H_1} P_c(H_1) F(|\vec p|\leq M) V_2(\cdot-s\vec{e_1}) 
\calg_{-\vec{e_1}}(s)
e^{-i(s-\tau)H_2} P_c(H_2) \nn \\
&& \mbox{\hspace{1.5in}} V_1(\cdot+\tau \vec{e_1}) \calge1(\tau) 
U(\tau)\psi_0 \, d\tau ds \nn
\eea
by means of Lemmas~\ref{lem:lowvel} and~\ref{lem:dual}. Indeed, one 
has
\bea
&& \|J_c^{2, \text{ low }}\|_{L^2 + L^\infty} \leq  \|J_c^{2,\text{ 
low}}\|_{L^2} \nn \\
&&\les \frac{AM}{\delta t} \int^t_{t-A}\int^s_{s-B}
\Big \|\calg_{-\vec{e_1}}(s) e^{-i(s-\tau)H_2} P_c(H_2) \calge1(\tau) 
V_1  U(\tau)\psi_0\Big \|_{L^2} \; d\tau d s \nn \\
&&\les \frac{AM}{\delta t} \int^t_{t-A}\int^s_{s-B}\|V_1
U(\tau)\psi_0\|_{L^2} \, d\tau ds \lesssim \frac{AM}{\delta 
t}\|V_1\|_{L^2\cap L^\infty}
\int^t_{t-A}\int^s_{s-B} \|U(t)\psi_0\|_{L^2+L^\infty} \; d\tau ds 
\nn \\
&&\les \frac{A^2 BM}{\delta t} \|V_1\|_{L^2\cap L^\infty}\,C_0 
\langle t\rangle^{-\frac{n}{2}} \|\psi_0\|_{L^1\cap L^2}. \nn
\eea
Since $t>> A, M, \|V_1\|_{L^2\cap L^\infty}$ we obtain that 
\begin{equation}
\|J_c^{2, \text{ low }}\|_{L^2+L^\infty} \leq 10^{-2} C_0\langle
t\rangle^{-n/2} \|\psi_0\|_{L^1\cap L^2}.\label{eq:28}
\end{equation}

\subsection{Large velocity estimates via Kato smoothing}

The estimate for $J_c^{2, \text{ high }}$ will require the use of 
the Kato $\half$-smoothing estimate. We present a variant of it
in the following lemma. 
It differs from the original version (see \cite{CS1}, \cite{KY}, \cite{Sj}, 
\cite{Ve}) in so far as we consider 
$-\Lapl+V$ instead of~$-\Lapl$.
Moreover, we employ commutator methods. Both the statement and the 
proof are known, see e.g. \cite{CS2}, \cite{BK}, 
\cite{So}, \cite{Do}, so we will only sketch the argument. 
\begin{lemma}
\label{lem:kato1/2}
Let $H=-\Lapl+V$ , $\|V\|_\infty<\infty$. Then for all $T,R,M\ge 1$, 
\begin{equation}
\label{eq:kato1/2}
\sup_{B_R}\int^T_0 \|F(|\vec p|\geq M) e^{-it H} f\|_{L^2(B_R)}\,dt
\leq C(n,V)\,\frac{T R}{M^{\frac 12}}\|f\|_{L^2}.
\end{equation}
Here the supremum ranges over all balls $B_R$ of radius $R\ge 1$ and 
$C(n,V)$ is a
constant that depends only on $\|V\|_\infty$ and the dimension~$n$. 
\end{lemma}
\begin{proof} We first prove the following estimate, which will then 
imply the lemma.
Let $\psi(t)=e^{-itH}\psi_0$ with $H=-\Lapl+V$, 
$\|V\|_\infty<\infty$. Then for all $T>0$ and $0<\alpha$,  
\begin{equation}
\label{eq:kato_first}
\sup_{x_0\in\R^n}\int_0^T \int_{\R^n} 
\frac{|\nabla{\langle \nabla\rangle^{-\frac 12}}\psi(x,t)|^2}
{(1+|x-x_0|^\alpha)^{\frac{1}{\alpha}+1}} \, 
dxdt \le 
C_{\alpha,n}\,T(1+\|V\|_\infty)\|\psi_0\|_2^{2}
\end{equation}
The multiplier $ \nabla{\langle \nabla\rangle^{-\frac 12}}$ 
corresponds to the symbol $\xi{\langle \xi\rangle^{-\frac 12}}=
\xi(1+|\xi|^2)^{-\frac14}$. 
It suffices to prove this with $x_0=0$ fixed.
The proof is based on taking the commutator of $H$ with 
$m:=w(x)x\cdot\frac{\nabla}{\langle \nabla\rangle }$, where
\[ w(x) = (1+|x|^\alpha)^{-\frac{1}{\alpha}},\qquad \alpha>0.\]
One has, with $\psi=\psi(t)$ for simplicity, 
\bea
\frac{d}{dt} \la m\psi,\psi \ra &=& -i  \la [m,H] \psi, \psi \ra \nn 
\\
\int_0^T \la -[m,H]\psi(t), \psi(t) \ra \,dt &=& i \la 
m\psi(0),\psi(0) \ra-i\la m\psi(T),\psi(T) \ra \nn\\
\la -[m,H]\psi, \psi \ra  &=& -\la [m,V]\psi,\psi 
\ra 
+ \Big\la \partial_\ell(w(x)\,x_j)\frac {\partial_j}{\la \nabla\ra} \psi,
\partial_\ell \psi \Big\ra 
+ \half \Big\la \triangle (w(x)\,x_j)\frac{\partial_j}{\la \nabla\ra } \psi,
\psi \Big\ra \label{eq:mHcomm} \\
&=& \int_{\R^n} w |{\nabla}{\la\nabla\ra^{-\frac 12}}\psi|^2\,dx + 
\int_{\R^n} (\partial_\ell 
w)(x) x_j{\partial_j}{\la\nabla\ra^{-\frac 12}}\psi\,
{\partial_\ell}{\la\nabla\ra^{-\frac 12}} \bar\psi\,dx 
\label{eq:positive}\\ &-&
\Big\la [\la\nabla\ra^{\frac 12},w]\frac{\nabla}{\la\nabla\ra }\psi, 
{\nabla}{\la \nabla\ra^{-\frac 12}}\psi \Big\ra\, - 
\Big\la [\la\nabla\ra^{\frac 12},(\partial_\ell 
w)(x) x_j]\frac {\partial_j}{\la\nabla\ra }\psi ,
{\partial_\ell}{\la\nabla\ra^{-\frac 12}}\psi \Big\ra\label{eq:commutat} \\ &+&
\half \Big\la \triangle (w(x)\,x_j)\frac{\partial_j}{\la \nabla\ra } \psi,\psi \Big\ra 
-\la [m,V]\psi,\psi \ra  \label{eq:remnan}
\eea
One now checks easily that the two terms in \eqref{eq:positive} 
satisfy
\bea
\int_{\R^n} w(x) |\la\nabla\ra^{-\frac 12}\nabla\psi(t,x)|^2\,dx &+& \int_{\R^n} 
(\partial_\ell w)(x) x_j\la\nabla\ra^{-\frac 
12}\partial_j\psi(t,x)\,\la\nabla\ra^{-\frac12}\partial_\ell 
\bar\psi(t,x)\,dx \nn \\ &\ge& 
\int_{\R^n} \frac{w(x)}{1+|x|^\alpha}\, |\la\nabla\ra^{-\frac 12}\nabla\psi(t,x)|^2\,dx 
= \int_{\R^n} 
\frac{|\la\nabla\ra^{-\frac 12}\nabla\psi(t,x)|^2}{(1+|x|^\alpha)^{\frac{1}{\alpha}+1}} \, dx. 
\label{eq:des}
\eea
Notice that \eqref{eq:des} is precisely the space integral in the 
desired lower bound from~\eqref{eq:kato_first}.

There are several ways to bound the commutators $A_1:=[\la\nabla\ra^{\frac 12},w]$ and 
$ A_2:=[\la\nabla\ra^{\frac 12},(\partial_\ell w)(x) x_j]$.
For example, one can use the Kato square root formula as in~\cite{JSS}. 
Alternatively, one can invoke the standard composition formula from 
$\Psi$DO~calculus, see~\cite{Taylor}. This gives
$ A_1 = T_{\{\la\nabla\ra^{\frac 12},w\}} + R_{a_1}$
where $\{\la\nabla\ra^{\frac 12},w\}$ is the Poisson bracket of the symbols 
$\la\nabla\ra^{\frac 12}$ and~$w$. Moreover,  $T_{\{\la\nabla\ra^{\frac 12},w\}}$ is an
associated $\Psi$DO operator and $R_{a_1}$ another $\Psi$DO operator with 
symbol $a_1\in S^{-\frac32}$. A similar expression holds
for~$A_2$. One checks that $|\{\la\nabla\ra^{\frac 12},w\}(\xi,x)|\les |\nabla w(x)|$
for all $x,\xi$. 
Therefore, the two terms in \eqref{eq:commutat} both satisfy
\bea
\Big| \big\la A_i \frac{\nabla}{\la\nabla\ra}\psi(t), \nabla \la\nabla\ra^{-\frac12}\psi(t) 
 \Big\ra \Big|
&\les&
\|\psi(t)\|_2 \Big\|\frac{\la\nabla\ra^{-\frac 12}\nabla\psi(t)}
{(1+|x|^{\alpha})^{\frac 1\alpha +1}}\Big\|_2 + \|\psi(t)\|_2 
\|R_{a_i}  \nabla \la\nabla\ra^{-\frac12} \psi(t) \|_2 \nn \\
&\le& C \|\psi(t)\|_2^2 + \frac14 \Big\|\frac{\la\nabla\ra^{-\frac 12}\nabla\psi(t)}
{(1+|x|^{\alpha})^{\frac 1\alpha +1}}\Big\|_2^2. \label{eq:commut}
\eea
Finally, the two terms in \eqref{eq:remnan} are bounded by
\begin{equation}
\label{eq:einf}
(1+ \|V\|_{\infty})\|\psi(t)\|_{{2}}^{2} \le (1+ \|V\|_{\infty})\|\psi_{0}\|_{{2}}^{2}. 
\end{equation}
In the above estimate we have used the boundedness of the multipliers 
$m$ and $\triangle (w(x)\,x_j)\frac{\partial_j}{\la \nabla\ra }$  on $L^{2}$.
Integrating \eqref{eq:des}, \eqref{eq:commut}, and~\eqref{eq:einf} in time,  
inserting the resulting bounds into~\eqref{eq:mHcomm}, and finally 
using that
\[ |\la m\psi(0),\psi(0) \ra| + |\la m\psi(T),\psi(T) \ra| \le 
2\|\psi_0\|^{2}_2,\]
one obtains \eqref{eq:kato_first}. To pass from~\eqref{eq:kato_first} 
to~\eqref{eq:kato1/2},
let $\chi_R$ be a smooth cut-off to the ball $B_R$, so that 
$\widehat{\chi_R}$ has compact support in a ball
of size~$\sim R^{-1}$.  Then, by~\eqref{eq:kato_first} 
with~$\alpha=1$ 
\bea 
&& \int_0^T \|F(|\vec p|\ge M)e^{-itH}\psi_0\|_{L^2(B_R)}^2 \,dt \nn 
\le 
\int_0^T \|\chi_R F(|\vec p|\ge M) e^{-itH}\psi_0\|_{L^2}^2 \,dt \nn \\
&\les& \int_0^T \|F(|\vec p|\ge M) \chi_{R} e^{-itH} 
\psi_0\|_{L^2}^2 \,dt +  
\int_0^T \|[\chi_R, F(|\vec p|\ge M)]\|^{2}_{2\to 2}
\|e^{-itH} \psi_0\|_{L^2}^2 \,dt \nn \\
&\les& M^{-1} \int_0^T 
\|\nabla \la\nabla\ra^{-\frac 12}F(|\vec p|\ge M) \chi_{R} e^{-itH} 
\psi_0\|_{L^2}^2 \,dt + T (MR)^{-2} \|\psi_{0}\|_{L^{2}}^{2}\nn \\&\les&
 M^{-1} \int_0^T 
\|F(|\vec p|\ge M) \chi_{R}
\nabla \la\nabla\ra^{-\frac 12} e^{-itH} 
\psi_0\|_{L^2}^2 \,dt +  M^{-1} \int_0^T \| [\nabla 
\la\nabla\ra^{-\frac 12}, \chi_{R}]\|_{2\to 2}^{2} 
\|e^{-itH} \psi_0\|_{L^2}^2 \nn\\ &+&
T (MR)^{-2} \|\psi_{0}\|_{2}^{2} \les 
M^{-1} R^{2} \int_0^T \int_{\R^{n}}\frac {|\nabla \la\nabla\ra^{-\frac 12} 
 e^{-itH}\psi_{0}|^{2}}{(1+|x|)^{2}} dx\,dt + T M^{-1} R^{-1} 
\|\psi_{0}\|_{L^{2}}^{2}\nn\\ &\le & C(n,V)\, T M^{-1} R^{2}\, 
\|\psi_{0}\|_{L^{2}}^{2}.\nn
\eea
The lemma follows.
\end{proof}

\noindent We then have for $J_c^{2, \text{ high }}$
\bea
&& \qquad \|J_c^{2, \text{ high }} \|_{L^2+L^\infty} \leq \|J_c^{2, 
\text{ high}}\|_{L^2} \label{eq:rep1} \\
&& \les \int^t_{t-A}\int^s_{s-B} \|F(|\vec p|\geq M) V_2 
(\cdot-s\vec{e_1}) \calg_{-\vec{e_1}}(s)
e^{-i(s-\tau) H_2} P_c(H_2)  V_1(\cdot+\tau \vec{e_1}) \calge1(\tau)
 U(\tau) \psi_0\|_{L^2} \, d\tau ds \nn \\
&&\les \int^t_{t-A}\int^s_{s-B} \|F(|\vec p|\geq M) 
\calg_{-\vec{e_1}}(s) V_2 
e^{-i(s-\tau) H_2} P_c(H_2) \calge1(\tau) V_1
 U(\tau) \psi_0\|_{L^2} \, d\tau ds. \nn 
\eea
We need to commute the cutoff $F(|\vec p|\geq M)$ past the first two 
terms following it. 
As for the Galilei transform, one has
\[
\| F(|\vec p|\geq M) \calg_{-\vec{e_1}}(s) f \|_2 = \| F(|\vec 
p|\geq M) e^{i{x_1}} e^{-ip_1s} f\|_2 = \| F(|\vec p|\geq M) e^{i{x_1}} f \|_2 = 
\| F(|\vec p -  \vec{e_1}|\geq M) f\|_2. 
\]
Therefore, applying this to the final line of~\eqref{eq:rep1}, one 
obtains
\bea
&& \|J_c^{2, \text{ high }} \|_{L^2+L^\infty} \les \int^t_{t-A}
\int^s_{s-B} \|V_2 F(|\vec p -  \vec{e_1}|\geq M) 
e^{-i(s-\tau) H_2} P_c(H_2)  
 \calge1(\tau) V_1  U(\tau) \psi_0\|_{L^2} \, d\tau ds \nn \\
&& \qquad + \int^t_{t-A}\int^s_{s-B} \big \|[V_2, F(|\vec p -  
\vec{e_1}|\geq M)]\; e^{-i(s-\tau) H_2} P_c(H_2)  \calge1(\tau) V_1  
U(\tau) \psi_0\big \|_{L^2} \, d\tau ds.
\label{eq:rep2}
\eea
The final two terms in \eqref{eq:rep2} are now easily estimated by 
means of Lemma~\ref{lem:kato1/2} (integrating in the 
variable~$s-\tau$) and~\eqref{eq:3comm}, respectively. To apply 
Lemma~\ref{lem:kato1/2} choose a ball $B_R$ around the origin that 
contains $\sppt(V_2)$. The conclusion is that 
\bea
\|J_c^{2, \text{ high }} \|_{L^2+L^\infty} &\les& 
\|V_2\|_{L^\infty} \frac{BR}{M^{\frac 12}}\int^t_{t-B-A}
\|P_c(H_2) \calge1(\tau)V_1 U(\tau)
\psi_0\|_{L^2} \, d\tau \nn \\
&& \qquad +  M^{-1} \int^t_{t-A}\int^s_{s-B} \big \| e^{-i(s-\tau) 
H_2} P_c(H_2) \calge1(\tau) V_1  U(\tau) \psi_0\big \|_{L^2} \, d\tau 
ds \nn \\ 
&\les& \Bigl(\|V_2\|_{L^\infty} \frac{BR}{M^{\frac 12}} B + M^{-1} B^2 
\Bigr) \|V_1\|_{L^1\cap L^\infty} C_0 \langle t\rangle^{-n/2} 
\|\psi_0 \|_{L^1\cap L^2}. \nn 
\eea
Thus for a sufficiently large constant $M$ one obtains
\begin{equation}
\|J_c^{2, \text{ high }}\|_{L^2+L^\infty} \leq 10^{-2} C_0 \langle
t\rangle^{-n/2}\|\psi_0\|_{L^1\cap L^2}. \label{eq:29}
\end{equation} 
Combining \eqref{eq:21}, \eqref{eq:23}, \eqref{eq:24}, 
\eqref{eq:24'}, \eqref{eq:28}, \eqref{eq:29}
yields
\begin{equation}
\label{eq:piece1}
\|\chi_1(t,\cdot) P_c U(t) \psi_0\|_{L^2+L^\infty} \leq 
\frac{C_0}{10} \langle t\rangle^{-n/2} \|\psi_0\|_{L^1\cap L^2}.
\end{equation}
This is precisely the desired bound on the first channel. Since the 
second channel can be treated
by the very same method, see the discussion 
following~\eqref{eq:DU2mov}, only the third channel remains.

\subsection{The third channel}

By \eqref{eq:DU3}, 
$$
\chi_3\, U(t)\psi_0 = \chi_3\, e^{it\Laplace} \phi_0- i\chi_3\,
\int^t_0 e^{i(t-s)\Laplace} (V_1 + V_2 (\cdot-s\vec{e_1})) U(s) 
\psi_0 \,ds
$$
where $\chi_3=\chi_3(t,x)$. 
By the standard decay estimates for $e^{it\Laplace}$,
$$ 
\|e^{it\Laplace} \psi_0 \|_{L^2 + L^\infty} \les \langle
t\rangle^{-n/2} \| \psi_0\|_{L^1\cap L^2}.
$$
Thus
\[
\|\chi_3\, U(t)\psi_0\|_{L^2 +L^\infty} \les \langle t\rangle^{-n/2}
\|\psi_0\|_{L^1\cap L^2}  + \Bigl\|\int^t_0\chi_3\, e^{i(t-s) \Laplace}
(V_1 + V_2(\cdot-s \vec{e_1}))U(s) \psi_0\, ds\Bigr\|_{L^2+L^\infty}.
\]
By Lemma~\ref{lem:dual}  and the bootstrap assumption~\eqref{eq:boot},
\bea
&& \|(V_1 + V_2(\cdot - s\vec e_1))U(s) \psi_0 \|_{L^1\cap L^2} \leq
(\|V_1\|_{L^1\cap L^\infty} + \| V_2\|_{L^1\cap 
L^\infty})\|U(s)\psi_0\|_{L^2+L^\infty} \nn \\
&&  \les C_0\,(\|V_1\|_{L^1\cap L^\infty} +
\|V_2\|_{L^1\cap L^\infty}) \langle t\rangle^{-n/2}
\|\psi_0\|_{L^1\cap L^2}. \nn
\eea
Therefore, by the same calculation that lead to~\eqref{eq:At-A}, 
\bea
\|\chi_3\, U(t) \psi_0\|_{L^2+L^\infty} &\les& \langle t\rangle^{-n/2}
\|\psi_0\|_{L^1\cap L^2} + C_0\,A^{-n/2 + 1}\langle t\rangle^{-n/2}
\|\psi_0\|_{L^1\cap L^2} \nn \\
&& \qquad + \Bigl\|\int^t_{t-A} \chi_3\, e^{i(t-s)\Laplace} (V_1 + 
V_2(\cdot-s\vec{e_1}))U(s)\psi_0\, ds \Bigr\|_{L^2+L^\infty}. \nn 
\eea
Since $A$ is sufficiently large we only need to estimate the last
term, which we can split as follows:
$$
\int^t_{t-A} \chi_3\, e^{i(t-s)\Laplace} V_1 U(s) \psi_0\, ds +
\int^t_{t-A} \chi_3\, e^{i(t-s) \Laplace} V_2(\cdot - s\vec{e_1})U(s) 
\psi_0\, ds.
$$
It suffices to consider the first term, which is split further by 
means of the projections $P_b(H_1)$
and $P_c(H_1)=Id-P_b(H_1)$:
\bea
 \int^t_{t-A} \chi_3\, e^{i(t-s)\Laplace} V_1 U(s) \psi_0 \, ds &=&
 \int^t_{t-A} \chi_3\, e^{i(t-s)\Laplace} V_1  P_b(H_1)U(s) \psi_0 \, 
ds \nn \\
&& + \int^t_{t-A} \chi_3\, e^{i(t-s)\Laplace} V_1 P_c(H_1)U(s) \psi_0 
\, ds =: \calL_1 + \calL_2. 
\eea
Since we are assuming asymptotic orthogonality of $U(t)\psi_0$ to the 
bound states, see
Definition~\ref{def:asymp}, Proposition~\ref{prop:bdstates} implies 
that
\bea
\|\calL_1\|_{L^2+L^\infty} &\leq& \|\calL_1\|_{L^2}\leq \int^t_{t-A}
\|V_1\|_{L^2\cap L^\infty} \|P_b(H_1) U(s) \psi_0\|_{L^2} \, ds \nn \\
&\les&  A\|V_1\|_{L^1\cap L^\infty} e^{-\alpha t/2}\|\psi_0\|_{L^2} 
\les
\langle t\rangle^{-n/2} \|\psi_0\|_{L^1\cap L^2}.\label{eq:31}
\eea
It is now natural to expand $\calL_2$ further, cf.~\eqref{eq:Jc12def}:
\bea
\calL_2 &=& \chi_3\, \int^t_{t-A} e^{i(t-s)\Laplace} V_1 P_c(H_1) 
e^{-isH_1} \psi_0 \, ds \label{eq:int1}\\
&& -i\chi_3\, \int^t_{t-A} \int^s_0 e^{i(t-s)\Laplace} V_1 P_c(H_1)
e^{-i(s-\tau) H_1} V_2(\cdot - \tau \vec{e_1})
U(\tau) \psi_0 \, d\tau ds. \label{eq:int2}
\eea 
The first integral~\eqref{eq:int1} is controlled by means of the 
decay estimate for the 
evolution $e^{-itH_1} P_c(H_1)$: 
\bea 
&& \Bigl \|\chi_3\,\int^t_{t-A} e^{i(t-s)\Laplace} V_1 P_c(H_1) e^{-i 
sH_1} \psi_0 \, 
ds \Bigr \|_{L^2+L^\infty} 
\les \|V_1\|_{L^2\cap L^\infty} \int^t_{t-A} \|e^{-isH_1} P_c (H_1)
\psi_0 \|_{L^2+L^\infty}\, ds \nn \\
&& \les \|V_1\|_{L^2\cap L^\infty} A \, \langle t\rangle^{-n/2}
\|\psi_0\|_{L^1\cap L^2} \les\langle t 
\rangle^{-n/2}\|\psi_0\|_{L^1\cap L^2}\label{eq:32}.
\eea
For the second integral~\eqref{eq:int2} one also invokes the bootstrap 
assumption~\eqref{eq:boot}, as well as Lemma~\ref{lem:dual}:
\bea
&& \Bigl \|\chi_3\,\int^t_{t-A} \int^{s-B}_0 e^{i(t-s)\Laplace} V_1 
e^{-i(s-\tau)H_1}P_c(H_1)  V_2(\cdot - \tau \vec{e_1}) U(\tau) \psi_0 
\,d\tau ds\Bigr \|_{L^2+L^\infty} \nn \\
&& \les\|V_1\|_{L^2\cap L^\infty}
\int^t_{t-A}\int^{s-B}_0 \Bigl \|e^{-i(s-\tau)H_1} P_c(H_1) 
V_2(\cdot-\tau\vec{e_1}) U (\tau)\psi_0 \,d\tau ds\Bigr\|_{L^2+ 
L^\infty} \nn \\
&& \les\|V_1\|_{L^2\cap L^\infty} \|V_2\|_{L^1\cap L^\infty} 
\int^t_{t-A}\int^{s-B}_0\frac{1}{\langle
s-\tau\rangle^{n/2}}\|U(\tau)\psi_0\|_{L^2+L^\infty} \,d\tau ds \nn \\
&& \les\|V_1\|_{L^2\cap L^\infty} \|V_2\|_{L^1\cap L^\infty} \left\{ 
C_0\,\int^t_{t-A} \int^{s-B}_B \frac{d\tau}{\langle s-\tau 
\rangle^{n/2}}
\frac{ds}{\langle \tau\rangle^{n/2}}\|\psi_0\|_{L^1\cap L^2} \,d\tau 
ds \right. \nn \\
&& \left. \mbox{\hspace{2in}} + \int^t_{t-A} \int^{B}_0 
\frac{d\tau}{\langle s-\tau \rangle^{n/2}}
    \|\psi_0\|_{ L^2} \,d\tau ds \right\}\nn \\
&& \les \|V_1\|_{L^2\cap L^\infty} \|V_2\|_{L^1\cap L^\infty} \Bigl( 
AB^{-n/2+1} C_0\, \|\psi_0\|_{L^1\cap L^2} +  AB \la 
t\ra^{-\frac{n}{2}} \|\psi_0\|_{ L^2} \Bigr) \nn \\
&& 
\leq 10^{-2} C_0\|\psi_0\|_{L^1\cap L^2} \nn
\eea
since  $B>>A$, and provided $C_0$ is large.
As in the proof of the decay estimates for the first channel, we 
split the remaining expression, 
which is the most difficult part, into low and high momenta:
$$
\chi_3\, \int^t_{t-A} \int^s_{s-B} e^{-i(t-s)\Laplace} V_1
e^{-i(s-\tau)H_1}P_c(H_1) V_2(\cdot - \tau\vec{e_1}) U(\tau) \psi_0\, 
d\tau ds = 
\calL_1^{\text{ low}}+\calL_1^{\text{ high}}
$$
where we define
\bea 
\calL_1^{\text{ low}} &=& \chi_3\, \int^t_{t-A}\int^s_{s-B}
e^{-i(t-s)\Laplace} F(|\vec p|< M) V_1 e^{-i(s-\tau) H_1} P_c(H_1)
V_2(\cdot-\tau\vec{e_1}) U(\tau) \psi_0\, d\tau ds \nn \\
\calL_1^{\text{ high}} &=& \chi_3\, \int^t_{t-A} \int^s_{s-B}
e^{-i(t-s)\Laplace} F(|\vec p|\geq M) V_1
e^{-i(s-\tau)H_1} P_c (H_1) V_2 (\cdot-\tau\vec{e_1}) U(\tau)\psi_0 
\, d\tau ds. \nn
\eea
As in Lemma~\ref{lem:lowvel} one has
\begin{equation}
\label{eq:dist}
\sup_{|t-s|\leq A} \|\chi_3\, e^{-i(t-s)\Laplace} F(|\vec p|\leq M)
V_1 \|_{L^2 \to L^2} \les \frac{AM}{\delta t}.
\end{equation}
This is just a consequence of the following:
\begin{enumerate}
\item $\sppt(\chi_3(t, \cdot)) \cap \sppt (V_1) = \emptyset \text{\ \ 
for\ \ } t\geq t_0
= \left(\frac{C_0}{2}\right)^{2/n}$
\item $ |\vec \nabla \chi_3 (t,\cdot) | \les \frac{1}{\delta t}.$
\end{enumerate}
The first property holds since
$$
\chi_3(t,x) = 1 - \chi_1(t, x)- \chi_2 (t, x) = 0 
$$
on the set
$$
\{ x: |x| \leq \delta t\} \cup \{ x: |x-r\vec{e_1}| \leq \delta t\}.
$$
The bound \eqref{eq:dist} leads to the following estimate on 
$\calL_1^{\text{ low}}$
\bea
\|\calL_1^{\text{ low}} \|_{L^2+ L^\infty} &\leq&  
\|\calL_1^{\text{low}} \|_{L^2} \nn \\
&\les&  \frac{AM}{\delta t} \|V_1\|_{L^\infty}\|V_2\|_{L^2\cap 
L^\infty} \int^t_{t-A} \int^s_{s-B} \|U(\tau) \psi_0\|_{L^2+L^\infty} 
\, d\tau ds \nn \\ 
&\les& \frac{A^2 MB}{\delta t} \|V_1\|_{L^\infty} \|V_2\|_{L^2\cap 
L^\infty}
C_0 \langle t\rangle^{-n/2} \|\psi_0\|_{L^1\cap L^2}. \nn
\eea
Since $A,M<B$ are large but universal constants and $t\geq t_0 = 
\left(\frac{C_0}{2}\right)^{2/n}$,
this implies that
\begin{equation}
\|\calL_1^{\text{ low}} \|_{L^2+L^\infty} \leq 10^{-2} C_0 \langle
t\rangle^{-n/2} \|\psi_0\|_{L^1\cap L^2} \|U(\tau)\psi_0\|_{L^2 + 
L^\infty}.\label{eq:33}
\end{equation}
Finally, the high energy part is treated by means of 
Lemma~\ref{lem:kato1/2} as before.
More precisely, one has
\bea 
\calL_1^{\text{high}} &=& \chi_3\, \int^t_{t-A} \int^s_{s-B} 
e^{-i(t-s)\Laplace} [F(|p|\geq M), V_1 ]
e^{-i(\tau- s) H_1} P_c (H_1) V_2(\cdot-\tau \vec{e_1}) 
U(\tau)\psi_0\, d\tau ds \nn \\
&& \qquad +\chi_3\, \int^t_{t-A} \int^s_{s-B} e^{-i(t-s) \Laplace} 
V_1 F(|p|\geq M)
 e^{-i(s-\tau) H_1} P_c(H_1) V_2(\cdot -\tau \vec{e_1}) U(\tau)\psi_0 
\, d\tau ds \nn \\
&=& \calL_{2,1}^{\text{ high}} + \calL_{2,2}^{\text{ high}}. \nn 
\eea
The commutator estimate, see \eqref{eq:3comm}, 
$$\|[F(|p|\geq M), V_1 ] \|_{L^2\to L^2} \les \frac{1}{M}\|\vec 
\nabla V_1\|_{L^\infty}, 
$$
implies that
\bea
\|\calL_{2,1}^{\text{ high}}  \|_{L^2+L^\infty} &\les& 
\frac{AB}{M} C_0\, \langle t\rangle^{-n/2} \|\psi_0\|_{L^1\cap L^2}
\les 10^{-2} C_0\, t^{-n/2} \|\psi_0\|_{L^1\cap L^2}.\label{eq:34}
\eea
By the smoothing estimate of Lemma~\ref{lem:kato1/2} applied to the 
integration variable $u=s-\tau$ and
with $R$ sufficiently large depending on the support of~$V_1$,
\bea 
&& \|\calL_{2,2}^{\text{ high}} \|_{L^2+L^\infty}
 \les \int^t_{t-A-B}
\int^B_0 \|V_1 F(|p|\geq M) e^{-iuH_1}
P_c(H_1) V_2(\cdot - \tau\vec{e_1}) U(\tau) \psi_0 \|_{L^2}\, du 
d\tau \nn \\
&&\les \frac{\|V_1\|_\infty}{M^{\frac12}}\, BR\, \|V_2\|_{L^1\cap 
L^\infty}
\int^t_{t-A-B} \|U(\tau) \psi_0\|_{L^2+L^\infty} \,d\tau 
\les \|V_1\|_{L^\infty} \|V_2\|_{L^1\cap L^\infty}\,
\frac{BR}{M^{\frac12}}\,  C_0\, \langle t\rangle^{-n/2} 
\|\psi_0\|_{L^1\cap L^2}. \nn 
\eea
Therefore, with $M$ large, 
\begin{equation}
\|\calL_{2,2}^{\text{ high}}\|_{L^2+L^\infty} \leq 10^{-2}  C_0 
\|\psi_0\|_{L^1\cap L^2}.
\label{eq:35}
\end{equation}
Combining \eqref{eq:31} - \eqref{eq:35} we obtain
\[
\|\chi_3\, U(t) \psi_0\|_{L^2+L^\infty} \leq 10^{-2} C_0\,\langle
t\rangle^{-n/2} \|\psi_0\|_{L^1\cap L^2},
\]
which concludes the proof of \eqref{eq:boot/2}, and thus of 
Theorem~\ref{thm:main}.

\subsection{Perturbations of charge transfer models}
\label{subsec:perturb}

It is easy to see that the argument proving
Theorem~\ref{thm:main} also yields the following
result.

\begin{theorem}
\label{thm:main_perturb}
Consider the equation
\bea
&& \frac{1}{i} \partial_t \psi - \Lapl \psi + \sum^m_{\kappa =1}
V_\kappa(x - \vec v_\kappa t) \psi + V_0(t,x)\psi = 0 \label{eq:transfer_perturb} \\
&& \psi |_{t=0} = \psi_0,  x\in \R^n, \nn
\eea
where the charge transfer part is an in 
Definition~\ref{def:chargetrans} and
the perturbation satisfies
\[ \sup_{t} \|V_0(t,\cdot)\|_{1\cap\infty} < \eps.\]
Let $\tilde{U}(t)$ denote the propagator of the 
equation~\eqref{eq:transfer_perturb}. Then 
for any initial data $\psi_0 \in L^1\cap L^2$, which is 
asymptotically orthogonal to the bound
states of $H_j$ in the sense of Definition~\ref{def:asymp}
(with $\tilde U(t)$ instead of~$U(t)$), 
one has the decay estimates
\begin{equation}
\| \tilde U(t) \psi_0 \|_{L^2 + L^\infty} \lesssim
\langle t\rangle^{-n/2}\|\psi_0\|_{L^1\cap L^2} \nn
\end{equation}
provided $0<\eps<\eps_0$ is sufficiently small with $\eps_0$ 
independent of~$\psi_0$.  
\end{theorem}

The proof of this theorem is basically identical with 
the proof of Theorem~\ref{thm:main}. For example, consider
again the case of two potentials $m=2$. Then the Duhamel
formula relative to the first channel, see~\eqref{eq:DU1},
becomes
\bea
\chi_1 P_c (H_1) \tilde U(t) \psi_0 &=& \chi_1 e^{-it H_1} P_c(H_1)
\psi_0 -i\chi_1\int^t_0 e^{-i(t-s) H_1} P_c (H_1) 
V_2(\cdot-s \vec{e_1})\tilde U(s) \psi_0 \, ds \nn \\
&& \qquad -i\chi_1\int^t_0 e^{-i(t-s) H_1} P_c (H_1) V_0(s,\cdot) 
\tilde U(s) \psi_0 \, ds. \label{eq:stor}
\eea
Impose the bootstrap assumption~\eqref{eq:boot}
on the evolution~$\tilde U$. Then the new 
term~\eqref{eq:stor} satisfies
\bea
\|\eqref{eq:stor}\|_{2+\infty} &\le& \int_0^t C\,\la t-s\ra^{-\frac32}
\sup_{0\le s\le t} \|V_0(s,\cdot)\|_{1\cap\infty} C_0\, \la s\ra^{-\frac32}
\|\psi_0\|_{1\cap2} \, ds \nn \\
&\le& C\,\eps\,C_0 \la t\ra^{-\frac32} \|\psi_0\|_{1\cap 2}. \nn
\eea
Modifying the other two channels \eqref{eq:DU2stat} and~\eqref{eq:DU3}
in the same way, one concludes that the only change
from the  estimates in the previous section is the addition
of a term~$C\,\eps\,C_0 \la t\ra^{-\frac32} \|\psi_0\|_{1\cap 2}$. 
Choosing~$\eps$ small, one again arives at the 
improved bootstrap assumption~\eqref{eq:boot/2}
for the evolution~$\tilde U(t)$. 

\section{Asymptotic completeness}
\label{sec:asymp}

Our next goal is to provide the following version of the asymptotic
completeness for the charge transfer model (see \cite{Ya2}, \cite{Gr}).
In this section we require more regularity of the potential,
see Lemma~\ref{le:woper} below. 

\begin{theorem}
\label{thm:Assymp}
Let, as before,  $u_1,\ldots, u_m$ and $w_{1},\ldots, w_\ell$ be the eigenfunctions 
of $H_1 =- \Laplace + V_{1}(x)$ and $H_2=-\Laplace + V_{2}(x)$, respectively, 
corresponding to the negative eigenvalues $\lambda_1,
\cdots, \lambda_m$ and $\mu_1, \cdots, \mu_\ell$.
Then for any initial data $\psi_0 \in  L^2$ the solution 
$U(t)\psi_{0}$ of the 
charge transfer problem \eqref{eq:2pot} can be written in the form
\begin{equation}
\nn
U(t) \psi_0 = \sum^m_{r=1} A_r e^{-i\lambda_rt} u_r + 
\sum^\ell_{s=1}B_s e^{-i\mu_s t} \calg_{-\vec e_{1}} (t) w_s + e^{-it\Laplace} \phi_{0}
+ \mathcal {R}(t),
\end{equation}
for some choice of the constants $A_{r}, B_{s}$ and the function 
$\phi_{0}$. The remainder term ${\mathcal {R}}(t)$ satisfies the 
estimate
\begin{equation}
\nn
\|\mathcal {R}(t)\|_{L^{2}}\longrightarrow 0, \quad \text{as}\,\, 
\,t\to \infty
\end{equation}
\end{theorem}
\begin{proof} We first assume that $\psi_0 \in L^1\cap L^2$. The $L^1$-assumption 
can be removed by $L^2$-continuity of the wave operators that appear below. 
Decompose
$$
\psi(t):=U(t) \psi_0 = P_b(H_1) U(t) \psi_0 +  P_b (H_2,t) U(t) \psi_0
+ R(t).
$$
By construction we clearly have 
\begin{align}
& P_b(H_2,t) U(t) \psi_0 + R(t) \in \Ran(P_c (H_1)),\label{eq:dec2}\\ 
&P_b(H_1) U(t) \psi_0 + R(t) \in \Ran(P_c(H_2,t))\nn
\end{align}
We further write
$$
P_b(H_1) U(t) \psi_0 = \sum^m_{r =1} e^{-i\lambda_r t}
a_r (t) u_r (x)
$$
for some choice of unkown functions $a_{r}(t)$.
Using \eqref{eq:dec2} we obtain, similar to  \eqref{eq:ODE},
$$
\dot a_{r}+ i\,\langle V_2 (\cdot - t\vec
e_{1})\psi(t), u_r\rangle = 0 \text{\ \ for all\ \ }1\le r\le m.
$$
The exponential localization of $u_r$ implies that  
$
|\langle V_2 (\cdot - t\vec e_1) \psi(t), u_{r} \rangle |
\les e^{-\alpha t}$.
Therefore, $a_r (t)$ has a limit  
$
a_r(t) \to A_{r},\quad \text{as}\,\,t\to +\infty
$ 
and thus, 
\begin{equation}
\label{eq:tired}
\Big\| P_b (H_1) U(t) \psi_0 - \sum^m_{r=1} A_r e^{-i\lambda_r t}u_r \Big\|_{L^2}\to 0, 
\quad t\to +\infty. 
\end{equation}
We next define the functions $v_r  = \underset{ t\to +\infty}{\lim}\,  U(t)^{-1}
e^{-i\lambda_r t} u_r.$
The existence of $v_r$ is guaranteed by the existence of the wave
operators (see \cite{Gr} and for an independent proof, Lemma~\ref{le:woper}) 
$$
\Omega^1_{-} = \underset{ t\to +\infty}{s-\lim}\, U(t)^{-1} e^{-itH_{1}}
P_b(H_1)$$  
In this notation,  $v_r = \Omega^1_{-} u_r$.
In particular,
\begin{equation}
\label{eq:1chan}
\Big\|U(t) \big (\sum^m_{r=1} A_{r} v_r\big ) - \sum^m_{r=1} A_r e^{-i\lambda_r
t}u_r \Big\|_{L^2}\to 0, \quad t
\to + \infty.
\end{equation}
We then infer from \eqref{eq:tired} that
\begin{equation}
\label{eq:pbh1}
\Big\|U(t) \big (\sum^m_{r=1} A_{r} v_r\big ) - P_b (H_1) U(t) \psi_0
\Big\|_{L^2}\to 0, \quad t
\to + \infty.
\end{equation}
Similar arguments apply to $P_b(H_2,t) U(t) \psi_0$. More precisely,
we write 
$$
U(t) \psi_0 = P_b (H_2,t) U(t) \psi_0 + \Gamma(t)
=\calg_{-\vec e_{1}}(t)P_b (H_2) \calg_{\vec e_{1}}(t) U(t)
\psi_0 + \Gamma(t).
$$
Therefore,
\begin{equation}
\label{eq:2dec}
 \calg_{\vec e_{1}}(t)U(t) \psi_0 = P_b (H_2) \calg_{\vec e_{1}}(t)
 U(t) \psi_0 + \calg_{\vec e_{1}}(t)\Gamma(t)
\end{equation}
Recall that the function
$$
\tilde \psi(t) = \calg_{\vec e_{1}}(t) U(t) \psi_0$$
is a solution of the problem
\begin{equation}
\label{eq:anduh}
\frac{1}{i} \partial_t \tilde \psi - \Laplace\tilde \psi + V_2(x) 
\tilde\psi + V_1(x +
t\vec e_1) \tilde\psi = 0,\qquad
\tilde\psi|_{t=0} = \calg_{\vec e_{1}}(0)\psi_{0}
\end{equation}
According to \eqref{eq:2dec},
$\tilde\psi(t) = P_b (H_2) \tilde \psi(t) +  \Gamma_1 (t)$,
where $\Gamma_1 (t) = \calg_{\vec e_{1}}(0) \Gamma(t)$.
In particular, 
$$
\Gamma _1(t) \in P_c (H_2) L^{2}.
$$
Decompose 
$$
P_b(H_2) \tilde\psi(t) = \sum^\ell_{s=1} b_s (t) e^{-i \mu_{s}t} w_s
$$ 
for some choice of unknown functions $b_{s}(t)$.
After substituting the decomposition in \eqref{eq:anduh} we obtain
the equations
$$
\dot b_s (t) + i \,\langle V_1 (\cdot + t\vec e_1) \tilde\psi , 
w_s\rangle = 0 \text{\ \ for all\  \ }1\le s\le\ell.
$$
Using exponential localization of $w_s$ we conclude the existence
of the limit $b_s(t) \to B_{s}$ as $t \to +\infty$.
Thus 
$\| P_b (H_2) \tilde\psi(t) - \sum ^\ell_{s=1} B_{s} e^{-i\mu_s t} w_s
\|_{L^2} \to 0,  \quad t\to \infty $.
Equivalently, after applying $\calg_{-\vec e_{1}}(t)$, we have 
\begin{equation}
\label{eq:asdec}
\Big\|  P_b (H_2,t) U(t) \psi_0 - \sum^\ell_{s=1} B_{s}
e^{-i\mu_j t} \calg_{-\vec e_{1}}(t) w_s \Big\|_{L^2}
\to 0.
\end{equation}
We now invoke the existence of the wave operator (see \cite{Gr} and 
Lemma~\ref{le:woper})
$$
\Omega^-_2 = \underset {t\to +\infty }{s-\lim }\, U(t)^{-1} \calg_{-\vec e_{1}}(t) 
e^{-itH_2} P_b (H_2)
$$
which allows us to define $\omega_s := \Omega^-_2 w_s$. Moreover,
\begin{equation}
\label{eq:2chan}
\Big\|U (t) \big (\sum_{s=1}^{\ell} B_{s}\omega_s\big ) - 
\sum_{s=1}^{\ell } B_{s} e^{-i\mu_s t} \calg_{-\vec e_{1}} (t) w_s
\Big\|_{L^2} \to 0,  \quad t\to + \infty.
\end{equation}
It then follows from \eqref{eq:asdec} that
\begin{equation}
\label{eq:pbh2}
\|  P_b (H_2,t) U(t) \psi_0 - U(t) \big (\sum^\ell_{s=1} B_s
\omega_s\big ) \|_{L^2} \to 0, \quad t \to + \infty.
\end{equation}
We now define the function 
\begin{equation}
\label{eq:free}
\phi := \psi_0 - \sum^m_{r=1} A_r v_r - \sum^\ell_{s=1} B_s \omega_s,
\end{equation}
which will lead to the initial data $\phi_{0}$ for the free channel.
We have that
$$
P_b(H_1) U(t) \phi = P_b(H_1) U(t) \psi_0 -
P_b (H_1) U(t) \big (\sum^m_{r=1} A_{r} v_r\big ) - P_b (H_1) U(t)
\big (\sum^\ell_{s=1} B_s \omega_s\big ). 
$$
It follows from \eqref{eq:pbh1} and the identity $P^{2}_{b}(H_{1}) = 
P_{b} (H_{1}) $  that 
\begin{equation}
\label{eq:again1}
\Big\| P_b(H_1) U(t) \psi_0 -
P_b (H_1) U(t) \big (\sum^m_{r=1} A_{r} v_r\big ) \Big\|_{L^{2}} \to 0.
\qquad \text{as}\,\, t\to +\infty
\end{equation}
Furthermore, 
\begin{equation}
\label{eq:again2}
P_b(H_1) \sum^\ell_{s=1} B_s e^{-i\mu_s t} \calg_{-\vec e_{1}}(t) w_j
=\sum^m_{r=1} \sum^\ell_{s=1} B_s e^{-i\mu_s t} \langle
 \calg_{-\vec e_{1}}(t) w_j , u_r \rangle \, u_{r}\to 0
\end{equation}
in the $L^{2}$ sense as $t\to +\infty $, due to the exponential localization of the 
eigenfunctions $u_{r}$.
One concludes from~\eqref{eq:again1}, \eqref{eq:2chan}, and~\eqref{eq:again2} that 
$
\|P_b (H_1) U(t) \phi \|_{L^2} \to 0$.
Similarly, we can show that 
$\|P_b (H_2,t) U(t) \phi\|_{L^2} \to 0$.
Thus, $U(t) \phi$ is asymptotically orthogonal to the bound states of
$H_1$ and $H_2$ and therefore, according to Theorem \ref{thm:main},  satisfies 
the estimate
\begin{equation}
\label{eq:wavedec}
\|U(t) \phi\, \|_{L^2+ L^\infty} \les \langle t \rangle^{-\frac n2}
\|\phi\|_{L^1\cap L^2}.
\end{equation}
In order to be able to apply estimate \eqref{eq:wavedec} one needs to 
verify that $\phi\in L^{1}\cap L^{2}$. By the assumption of the theorem
the function $\psi_{0}\in L^{1}\cap L^{2}$. Thus it remains to check 
this property for the functions $v_{r},\, r=1,..,m$ and $\omega_{s},\, 
s=1,..,\ell$. Since $v_{r}= \Omega_{1}^{-} u_{r}$ and $\omega_{s}= 
\Omega_{2}^{-}w_{s}$ the $L^{2}$ property follows immediately.
The $L^{1}$ property, on the other hand, 
is guaranteed by Lemma~\ref{le:woper} below.
Assuming this lemma for the moment, 
we now consider the expression
$$
e^{-it\Laplace}U(t)\phi = \phi -i \int^t_0 e^{-is\Laplace} \left(V_1(x) +
V_2 (x- s \vec e_1)\right) U(s) \phi\, ds.
$$
The estimate 
\begin{align*}
\int^\infty_t \|e^{-is\Laplace} \left( V_1(x) + V_2(x - t\vec 
e_1)\right)&
U(s) \phi\|_{L^2}
\les (\| V_1\|_{L^2\cap L^\infty} + \|V_2\|_{L^2\cap
L^\infty})\int^\infty_t \| U(s) \phi \|_{L^2+L^\infty} ds\\
&\les t^{-\frac n2+ 1} \left(\|V_1\|_{L^2\cap L^\infty} + \|V_2\|_{L^2\cap
L^\infty}\right)\to 0, \quad \text{as}\,\,\,\, t\to +\infty
\end{align*}
allows us to show the existence of the limit 
$$
\phi_{0}:=\underset{ t\to\infty}{\lim} e^{it\Laplace} U(t) \phi.
$$
It follows that
\begin{equation}
\label{eq:freech}
\|U(t)\phi - e^{it\Laplace} \phi_0\|_{L^2} \to 0, \quad  t\to + \infty.
\end{equation}
Combining \eqref{eq:1chan}, \eqref{eq:2chan}, \eqref{eq:free}, and 
\eqref{eq:free} we infer that 
$$
\Big\|U(t)\psi_0 - \sum^m_{r=1} A_r e^{-i\lambda_r t} u_r -
\sum^\ell_{s=1} B_{s} e^{-i\mu_s t}\calg_{-\vec e_{1}}(t) w_s
 - e^{it\Laplace} \phi_0 \Big\|_{L^2} \to 0,\quad  \text{ 
as }\,\,\,\,t\to + \infty,
$$
as claimed.
\end{proof}

\begin{lemma}
Assume that the potentials $V_{1}(x), V_{2}(x)$ satisfy the following 
conditions 
\begin{equation}
\label{eq:assV}
\sum_{0\le|\gamma|\le n+2}
|\pr_{x}^{\gamma} V_{\kappa}(x)|\le \frac{c(\kappa) }
{\la x\ra^{3[\frac n2] +5}},
\end{equation}
for some positive constants $c(\kappa),\, \kappa=1,2$.
Then the range of the wave operators $\Omega_{1,2}^{-}$ is contained in 
the space of $L^{1}$ functions.
\label{le:woper}
\end{lemma}
\begin{proof}
Without loss of generality we only consider the wave operator 
$$
\Omega_{1}^{-} = \underset{t\to +\infty}{s-\lim}\, U(t)^{-1} e^{-itH_1} 
P_{b}(H_{1}).
$$
Therefore, for an arbitrary $L^{2}$ function $f$
$$
\Omega_{1}^{-} f = \sum_{r=1}^{m} f_{r}\underset{t\to +\infty}{\lim} U(t)^{-1} 
e^{-it H_{1}} u_{r},
$$
where $P_{b}(H_{1}) f =\sum_{r=1}^{m} f_{r} u_{r}$ for some 
constants $f_{r}$. 
It follows from the Duhamel formula for the equation \eqref{eq:2pot} 
generating the evolution operator $U(t)$ that 
\begin{align}
U(t)^{-1} e^{-itH_{1}} u_{r} = &u_{r} + i \int_{0}^{t} U(s)^{-1} V_{2}(\cdot - 
s\vec e_{1}) e^{-isH_{1}} u_{r}\, ds \nn\\ = &u_{r} + i \int_{0}^{t} U(s)^{-1} V_{2}(\cdot - s\vec e_{1}) e^{-i\lambda_{r} s} u_{r}\, ds \label{eq:duhamel}
\end{align}
with \eqref{eq:duhamel} follows since $u_{r}$ is an eigenfunction of 
$H_{1}$ corresponding to an eigenvalue $\lambda_{r}$. The function 
$u_{r}$ is exponentially localized in $L^{2}$ together with its 
$n+2$ derivatives \footnote
{The localization of higher derivatives of $u_{r}$ follows from the localization
of $u_{r}$ stated in \eqref{eq:local} and the equation 
$-\Laplace u_{r} + V_{1}(x) u_{r} = \lambda_{r} u_{r}$  with
potential $V_{1}(x)$ which is bounded together with all its derivatives 
of order $\le (n+2)$.}
$$
\sum_{0\le |\gamma|\le n+2}
\int_{\R^{n}}e^{2\alpha |x|} |\partial_{x}^{\gamma}u_{r}(x)|^{2}\,
dx \le C
$$
for some positive constant $\alpha$ appearing in \eqref{eq:local}.
It follows from \eqref{eq:assV} that the function
$$
G_{r}(s,x): = e^{-i\lambda_{r} s} V_{2}(x-s\vec e_{1}) u_{r} (x)
$$
has the property that for any $k\ge 0$ and multi-index $\gamma$,
$0\le |\gamma|\le n+2$
$$
\|\la x\ra^{k}\partial^{\gamma}_{x} G_{r}(s,\cdot)\|_{L^{2}_{x}}\le 
c(r,|\gamma |, k) \la s\ra^{-3[\frac n2]-5}. 
$$
To prove the desired conclusion it would then suffice to show that 
there exists a positive constant $k$ such that for any function $g(x)$
\begin{equation}
\label{eq:invariance}
\|\la x\ra^{[\frac n2] +1} U(t) g  \|_{L^{2}} \les 
\la t\ra^{3[\frac n2]+3}  \sum_{|\beta|\le n+2}
\|\la x\ra ^{k} \partial_{x}^{\beta} g\|_{L^{2}},\qquad \forall t\ge 0.
\end{equation}
We note here that $U(t)g$ denotes the solution of the problem
\begin{equation}
\label{eq:again}
\frac 1i \partial_{t} \psi - \Laplace \psi + V(t,x) \psi  = 0, \qquad 
V(t,x) = V_{1}(x) + V_{2}(x-t\vec e_{1})
\end{equation}
evaluated at time $t$ with initial data function $g$ given at time $t=0$. 
The required estimate should, however, involve the expression $U(t)^{-1} g$ which stands 
for the solution of the above problem at time $t=0$ with initial 
data $g$ given at time $t$. Yet it is not difficult to see that the 
two problems are almost equivalent. 
Therefore, for simplicity we shall prove estimate~\eqref{eq:invariance}. 
We note that the estimates of the type \eqref{eq:invariance} for problems with 
time independent potentials are well-known. They have been proved 
in the paper by Hunziker \cite{Hu}. In the time-dependent case the 
argument is essentially the same.
More precisely, define the functions
$$
\Phi_{j,|\gamma|}(t):=\sum_{j'=0}^{j}\sum_{|\gamma'|=0}^{|\gamma|}
\| \la x\ra^{j'} \pr_{x}^{\gamma'} U(t) g\|_{L^{2}}
$$
for any index $j\ge 0$ and any multi-index $\ga$.
Using equation \eqref{eq:again} we obtain that
$$
\frac d{dt} \| \la x\ra^{j} \pr_{x}^{\ga} U(t) g\|_{L^{2}}^{2} = i 
(-1)^{|\gamma|}\la \big [\Laplace - V(t,x),  
\la x\ra^{j}\pr_{x}^{\ga}\la x\ra^{j}\pr_{x}^{\ga}\big ] U(t)g, 
U(t)g\ra 
$$
Computing the commutator we obtain the recurrence relation
\begin{align*}
\Phi_{j,|\gamma|}(t)\les& \Phi_{j,|\gamma|}(0) +\la t\ra^{2}
\sum_{|\gamma'|\le 2|\gamma|} 
\Big \|\frac{\la x\ra}{\la t\ra} \pr_{x}^{\gamma'} V \Big\|_{L^{\infty}_{t,x}} 
 \sup_{0\le \tau \le t} \Phi_{j-1,|\gamma|+1}(\tau)\le\\ 
&C(V) \bigg (\sum_{k=0}^{j-1} \la t\ra^{2k}\Phi_{j-k,|\gamma|+k}(0) +
\la t\ra^{2j} \Phi_{0,|\gamma|+j}(\tau)\bigg ),
\end{align*}
where $C(V)$ is a constant depending on 
\begin{equation}
\label{eq:potV}
\sum_{|\gamma'|\le 2(|\gamma|+j-1)} 
\Big\|\frac{\la x\ra}{\la t\ra} \pr_{x}^{\gamma'} V \Big\|_{L^{\infty}_{t,x}} 
\end{equation}
In addition, differentiating the equation \eqref{eq:again}
$(|\gamma|+j)$ times with respect to $x$ and using the standard $L^{2}$ estimate,
we have 
$$
\Phi_{0,|\gamma|+j}(\tau)\le C(V) (1+|\tau|^{|\gamma|+j}) 
\Phi_{0,|\gamma|+j}(0)
$$
Therefore,
$$
\Phi_{j,|\gamma|}(t) \le C(V) (1+|t|)^{3j+|\gamma|} 
\Phi_{j,|\gamma|+j}(0).
$$
Now setting $j=[\frac n2] +1 $ and $|\gamma| =0$ we obtain the 
desired estimate \eqref{eq:invariance} with $k=[\frac n2] +1$.
Observe that the assumption \eqref{eq:assV} controls the constant $C(V)$ 
in \eqref{eq:potV} for the potential 
$V(t,x) = V_{1}(x) + V_{2}(x-t\vec e_{1})$.
\end{proof}

\section{$L^{\infty}$ estimates}
\label{sec:infty}

In this section we develop a general simple scheme allowing to convert 
the $L^{2}+L^{\infty}$ estimates developed above into the true 
dispersive $L^{\infty}$ estimates.

\begin{prop}
\label{pr:transit}
Let $\psi$ be a solution of the Schr\"odinger equation 
\begin{equation}
\nn
\frac 1i \partial_{t} \psi - \Laplace \psi + V(t,x) \psi = 0,\qquad
\psi |_{t=0} = \psi_{0}
\end{equation}
with a time-dependent potential $V(t,x)$ satisfying the condition
\begin{equation}
\sup_{t} \|V(t,\cdot)\|_{L^{1}\cap L^{2}} <\infty,\qquad
\sup_{t} \|\hat V(t,\cdot)\|_{L^{1}} <\infty.
\label{eq:assVinf}
\end{equation}
Here $\hat V(t,\cdot)$ denotes the Fourier transform of $V$ only with 
respect to the spatial variable $x$.
Assume that $\psi$ obeys the estimate 
\begin{equation}
\label{eq:estp}
\|\psi(t)\|_{L^{2}+L^{\infty}}\les \la t\ra^{-\frac n2}\|\psi_{0}\|_{L^{1}\cap L^{2}}.
\end{equation}
Then $\psi$ also satisfies the $L^{\infty}$ estimate
\begin{equation}
\label{eq:espi}
\|\psi(t)\|_{L^{\infty}}\les |t|^{-\frac n2} \|\psi_{0}\|_{L^{1}\cap L^{2}}.
\end{equation}
\end{prop}
\begin{proof}
For simplicity consider the case of dimension $n=3$. 
The $n>3$ dimensionial case can be treated by considering 
$k=[\frac n2] + 2$ terms in the Duhamel expansion and 
repeatedly exploiting the cancellation property below 
(see \cite{JSS}).

\noindent We have
\begin{align}
& \psi(t) = e^{it {\Laplace}}\psi_{0} - i \int_{0}^{t} 
e^{i(t-s){\Laplace}} V(s,\cdot) \psi(s)\, ds \nn\\
&= e^{it  {\Laplace}}\psi_{0} - i \int_{0}^{t} 
e^{i(t-s){\Laplace}} V(s,\cdot) e^{is {\Laplace}}\psi_{0}\, ds
- \int_{0}^{t}\int_{0}^{s}e^{i(t-s){\Laplace}} V(s,\cdot) 
e^{i(s-\tau){\Laplace}} V(\tau,\cdot) \psi(\tau)\, d\tau\, ds.
\label{eq:duha}
\end{align}
We recall the cancellation property 
\begin{equation}
\label{eq:canc}
\sup_{s}\|e^{is{\Laplace}} V(t,\cdot) e^{-is{\Laplace}} 
f\|_{L^{p}}\les \|\hat V(t,\cdot)\|_{L^{1}} \|f\|_{L^{p}}, \quad\forall 
p\in [1,\infty ]
\end{equation}
which was used in \cite{JSS}. This property can be checked directly for pure 
exponentials, and then one writes the potential as superposition of those.  
We can then easily estimate the first two terms in~\eqref{eq:duha} by the 
desired bound $|t|^{-\frac n2} \|\psi_{0}\|_{L^{1}}$. For the second one, 
split the integration according to $0<s<1$, $1<s<t-1$, and~$t-1<s<t$.  
For the values $t\ge 2$ we split the term 
\bea
&&\int_{0}^{t}\int_{0}^{s}e^{i(t-s){\Laplace}} V(s,\cdot) 
e^{i(s-\tau){\Laplace}} V(\tau,\cdot) \psi(\tau)\, d\tau\, 
ds = \nn \\
&& \Big (\int_{0}^{t-2}\int_{\tau+1}^{t-1} + 
\int_{0}^{t-2}\int_{\tau}^{\tau +1} +\int_{0}^{t-2}\int_{t-1}^{t}
+ \int_{t-2}^{t}\int_{\tau}^{t}\Big )
e^{i(t-s){\Laplace}} V(s,\cdot) 
e^{i(s-\tau){\Laplace}} V(\tau,\cdot) \psi(\tau)\, ds\, d\tau. \label{eq:unterteil}
\eea
With the help of the assumption \eqref{eq:estp} we estimate
\bea
\|\int_{0}^{t-2}\int_{\tau +1}^{t-1} &&
e^{i(t-s){\Laplace}} V(s,\cdot) 
e^{i(s-\tau){\Laplace}} V(\tau,\cdot) \psi(\tau)\, ds\, 
d\tau\|_{L^{\infty}} \nn \\ 
&& \les \int_{0}^{t}\int_{\tau }^{t} 
\frac {\sup_{s} \|V(s,\cdot)\|_{L^{1}}}{\la t-s\ra^{\frac 32} }
\frac {\sup_{\tau} \|V(\tau,\cdot)\|_{L^{1}\cap L^{2}}}{\la s-\tau\ra^{\frac 32}}
\la \tau\ra^{-\frac 32} ds\,d\tau  \|\psi_{0}\|_{L^{1}\cap 
L^{2}} \nn \\
&& \les \la t\ra ^{-\frac 32} \|\psi_{0}\|_{L^{1}\cap L^{2}}. \label{eq:einfach}
\eea
Furthermore, in view of \eqref{eq:canc}, the decay of the free evolution, 
and~\eqref{eq:estp}, 
\begin{align*}
\Big\|\Big (\int_{0}^{t-2}\int_{\tau }^{\tau+1} + \int_{0}^{t-2}\int_{t-1 
}^{t}\Big )&
e^{i(t-s){\Laplace}} V(s,\cdot) 
e^{i(s-\tau) {\Laplace}} V(\tau,\cdot) \psi(\tau)\, ds\, 
d\tau \Big\|_{L^{\infty}}\\ 
& \les \int_{0}^{t} \sup_{s}\|\hat 
V(s,\cdot)\|_{L^{1}} \frac {\sup_{\tau}\|V(\tau,\cdot)\|_{L^{1}}}{\la 
t-\tau\ra^{\frac 32}} \la \tau \ra^{-\frac 32}\, d\tau 
\|\psi_{0}\|_{L^{1}\cap L^{2}}\\ 
& \qquad\qquad \les \la t\ra^{-\frac 32} 
\|\psi_{0}\|_{L^{1}\cap L^{2}}.
\end{align*}
Finally,
\begin{align}
\Big\|\int_{t-2}^{t}\int_{\tau}^{t}&
e^{i(t-s){\Laplace}} V(s,\cdot) 
e^{i(s-\tau){\Laplace}} V(\tau,\cdot) \psi(\tau)\, ds\, d\tau 
\Big\|_{L^{\infty}} \nn \\ 
& \les \la t\ra^{-\frac 32} \sup_{s}\|\hat V(s,\cdot)\|_{L^{1}} 
\sup_{\tau}\|V(\tau,\cdot)\|_{L^{1}\cap L^{2}}
\int_{t-2}^{t}\int_{\tau}^{t} |t-\tau|^{-\frac 32} ds\,d\tau \, 
\|\psi_{0}\|_{L^{1}\cap L^{2}}\nn \\ 
&\qquad \qquad \les \la t\ra^{-\frac 32}
\|\psi_{0}\|_{L^{1}\cap L^{2}}.\label{eq:trip}
\end{align}
The remaining estimate in the case $t < 2$ can be carried out precisely 
in the same fashion as in \eqref{eq:trip} with an integral lower limit $t-2$ 
replaced by $0$.
\end{proof}
\section{Estimates for the inhomogeneous equation}

Inspection of the proof of Theorem \ref{thm:main} indicates that we can 
extend the $L^{\infty}$ estimates also to the case of the 
inhomogeneous equation
\begin{equation}
\label{eq:inhom}
\frac 1i \partial_{t }\psi - \Laplace \psi + \sum_{\kappa =1}^{m} 
V(x-\vec v_{\kappa} t) \psi = F(t,x), \qquad \psi |_{t=0}=\psi_{0}.
\end{equation}
This will be particularly useful in the forthcoming application~\cite{RSS} to 
the question of stability of multi-soliton states.
Again, for simplicity consider the case of two potentials, $m=2$. Then 
we have the following 

\begin{theorem}
\label{thm:inhom}
Let $V_\kappa$ be as in Definition~\ref{def:chargetrans} and
let $\psi$ be a solution of equation \eqref{eq:inhom} with the 
right-hand side~$F$ satisfying the condition
\begin{equation}
\label{eq:condF}
\||F\|| := \sup_{t}\Big \{\int_{0}^{t} 
\|F(\tau,\cdot)\|_{L^{1}}\,d\tau  + \la t\ra^{\frac n2+1} 
\|F(t,\cdot)\|_{L^{2}}\Big \}<\infty. 
\end{equation}
Assume that the following decay of the projections of $\psi$ onto the bound states 
corresponding to the Hamiltonians $H_{1}$ and $H_{2}$:
\begin{equation}
\label{eq:decbound}
\|P_{b}(H_{1} )\psi (t,\cdot)\|_{L^{2}} + \|P_{b}(H_{2},t) 
\psi(t,\cdot)\|_{L^{2}}\le B\la t\ra ^{-\frac n2},
\end{equation}
cf.\ Definition \ref{def:asymp}, holds for some positive constant $B$. 
Then we have the $L^{2}+L^{\infty}$ estimate
\begin{equation}
\label{eq:decinf}
\|\psi(t,\cdot)\|_{L^{2} + L^{\infty}}\les  \la t\ra^{-\frac n2} \Big 
(\|\psi_{0}\|_{L^{1}\cap L^{2}} + \||F\|| + B\Big )
\end{equation}
\end{theorem}
\begin{proof}
The proof requires repeating the arguments from Sections~\ref{sec:dec}. 
As in Theorem~\ref{thm:main}, we 
introduce the three channel decomposition and represent the solution via 
the respective Duhamel formulas. The novel terms generated by the 
inhomogeneous term $F$ can be written in the form
\begin{align*}
&I_{1}(t)=\chi_{1}\int_0^{t} e^{-i(t-s)H_{1}}P_{c}(H_{1}) F(s,\cdot)\, 
ds\\
&I_{2}(t)=\chi_{2}\calg_{-\vec e_{1}}(t)\int_0^{t} 
e^{-i(t-s)H_{2}}P_{c}(H_{2}) \calg_{\vec e_{1}}(s) F(s,\cdot)\, ds\\
&I_{3}(t)=\chi_{3}\int_0^{t} e^{-i(t-s)\Laplace} F(s,\cdot)\, ds,
\end{align*}
which should be compared with \eqref{eq:DU1}, \eqref{eq:DU2mov}, and~\eqref{eq:DU3},
respectively. 
Observe that by our assumption~\eqref{eq:decbound} the projections 
of~$\psi$ onto the bound states corresponding to~$H_{1}$ and~$H_{2}$ already
satisfy the desired decay.
The terms $I_{1}, I_{2}, I_{3}$ are estimated in a similar manner. 
Consider for instance the $L^{2}+L^{\infty}$ norm of $I_{1}(t)$:
\begin{align*}
\|I_{1}(t)\|_{L^{2}+L^{\infty}}&\les \int_{0}^{\frac t2} \la t\ra^{-\frac n2} 
\|e^{-i(t-s)H_{1}}P_{c}(H_{1}) F(s,\cdot)\|_{L^{2}+L^{\infty}}\, ds + 
\int_{\frac t2}^{t} \|e^{-i(t-s)H_{1}}P_{c}(H_{1}) 
F(s,\cdot)\|_{L^{2}} ds \\ &\les 
\la t\ra^{-\frac n2}\int_{0}^{t} \|F(t,\cdot\|_{L^{1}\cap L^{2}}\,ds + 
\int_{\frac t2}^{t}\|F(s,\cdot )\|_{L^{2}}\,ds\les \la t\ra^{-\frac n2}
\||F\||
\end{align*}
as desired.
\end{proof}

\section{Decay estimates: Systems with a single matrix potential}
\label{sec:compl}

\subsection{Three dimensions}

In the following sections we develop the decay estimates for 
the charge transfer model with matrix potentials. 
These systems arise in the study of stability of multi-soliton 
states on nonlinear Schr\"odinger equations~\cite{RSS}. 

We first
investigate the question of the decay estimates for problems 
with a single time-independent matrix potential.
We show that the approach developed by Rauch~\cite{R} for 
scalar equations can be naturally extended to systems. 
The method relies on meromorphic continuation of the resolvent of
the Hamiltonian across the spectrum and requires exponential 
decay of the potential at infinity. 
Then, following the ideas we have developed in the previous sections
for the scalar charge transfer models, we establish the 
decay estimates for solutions of the matrix charge transfer problem
satisfying appropriate "asymptotic orthogonality" conditions.
Our approach to dispersive estimates for systems is more direct
than the one of Cuccagna~\cite{Cuc}, who extended Yajima's $L^p$-boundedness
result for the wave operators to the case of systems.
On the other hand, \cite{Cuc} does not require exponential decay of
the potentials as we do. 

\noindent 
In this section we consider the case of three dimensions, $n=3$.
The higher dimensional situation is discussed in the next section.

\noindent In this section $H=-\Lapl+\mu$ where $\mu>0$ and $H$ is an operator 
on~$L^2(\R^3)$. Set
\[ B = \bm 0 & H \\ -H & 0 \endm,\quad V= \bm 0 & V_1 \\ -V_2 & 0 
\endm,\quad  A = B+V=\bm 0 & H+V_1 \\ -H-V_2 & 0 \endm.\] 
By means of the matrix $J=\bm 0 & 1 \\ -1 & 0 \endm$ one can also 
write
\[ B = \bm H & 0 \\ 0 & H \endm J,\quad A = \bm H+V_1 & 0 \\ 0 & 
H+V_2 \endm J.\]
Since $B^*=-B$ it follows that $\spec(B)\subset i\R$. 
One checks that for $\Re z\ne0$
\bea
(B-z)^{-1} &=& -(B+z)\bm (H^2+z^2)^{-1} & 0 \\ 0 & (H^2+z^2)^{-1} 
\endm  \nn \\
&=& -\bm (H^2+z^2)^{-1} & 0 \\ 0 & (H^2+z^2)^{-1} \endm (B+z) 
\label{eq:Bspec} \\
(A-z)^{-1} &=& (B-z)^{-1} - (B-z)^{-1} W_1 \Bigl[1+W_2 J (B-z)^{-1} 
W_1\Bigr]^{-1} W_2 J (B-z)^{-1},
\label{eq:grill}
\eea
where \eqref{eq:grill} also requires the expression in brackets to be 
invertible, and with
\[
 W_1 = \bm |V_1|^{\half} & 0 \\ 0 & |V_2|^{\half} \endm,\quad W_2 = 
\bm |V_1|^{\half}\sign(V_1) & 0 \\ 0 & |V_2|^{\half}\sign(V_2) \endm.
\]
It follows from \eqref{eq:Bspec} that $\spec(B) = 
(-i\infty,-i\mu]\cup [i\mu, i\infty)\subset i\R$. 
We will assume that there exists $\eps_0>0$ so that for all $x\in\R^3$
\begin{equation} \label{eq:Vdec}
|V_j(x)|\le Ce^{-\eps_0|x|} \text{\ \ for\ \ }j=1,2.
\end{equation}
Then it follows
from Weyl's theorem, see Theorem~XIII.14 in~\cite{RS4}, and the 
representation~\eqref{eq:grill} for 
the resolvent of~$B$, that 
$\spec_{ess}(A)=\spec_{ess}(B)=(-i\infty,-i\mu]\cup [i\mu, 
i\infty)\subset i\R$. This was observed by Grillakis~\cite{Grill}. 
Moreover, \eqref{eq:grill} implies via the analytic Fredholm 
alternative that $(A-z)^{-1}$
is a meromorphic function 
in~$\Compl\setminus (-i\infty,-i\mu]\cup [i\mu, i\infty)$. 
We need to make several further assumptions on~$A$ 
that are motivated by applications to NLS, but we first
switch to a different (but equivalent) way of writing
systems.
Let $\psi$ be a solution of some NLS. Studying stability questions 
for~$\psi$, as we do in~\cite{RSS}, leads to a system for the 
variation~$\delta\psi$
 that can be written either for
the vector of real and imaginary parts of~$\delta\psi$,  or for the 
vector 
$\binom{\delta\psi}{\overline{\delta\psi}}$.  
For the purposes of obtaining~\eqref{eq:grill}  
it was convenient for us to follow the former
convention, but from now on it will be advantageous to use the 
latter. 
It is clear that these two different ways of writing 
the system are equivalent by means of a change
of coordinates in~$\C^2$. More precisely, one has the following lemma 
whose proof is left to the reader.

\begin{lemma}
\label{lem:sys_equiv}
Let 
\[ \binom{\psi_1}{\psi_2} = \bm 1 & i\\ 1 & -i \endm 
\binom{f_1}{f_2}.\]
Then, with $H=-\Lapl+\mu$, 
\begin{equation}
\label{eq:s_1} \partial_t \binom{f_1}{f_2} + \bm 0 & -H-V_1 \\ H+V_2 
& 0 \endm \binom{f_1}{f_2} =0
\end{equation}
holds if and only if
\begin{equation}
\label{eq:s_2}
\frac{1}{i}\partial_t \binom{\psi_1}{\psi_2} + \bm H+U & -W \\ W & 
-H-U \endm \binom{\psi_1}{\psi_2}
=0
\end{equation}
where $U=\half(V_1+V_2)$, $W=\half(V_1-V_2)$. 
\end{lemma}

Abusing notation, we will from now on write
\begin{equation}
\label{eq:Bdef}
 B = \bm H & 0 \\ 0 & -H  \endm,\quad V= \bm U & -W \\ W & -U 
\endm,\quad  A = B+V=\bm H+U & -W \\ W & -H-U \endm.
\end{equation}
with $U,W$ real-valued. We now specify the spectral assumptions
we impose on~$A$. 

\begin{defi}
\label{def:spec_ass} Let $A=B+V$ with $V$ exponentially decaying.
We call the operator $A$ on $\Hil:=L^2(\R^3)\times L^2(\R^3)$ {\em 
admissible} provided
\begin{itemize}
\item $\spec(A)\subset \R$ and $\spec(A)\cap (-\mu,\mu) = 
\{\omega_\ell\:|\: 0\le\ell \le M\}$,
      where $\omega_0=0$ and all $\omega_j$ are distinct eigenvalues. 
There are no eigenvalues
      in $\spec_{ess}(A)$.  
\item For $1\le \ell\le M$, 
$L_\ell:=\ker(A-\omega_\ell)^2=\ker(A-\omega_\ell)$, and 
$\ker(A)\subsetneq \ker(A^2)=\ker(A^3)=:L_0$. Moreover, these spaces 
are finite dimensional.
\item The ranges $\Ran(A-\omega_\ell)$ for $1\le\ell\le M$ and 
$\Ran(A^2)$ are closed.
\item The spaces $L_\ell$ are spanned by exponentially decreasing 
    functions in $\Hil$ (say with bound $e^{-\eps_0|x|}$).
\item All these assumptions hold as well for the adjoint $A^*$. We 
denote the corresponding (generalized)
eigenspaces by $L_\ell^*$.
\item The points $\pm \mu$ are not resonances of $A$ (we will define 
a resonance below, see Remark~\ref{rem:reson}).
\end{itemize}
\end{defi}
Definition \ref{def:spec_ass} is motivated by applications to the 
questions of stability of soliton solutions of NLS 
(see e.g. \cite{BuPe1}, \cite{RSS}).
 
\noindent Henceforth it will be automatically assumed that $A$ is 
admissible.

\begin{lemma}
\label{lem:split} There is a direct sum decomposition
\begin{equation}
\label{eq:split}
\Hil = \sum_{j=0}^M L_j + \Bigl(\sum_{j=0}^M L_j^*\Bigr)^\perp
\end{equation}
This means that the individual summands are linearly independent. 
The decomposition \eqref{eq:split} is invariant under~$A$. 
Let~$P_c$ denote the projection
onto $\Bigl(\sum_{j=0}^M L_j^*\Bigr)^\perp$ which is induced by the 
splitting~\eqref{eq:split}, 
and set $P_b=Id-P_c$. Then $AP_c=P_cA$, and there exist numbers 
$c_{ij}$ so that
\begin{equation}
\label{eq:Pbexpl}
 P_b f = \sum_{i,j} \phi_j c_{ij} \la f,\psi_i \ra \text{\ \ for all\ 
\ }f\in\Hil
\end{equation}
where $\phi_j,\psi_i$ are exponentially decreasing functions.
\end{lemma}
\begin{proof}
One has
\begin{equation}
\label{eq:randurch} 
\Bigl(\sum_{j=0}^M L_j^*\Bigr)^\perp = \Ran(A^2)\cap \bigcap_{j=1}^M 
\Ran(A-\omega_j)
\end{equation}
by the assumption that the ranges are closed. We verify first that 
the summands in~\eqref{eq:split}
are linearly independent. To this end suppose
\begin{equation}
\label{eq:fdef}
 f=\sum_{j=0}^M f_j,\quad f_j\in L_j,\;0\le j\le M, \quad f\in 
\Bigl(\sum_{j=0}^M L_j^*\Bigr)^\perp.
\end{equation}
Then by~\eqref{eq:randurch}, $f=A^2 g_0$ and $f=(A-\omega_j) g_j$ 
for some $g_0\in \Dom(A^2)$, $g_j\in\Dom(A)$. Note that for $j\ge1$, 
$f_j\in \ker(A-\omega_j) \Longrightarrow f_j \in \Dom(A^s)$ for all 
$s\ge1$. 
Thus \eqref{eq:fdef} implies that
\[ A^2 f = \sum_{j=1}^M A^2 f_j, \text{\ \ and\ \ } A^2 \prod_{j=1}^M 
(A-\omega_j) f =0.\]
Therefore, by our assumption on (generalized) eigenspaces 
\[ (A-\omega_\ell)^2 \prod_{j\ne\ell} (A-\omega_j) A^2 g_\ell =0 
\Longrightarrow (A-\omega_\ell)\prod_{j\ne\ell}(A-\omega_j)A^2 
g_\ell =0,\]
and similarly,
\[ A^4 \prod_{j=1}^M (A-\omega_j) g_0 =0 \Longrightarrow A^2 
\prod_{j=1}^M (A-\omega_j) g_0 =0.\]
Hence $f$ satisfies the equations
\[ \prod_{j=1}^M (A-\omega_j) f =0, \quad A^2 
\prod_{j\ne\ell}(A-\omega_j) f =0 \text{\ \ for all\ \ }
1\le \ell\le M.\]
Using that $L_j\cap L_k=\emptyset$ for $j\ne k$, one concludes 
inductively that $\prod_{j>j_0} (A-\omega_j) f=0$ for any $j_0$. 
Indeed, setting $A_0:=A^2$ and $A_j:= A-\omega_j$, we just showed 
that
$\prod_{j\ne\ell} A_j f =0$ for all $0\le \ell\le M$. Thus also 
$\prod_{j\ne\ell,k} A_j f =0$
for any $\ell\ne k$. Continuing inductively one finds that $f=0$. By 
the linear independence of the 
spaces $\{L_j\}_{j=0}^M$ this implies that $f_0=f_1=\ldots = f_M=0$. 
Thus,
\[ \Bigl(\sum_{j=0}^M L_j\Bigr) \cap \Bigl(\sum_{j=0}^M 
L_j^*\Bigr)^\perp = \{0\}.\]
Dually, one finds that
 \[ \Bigl(\sum_{j=0}^M L_j\Bigr)^\perp \cap \Bigl(\sum_{j=0}^M 
L_j^*\Bigr) = \{0\}\]
which means that the right-hand side of \eqref{eq:split} is dense 
in~$\Hil$.
Since the sum of a closed space with a finite dimensional space is 
again closed, 
we have proved that~\eqref{eq:split} is indeed a direct sum 
decomposition.
Hence every $f\in\Hil$ can be written as $f=g+h$ where 
$g\in\sum_{j=0}^M L_j$, and 
$h\in \Bigl(\sum_{j=0}^M L_j^*\Bigr)^\perp$. We set $P_c f = h$, $P_b 
f=g$. 
Since $AL_j\subset L_j$ for $0\le j\le M$, and $A\Bigl( 
\Bigl(\sum_{j=0}^M L_j^*\Bigr)^\perp\Bigr)
\subset \Bigl(\sum_{j=0}^M L_j^*\Bigr)^\perp$, one has $AP_c=P_c A$. 
Denote the orthogonal projection onto $\sum_{j=0}^M L_j^*$ by~$Q$. 
Then there is an isomorphism
\[ S:\left\{ \begin{array}{lll} Qf \mapsto P_b f \\ \sum_{j=0}^M 
L_j^* \to \sum_{j=0}^M L_j
             \end{array}
\right. \]
and $P_b = SQ$ by construction. As a map between finite dimensional 
spaces $S$ is given by
a matrix $\{c_{ij}\}$. By assumption, these finite dimensional 
subspaces of~$\Hil$ are
spanned by exponentially decaying functions, and~\eqref{eq:Pbexpl} 
follows. 
\end{proof}

Observe that \eqref{eq:grill} implies that there exists $\lambda>0$ 
large so that 
\[ \|(A-(i\lambda + z))^{-1} \| \le |\Im z|^{-1} \text{\ \ for all\ \ 
}\Im z\ge 0.\]
Unless indicated otherwise, $\|\cdot\|$ refers to the operator norm 
on~$\Hil$. 
This allows one to define the (quasi-bounded) semigroup $e^{itA}$ by 
means of the Hille--Yoshida theorem
and $\|e^{itA}\|\le e^{t\lambda}$. We now make one more requirement 
for $A$ to be admissible,
namely the {\em stability assumption} 

\begin{equation}
\label{eq:stab} 
\sup_{t\in\R}\|e^{itA} P_c\| < \infty.
\end{equation}

\noindent Under this assumption the semigroup is bounded on all 
functions that have no
component in the space~$L_0$ under the splitting~\eqref{eq:split}. 
Otherwise, it can
grow at most like~$t$.

\noindent Assumption \eqref{eq:stab} is related to the notion
of {\it {linear stability}} in the context of stability of soliton
solutions of NLS (see e.g. \cite{W1}, \cite{St}).

\noindent 
Our goal in this section is to prove the following theorem.

\begin{theorem}
\label{thm:vec_dec}

\vskip 1pc
\noindent \,

\begin{itemize}
\item
Let $A$ be admissible as in Definition~\ref{def:spec_ass} together 
with the stability
assumption~\eqref{eq:stab}. Then
\begin{equation}
\label{eq:vec_dec}
 \|e^{itA} P_c f\|_{L^2+L^\infty} \les \la t\ra^{-\frac32} 
\|f\|_{L^1\cap L^2}
\end{equation}
for all $t$.  
\item
If in addition the matrix potential satisfies $\|\hat{V}\|_{L^1}<\infty$, 
then one has the bound
\begin{equation}
\label{eq:vec_dec_infty}
 \|e^{itA} P_c f\|_{L^\infty} \les | t|^{-\frac32} 
\|f\|_{L^1\cap L^{2}}
\end{equation}
for all $t$.
\end{itemize}
\end{theorem}

The norms here are given by Definition \ref{def:2+infty} 
$$
\| f \|_{L^2 + L^\infty} := {\inf}_{f=h+g} (\| h \|_{L^2} + 
\|g\|_{L^\infty})
$$
with the dual space   $(L^2 + L^\infty)^* = L^1\cap L^2$.

The proof of Theorem~\ref{thm:vec_dec} will follow an old strategy of 
Rauch~\cite{R}, see
also~\cite{DMT} and~\cite{V},  that is
based on representing $e^{itA}$ as the inverse Laplace transform of 
the resolvent,
and then to move the contour across the spectrum. This only gives 
local $L^2$ decay,
and we then use an observation of Ginibre~\cite{Gin} to pass from 
this to the bound~\eqref{eq:vec_dec}.
In the following we prepare the proof of Theorem~\ref{thm:vec_dec} by 
means of
several technical lemmas. We will assume from now on that $A$ is as 
in Theorem~\ref{thm:vec_dec}.

\noindent
The next lemma is a version of von Neumann's ergodic theorem.

\begin{lemma}
\label{lem:ergodic}
For any $\omega\in\R$,
\begin{equation}
\label{eq:ergod_lim}
 \frac{1}{T} \int_0^T e^{-i\omega t} e^{itA} P_c \,dt \to 0 \text{\ \ 
as\ } T\to\infty
\end{equation}
in the strong sense.
\end{lemma}
\begin{proof} 
Let $\{\omega_j\}_{j=0}^M$ be as in Definition~\ref{def:spec_ass}. We 
first 
prove~\eqref{eq:ergod_lim} for $\omega$ equal to one of the 
$\omega_j$, $j\ne0$. Let
$f\in\Ran P_c$. Then for every $1\le j\le M$ there is $g_j\in 
\Dom(A)$ so that $f=(A-\omega_j)g_j$,
see~\eqref{eq:randurch}. Note that $g_j$ is unique up to addition of 
an element of~$L_j$. Thus
\begin{equation}
\label{eq:simp_int}
\frac{1}{T} \int_0^T e^{-i\omega_j t} e^{itA} f\,dt = \frac{1}{iT} 
\Big(e^{iT(A-\omega_j)}g_j-g_j\Big).
\end{equation}
It remains to show that we can choose $g_j$ so that~\eqref{eq:stab} 
applies. To do so, 
use~\eqref{eq:split} to write 
 $g_j=\sum_{k=0}^M h_k + \tilde{h}$ where $h_k\in L_k$ and 
$\tilde{h}\in (\sum_{k=0}^M L_k^*)^\perp=:\tilde{L}$. 
 We may assume that $h_j=0$. Then $Ag_j=\sum_{k=0}^M A h_k + 
A\tilde{h}$ and by $A$-invariance 
of~\eqref{eq:split}, $Ah_k \in L_k$, $A\tilde{h}\in \tilde{L}$. On 
the other hand,
\[ Ag_j =f+\omega_j g_j=\sum_{k=0}^M \omega_j\,h_k + \omega_j 
\tilde{h} + f. \]
The last two terms lie in $\Ran(P_c)=\tilde{L}$, whereas clearly 
$\omega_j\,h_k\in L_k$ for
all $0\le k\le M$. By uniqueness of such a representation, 
$Ah_k=\omega_j h_k$ for all $0\le k\le M$.
But then 
\bea 
h_k &\in& \ker(A-\omega_j)\cap \ker(A-\omega_k) = \{0\} \text{\ \ 
for all\ \ }k\ne j, \nn \\
h_0 &\in& \ker(A-\omega_j)\cap \ker(A^2) = \{0\}. \nn
\eea
Since $h_j=0$ by assumption, it follows that $g_j= \tilde{h} \in 
\Ran(P_c)$, as desired. 
Hence the right--hand side of \eqref{eq:simp_int} is $O(T^{-1})$ and 
we are done with this case.
Now consider the case $\omega=\omega_0=0$. We write $f=A^2 g$ for 
some $g\in\Dom(A^2)$.
By~\eqref{eq:split} one has as before $Ag=\sum_{k=0}^M h_k + 
\tilde{h}$ where $h_k\in L_k$ and 
$\tilde{h}\in \tilde{L}$. We claim that $-h_0+ Ag\in \Ran(P_c)$ and 
$Ah_0=0$. If so, then
one has $f=A(-h_0+ Ag)$ which implies that 
\begin{equation}
\nn
\frac{1}{T} \int_0^T e^{itA} f \,dt = \frac{1}{iT} 
\Big(e^{iTA}(-h_0+Ag)-(-h_0+Ag)\Big),
\end{equation}
cf.~\eqref{eq:simp_int}, and the right--hand side is again 
$O(T^{-1})$ by~\eqref{eq:stab}. 
To prove the claim, simply observe that
\[ f=A^2 g = \sum_{k=0}^M Ah_k + A\tilde{h} \Longrightarrow 
\sum_{k=0}^M Ah_k \in \tilde{L}.\]
Since $Ah_k\in L_k$ one concludes from Lemma~\ref{lem:split} that 
$Ah_k=0$ for all $0\le k\le M$.
If $1\le k\le M$, then $h_k\in L_0\cap L_k=\{0\}$ so that $h_k=0$. 
But this is precisely the claim,
and we are done with the case $\omega=0$ as well.  

\noindent It remains to consider $\omega\in\R\setminus 
\{\omega_j\}_{j=0}^M$. Firstly, if 
$\omega\in (-\mu,\mu)$, then $(A-\omega)^{-1}$ is a bounded 
operator on~$\Hil$ by assumption. Hence, 
\[ 
\frac{1}{T} \int_0^T e^{-i\omega t} e^{itA} P_cf\,dt = 
(A-\omega)^{-1}\frac{1}{iT} \Big(e^{iT(A-\omega)}P_cf-P_cf\Big),
\]
and the right--hand side is $O(T^{-1})$ by \eqref{eq:stab}. Finally, 
assume that $|\omega|\ge\mu$
so that $\omega$ belongs to the essential spectrum of~$A$. This case 
is different from
the previous ones because one cannot expect to obtain the $O(T^{-1})$ 
bound. Since 
\[ \sup_{T}\left\| \frac{1}{T} \int_0^T e^{-i\omega t} e^{itA} P_cf 
\,dt \right\| < \infty\]
for all $f$ by \eqref{eq:stab} it suffices to prove that
\[ \frac{1}{T} \int_0^T e^{-i\omega t} e^{itA} f \,dt \to 0 \]
for all $f$ belonging to a dense subspace of $\Ran(P_c)$. In view of 
the preceding, 
it therefore suffices to show that for all $|\omega|\ge\mu$, 
\begin{equation}
\label{eq:dense}
\overline{(A-\omega)\Ran(P_c)}=\Ran(P_c).
\end{equation}
By the $A$-invariance of $\Ran(P_c)$ it is clear that 
$\overline{(A-\omega)\Ran(P_c)}\subset\Ran(P_c)$.
By Lemma~\ref{lem:split} we can write a direct sum decomposition 
$\Hil=L'+\tilde L$ where
$L'=\sum_{j=0}^M L_j$ and $\tilde L$ is the space from above with 
finite codimension. 
Now  $\ker(A^*+\omega)=\{0\}$ implies that 
\begin{equation}
\label{eq:s1} 
\Hil = \overline{(A-\omega)\Hil} =   \overline{(A-\omega)\tilde{L}} 
+ (A-\omega){L'}.
\end{equation}
Note that by $A$-invariance of the splitting both summands here lie 
again in the spaces $\tilde L$
and $L'$, respectively. In particular, the sum on the right-hand side 
is again direct. 
Since $\ker(A-\omega)=\{0\}$, one has $\dim((A-\omega){L'})=\dim 
(L')<\infty$ so that 
$(A-\omega){L'}=L'$. By the directness of the sum in~\eqref{eq:s1} 
therefore 
\[ \overline{(A-\omega)\tilde{L}} = \tilde L,\]
as desired.
\end{proof}

\noindent 
The last of the abstract statements is a bound on the resolvents.

\begin{lemma}
\label{lem:resol} For all $z\in\Compl$ with $\Re z>0$ one has
\begin{equation}
\label{eq:Aresol}
 \|(iA-z)^{-1}P_c\| \le \frac{C}{\Re z},
\end{equation}
for some constant $C$.
\end{lemma}
\begin{proof} Since for $\Re z>0$,
\[ (iA-z)^{-1}P_c = - \int_0^\infty e^{-tz} e^{itA} P_c\,dt, \]
\eqref{eq:Aresol} is an immediate consequence of~\eqref{eq:stab} with
$C=\sup_{t\in\R}\|e^{itA} P_c\|$. 
\end{proof}

Following Rauch~\cite{R}, we now show how to continue the resolvents 
$(iA-z)^{-1}=-i(A+iz)^{-1}$ meromorphically across 
the spectrum as bounded operators in an exponentially weighted $L^2$ 
space. Set $(E_\eps f)(x):= e^{-\eps\rho(x)}f(x)$ where $\rho$ 
is smooth, $\rho(x)=|x|$ for large~$x$ and $\eps>0$ will be some 
small number (compared to $\eps_0$ in \eqref{eq:Vdec} above). Our 
goal now is to
prove that $E_\eps(iA-z)^{-1}E_\eps$ has a meromorphic continuation to 
the region which lies to the right
of the curve shown below. From now on, we set $\mu=1$, so that 
$H=-\Lapl+1$. Moreover, we denote the curve
depicted in the following figure by~$\Gamma_\eps$ and the open region 
to the right of it 
by~$\Reg_{\Gamma_\eps}$. 

\centerline{\hbox{\vbox{
\epsfxsize= 6.0 truecm
\epsfysize= 5.0truecm
\epsfbox{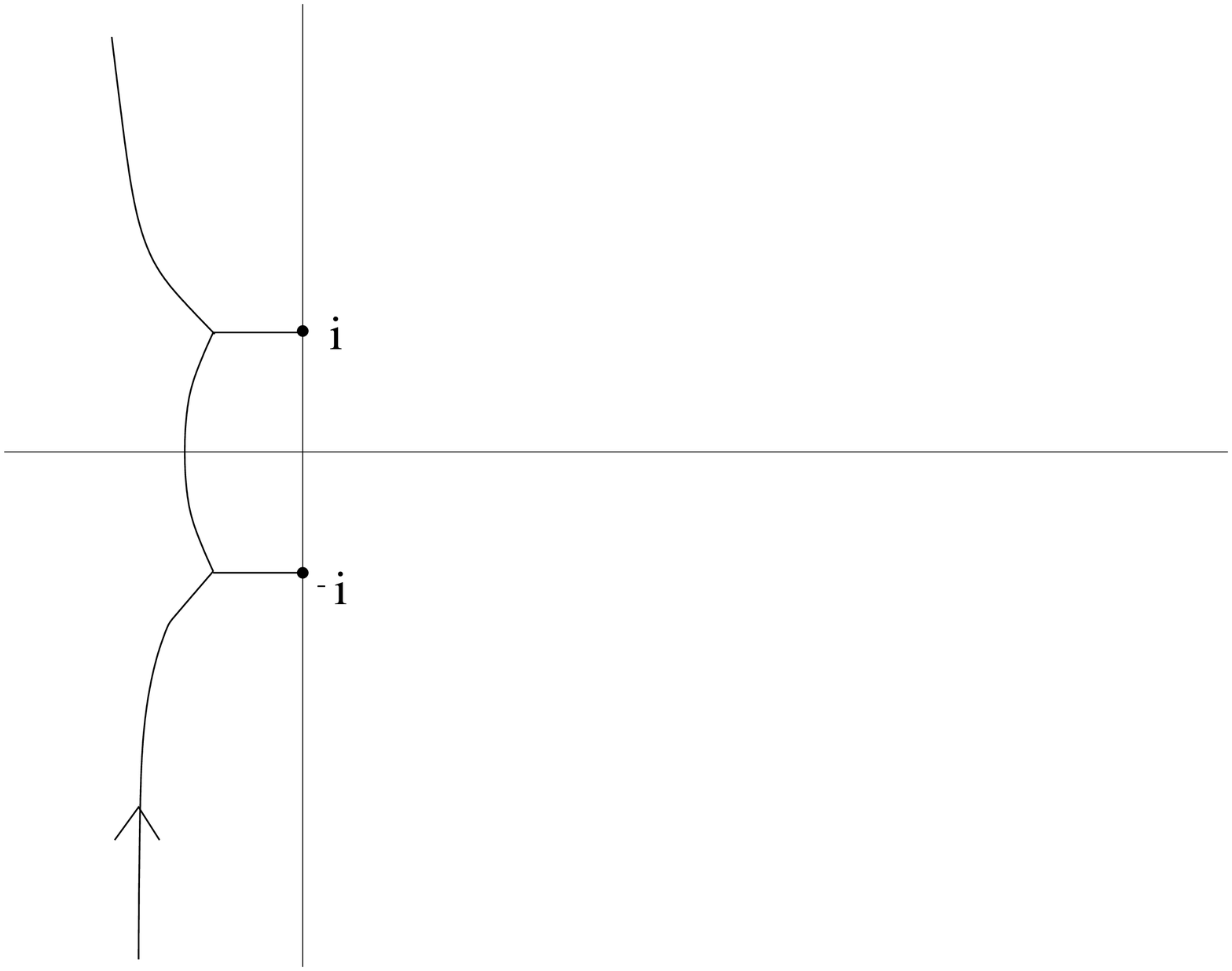}}}}

\noindent We start with the resolvent of~$B$. 
By~\eqref{eq:Bdef} one has
\begin{equation}
\label{eq:Bresol} 
(B+iz)^{-1} = \bm (H+iz)^{-1} & 0 \\ 
 0 & -(H-iz)^{-1} \endm.
\end{equation}
The resolvents
\begin{equation}
\label{eq:free_resol} 
(H\pm iz)^{-1}(x,y)=(-\Lapl+1\pm iz)^{-1}(x,y) = 
\frac{e^{-\sqrt{2(1\pm iz)}\,|x-y|}}{2\pi|x-y|}
\end{equation}
have singularities at the points $\pm i$, respectively. We will 
resolve these singularities 
by means of the transformations $U_1:\zeta \mapsto z=i-i\zeta^2$, and 
$U_2:\zeta\mapsto -i+i\zeta^2$. 
For the domain of $U_1$ we will choose either of the regions
\bea
\Reg_1=\Reg_1(\eps) &:=& \Big\{ \zeta\in\Compl\:|\: 
\Re\zeta>-\eps,\;\Im\zeta>-\eps,\;-\frac{\pi}{4}\le\Arg\,\zeta\le
\frac{3\pi}{4},\;|\zeta|\le r\Big\} \label{eq:reg1} \\
\Reg_2=\Reg_2(\eps) &:=& \Big\{ \zeta\in\Compl\:|\: 
\Re\zeta>-\eps,\;\frac{\pi}{4}\le\Arg\,\zeta\le \frac{3\pi}{4}\Big\}. 
\label{eq:reg2}
\eea
Here $\eps>0$ is a small parameter, and $r$ will be chosen to be on 
the order of~$1$. 
These two regions are shown in the following figures. Firstly 
region~$\Reg_1$, with $\zeta$ on the
left and~$z=U_1(\zeta)$ on the right,

\centerline{\hbox{\vbox{
\epsfxsize= 5.0 truecm
\epsfysize= 4.0truecm
\epsfbox{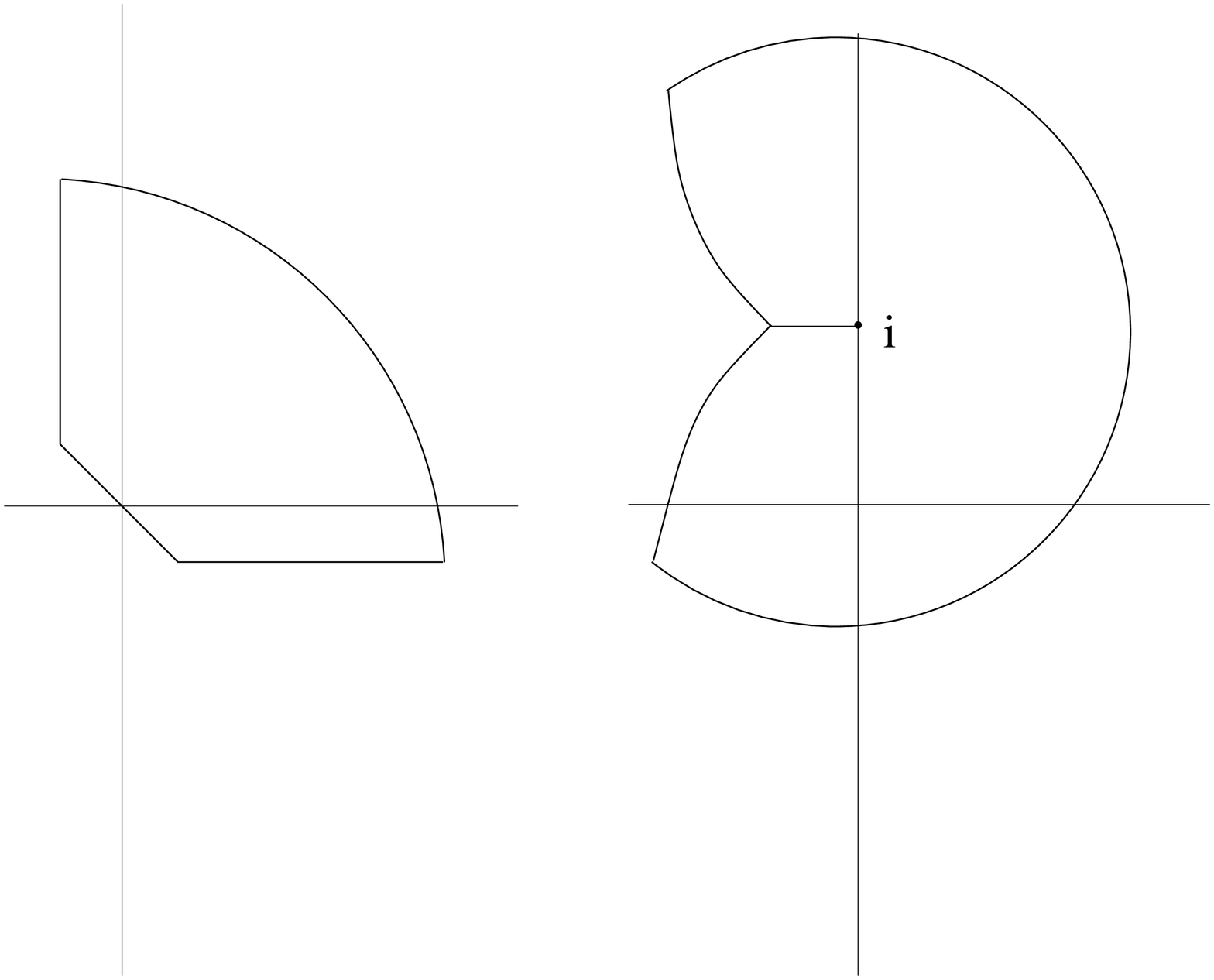}}}}

\noindent and secondly region~$\Reg_2$, with $\zeta$ on the
left and~$z=U_1(\zeta)$ on the right:

\centerline{\hbox{\vbox{
\epsfxsize= 5.0 truecm
\epsfysize= 4.0truecm
\epsfbox{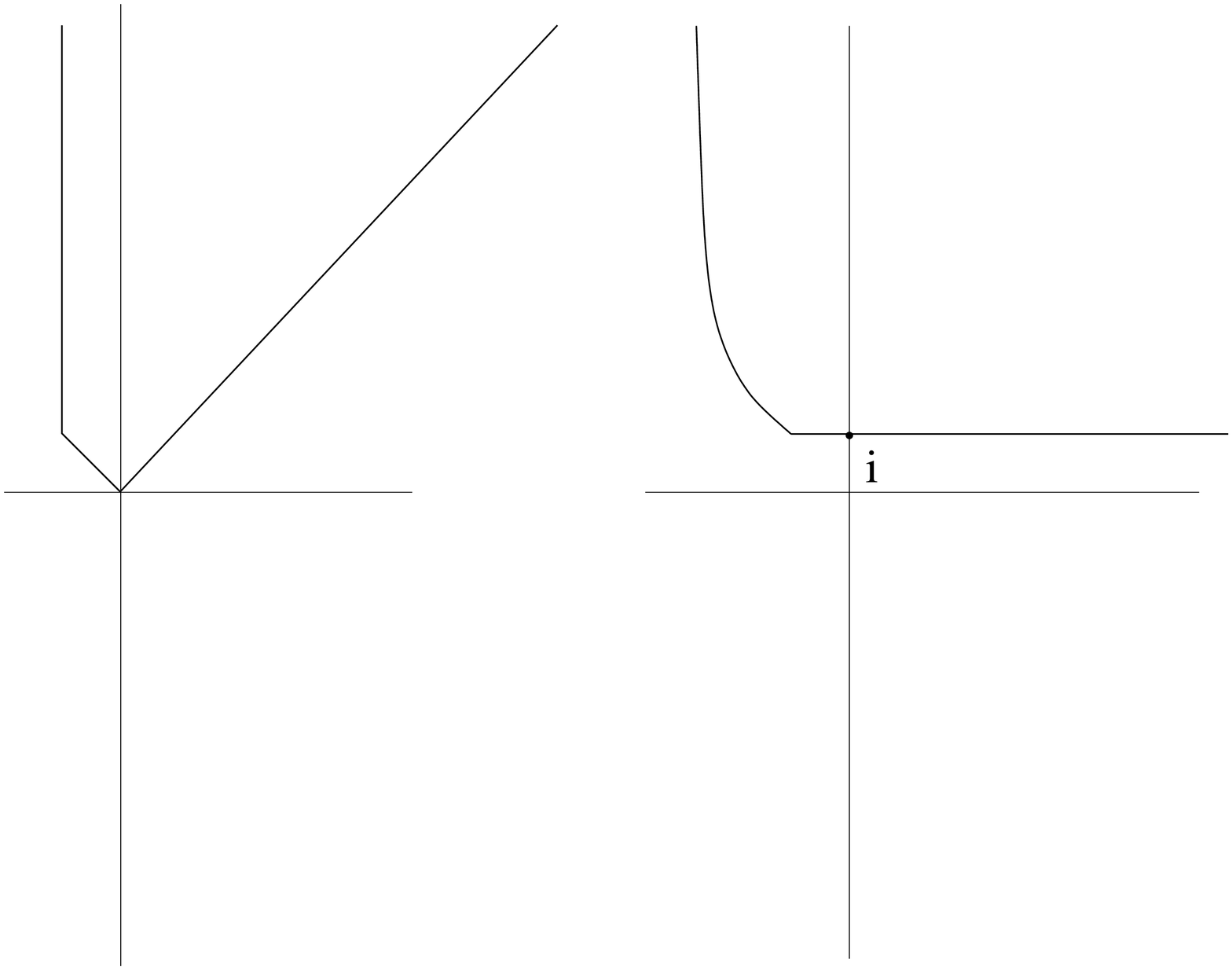}}}}

\noindent The curved (but not circular) pieces of the boundaries that 
are shown in the $z$-plane are
parabolic arcs which can be of course written down exactly. 
For the map $U_2$ we choose the domains $\Reg_1^*$ and $\Reg_2^*$ 
which are the reflections
of $\Reg_1$ and $\Reg_2$, respectively, across the real axis (i.e., 
conjugation; note that 
$U_2(\zeta)=\overline{U_1(\bar\zeta)}$).
One now defines $\Reg_{\Gamma_\eps}\cap\{z\in\Compl\:|\: \Re z\le 
1/2\}$ by means of
four separate regions, namely by $U_1(\Reg_1)$, $U_1(\Reg_2)$, 
$U_2(\Reg_1^*)$, $U_2(\Reg_2^*)$ with
$r$ sufficiently large, say $r=\sqrt{5}/2$ (by taking $\eps$ small we 
can still ensure that we move only
a very short distance into the negative half plane $\Re z<0$, at 
least on compact sets).
This choice of~$r$ is made to ensure that the two disks intersect the 
real axis without coming
close to each others centers.

\begin{lemma}
\label{lem:free_cont} 
Let $0<\eps$ be small. Then the weighted 
resolvent $E_\eps (iB-z)^{-1} E_\eps$,
originally defined on $\Re z>0$, admits an analytic continuation to 
the region $\Reg_{\Gamma_{\eps/4}}$
with values in the compact operators on~$\Hil$.
This continuation, which we denote by~$F_\eps(z)$, satisfies the 
estimate
$ \|F_\eps(z)\| \le C_\eps (1+|z|)^{-\half}$
for all $z\in\Reg_{\Gamma_{\eps/4}}$. Furthermore, $F_\eps(z)$ is
analytic on the union of the regions
$z=i-i\zeta^2$ and $z=-i+i\zeta^2$
for $|\zeta|<\frac{\eps}{4}$. 
\end{lemma}
\begin{proof} 
It was shown by Rauch that (recall that $H=-\Lapl +1$)
\[\zeta\mapsto E_\eps(H+i(i-i\zeta^2))^{-1}E_\eps = E_\eps 
(-\Lapl+\zeta^2)^{-1} E_\eps\]
has an analytic extension to the domain~$\Re \zeta>-\eps/4$ with a 
norm bound of $(1+|\zeta|)^{-1}$ 
(see equation~(2.12) in~\cite{R}). 
Similarly for $\zeta\mapsto 
E_\eps(H-i(-i+i\zeta^2))^{-1}E_\eps=E_\eps (-\Lapl+\zeta^2)^{-1} 
E_\eps$. 
This shows that both $E_\eps(H\pm iz)^{-1}E_\eps$ 
have analytic continuations to the region $\Reg_{\Gamma_{\eps/4}}$
satisfying the norm bound $(1+|z|)^{-\half}$ as $|z|\to\infty$ 
in~$\Reg_{\Gamma_{\eps/4}}$. The lemma follows.
\end{proof}

\noindent 
We can now formulate our meromorphic continuation result for the 
perturbed resolvent. 

\begin{cor}
\label{cor:pert_cont} Let $0<\eps$ be small depending on~$\eps_0$ 
from Definition~\ref{def:spec_ass}.
 Then for some $\delta>0$, 
the weighted resolvents $E_\eps (iA-z)^{-1} P_c E_\eps$,
originally defined on $\Re z>0$, admit a meromorphic continuation to 
the region $\Reg_{\Gamma_{\delta}}$
with values in the compact operators on~$\Hil$.
This continuation, which we denote by~$G_\eps(z)$, has only finitely 
many poles in~$\Reg_{\Gamma_\delta}$
all of which lie in~$\Re z\le 0$. 
Moreover, it has continuous limits on the two horizontal parts 
of~$\Gamma_\delta\setminus\{\pm i\}$ 
 from above and below. Finally, for large $|z|$ there are no poles 
and one has the estimate
\begin{equation}
\label{eq:Gdec}
 \|G_\eps(z)\| \le \frac{C_\eps}{\sqrt{|z|}}
\end{equation}
for all large $z\in\Reg_{\Gamma_{\delta}}$.
\end{cor}
\begin{proof} 
Since $P_c=Id-P_b$ and one has the explicit 
representation~\eqref{eq:Pbexpl} with 
exponentially decaying~$\phi_j$, it suffices to prove the same 
statement for 
$E_\eps (iA-z)^{-1} E_\eps=-iE_\eps(A+iz)^{-1}E_\eps$. By the resolvent identity,
\bea
 E_\eps (A+iz)^{-1} E_\eps &=& E_\eps(B+iz)^{-1} E_\eps - E_\eps 
(B+iz)^{-1}E_\eps\, E_\eps^{-1}VE_\eps^{-1}
E_\eps (A+iz)^{-1} E_\eps \nn \\
E_\eps (A+iz)^{-1} E_\eps &=& \Bigl(1+E_\eps (B+iz)^{-1} E_\eps \, 
E_\eps^{-2}V\Bigr)^{-1}\;E_\eps (B+iz)^{-1} E_\eps. \nn
\eea
Lemma~\ref{lem:free_cont} and the analytic Fredholm alternative 
imply that~$E_\eps (A+iz)^{-1} E_\eps$ has a meromorphic extension 
to~$\Reg_{\Gamma_{\eps/4}}$
which is analytic if~$|z|$ is large. 
Moreover, 
\begin{equation}
\label{eq:Geps}
 G_\eps(z) = \Bigl(1+iF_\eps(z) \, 
E_\eps^{-2}V\Bigr)^{-1}\;F_\eps(z)
\end{equation}
for all $z\in\Reg_{\Gamma_{\eps/4}}$ at which the inverse is defined. 
Hence~\eqref{eq:Gdec} follows from the previous lemma. 
The poles can only accumulate on the boundary~$\Gamma_{\eps/4}$ 
of~$\Reg_{\Gamma_{\eps/4}}$. If~$\delta<\frac{\eps}{4}$, then the 
poles of~$G_\eps$ can therefore only accumulate on the two horizontal 
pieces of~$\Gamma_\delta$. We claim that reducing~$\delta$ even 
further, one can eliminate any accumulation
on these horizontal pieces as well. To see this, ``resolve'' the 
singularities at the points~$\pm i$
by means of the maps $z=U_1(\zeta)=i-i\zeta^2$ and 
$z=U_2(\zeta)=-i+i\zeta^2$, respectively. Clearly, 
\begin{equation}
\label{eq:Geps_2}
 G_\eps(U_1(\zeta)) = \Bigl(1+F_\eps(U_1(\zeta)) \, 
E_\eps^{-2}V\Bigr)^{-1}\;F_\eps(U_1(\zeta)) 
\end{equation}
as long as $U_1(\zeta)\subset \Reg_{\Gamma_{\eps/4}}$. In particular, 
this holds for 
all~$\zeta\in\Reg(\eps/4)$. In view of the final statement of 
Lemma~\ref{lem:free_cont}, one
can apply the analytic Fredholm alternative to~\eqref{eq:Geps_2} with 
$|\zeta|<\frac{\eps}{4}$.
Thus $G_\eps(U_1(\zeta))$ is a meromorphic function 
for~$|\zeta|<\frac{\eps}{4}$. It follows that
for $\delta>0$ small, $G_\eps(U_1(\zeta))$ is analytic 
on~$|\zeta|<\delta$ with the possible
exception of a pole at $\zeta=0$. In the~$z$ variable this means that 
$G_\eps(z)$ has continuous
limits (from above and below) on the horizontal piece 
of~$\Gamma_\delta$ emanating from~$+i$, with a
possible pole at $z=i$. The same statements also hold at~$-i$. 
\end{proof}

We are now ready to prove the local $L^2$-decay.

\begin{prop}
\label{prop:localL2}
Let $A$ be as in Definition~\ref{def:spec_ass} together 
with~\eqref{eq:stab}. Then
for $\eps>0$ small one has
\begin{equation}
\label{eq:rauch}
\|E_\eps e^{itA}P_c E_\eps\| \les \la t \ra^{-\frac32}.
\end{equation}
for all $t$. 
\end{prop}
\begin{proof} Suppose $t>0$. By the inversion formula for the Laplace transform (to 
be justified below),  
\begin{equation}
\nn
E_\eps (e^{itA}-I)P_c E_\eps = \frac{-1}{2\pi i} 
\int_{a-i\infty}^{a+i\infty} e^{t\tau} 
\Bigl( E_\eps (iA-\tau)^{-1}P_c\, E_\eps  + E_\eps P_c E_\eps 
\tau^{-1}\Bigr) \, d\tau 
\end{equation}
for any $a>0$. The subtraction of the identity serves to make the 
integral absolutely convergent. 
Indeed, $(iA-\tau)(iA-\tau)^{-1}=Id$ implies that
\[ E_\eps(iA-\tau)^{-1}P_cE_\eps+\frac{1}{\tau}E_\eps P_c E_\eps = 
\frac{1}{\tau}E_\eps(iA-\tau)^{-1}P_cE_\eps E_\eps^{-1}iA E_\eps. \]
By means of a commutator calculation one checks that $ E_\eps^{-1}A 
E_\eps$ is a bounded
operator $W^{2,2}\times W^{2,2}\to\Hil$. Hence, by~\eqref{eq:Gdec}, 
the left-hand side decays like~$|\tau|^{-\frac32}$. We now deform the 
contour to $\Gamma_\delta$ as in Corollary~\ref{cor:pert_cont}. This 
gives
\bea
\nn
E_\eps (e^{itA}-I)P_c E_\eps &=& \frac{-1}{2\pi i} 
\int_{\Gamma_\delta} e^{t\tau} 
\Bigl( G_\eps(\tau)  + E_\eps P_c E_\eps \tau^{-1}\Bigr) \, d\tau \\
&& \qquad - \sum_{\tau_j\text{\ poles}}\Res\Bigl(e^{t\tau} 
\Bigl( G_\eps(\tau)  + E_\eps P_c E_\eps \tau^{-1}\Bigr); 
\tau=\tau_j\Bigr).  \nn
\eea
This can be simplified to
\begin{equation}
\label{eq:residues}
E_\eps e^{itA}P_c E_\eps = \frac{-1}{2\pi i} \int_{\Gamma_\delta} 
e^{t\tau} \, G_\eps(\tau) \, d\tau 
 - \sum_{\tau_j\text{\ poles}} \Res\Bigl(e^{t\tau} 
 G_\eps(\tau); \tau=\tau_j\Bigr). 
\end{equation}
The integral over $\Gamma_\delta$ converges by~\eqref{eq:Gdec}. 
Clearly,
\[ \sum_{\tau_j\text{\ poles}} \Res\Bigl(e^{t\tau} G_\eps(\tau); 
\tau=\tau_j\Bigr) = 
\sum_{\tau_j\text{\ poles}} p_j(t) e^{t\tau_j},\]
where $p_j(t)$ is a polynomial with coefficients given by compact 
operators 
of degree equal to the order of the pole minus one.
As far as the integral over $\Gamma_\delta$ is concerned, only the 
horizontal ``thermometers'' can
contribute much. ``Thermometer'' here refers to the fact that on the 
horizontal pieces 
one first moves to the right
on the lower edge, then around a loop encircling $\pm i$, and then to 
the left on the upper edge.
It is clear that the parabolic arcs give an exponentially small 
contribution. 
Denote the two thermometers around $i$ and $-i$ by $\T_1$ and $\T_2$, 
respectively. To determine
their contributions, we need to again use the $\zeta$-variables. 
Indeed, setting $z=U_1(\zeta)=i-i\zeta^2$
 we showed above that $G_\eps(z)$ is analytic for  $z=i-i\zeta^2$
and small $\zeta$ up to possibly a pole at~$0$.
Writing the Laurent expansion  $G_\eps(i-i\zeta^2) = 
\sum_{j=-N}^\infty B_j \zeta^{j}$, one
deduces from Lemma~\ref{lem:resol} that $N\le 2$. Also, our 
assumption on the absence of resonances
means that $B_{-1}=0$. Therefore, 
\begin{equation}
\label{eq:laurent}
 \frac{1}{2\pi i} \int_{\T_1}e^{t\tau}\, G_\eps(\tau)\,d\tau = 
\frac{1}{2\pi i} \int_{\T_1} e^{t\tau}\,\frac{B_{-2}}{\tau-i}\,d\tau
+ \frac{1}{2\pi i} \int_{\T_1} e^{t\tau}\,C(\tau)\,d\tau,
\end{equation}
where $\|C(\tau)\|\les |\tau-i|^{\half}$ so that the contribution of 
the second integral is 
$O(|t|^{-\frac32})$. The first is equal to the compact operator 
$e^{it}\,B_{-2}$. 
The conclusion is that~\eqref{eq:residues} now takes the form of the 
estimate
\begin{equation}
\label{eq:crux}
\Big\| E_\eps e^{itA}P_c E_\eps + \sum_{j} p_j(t) 
e^{it\eta_j} \Big\| \les |t|^{-\frac32}
\end{equation}
where we have only retained the purely imaginary poles 
$\tau_j=e^{i\eta_j}$ with distinct~$\eta_j$. 
In view of~\eqref{eq:stab}
it follows that for large~$t$, 
\[ \sup_{t\in\R}\Big\|\sum_{j} p_j(t) e^{it\eta_j} \Big\| \les 1 \]
which implies that necessarily each $p_j=\const$. Finally, the 
ergodic average lemma, 
Lemma~\ref{lem:ergodic}, 
shows that $p_j=0$. 
Hence we are left with
\[ \Big\| E_\eps e^{itA}P_c E_\eps \Big\| \les |t|^{-\frac32},\]
for large~$t$ as desired.
\end{proof}

\begin{remark}
\label{rem:reson}
According to Corollary~\ref{cor:pert_cont} (and its proof), 
$G_\eps(i-i\zeta^2)$
and $G_\eps(-i+i\zeta^2)$ can at most have poles of second order at 
$\zeta=0$. 
It is well-known that in generic situation there will be no 
singularity at the origin.
In fact, a second order pole means that $\pm i$ are eigenvalues, 
whereas a first
order pole is what one calls a resonance. For more on this, see for 
example~\cite{R}. 
\end{remark}

We now pass from local $L^2$-decay to global decay.

\begin{proof}[Proof of Theorem~\ref{thm:vec_dec}] We start with the first
statement in the theorem. By Duhamel's 
formula, and with $A,B,V$ as in~\eqref{eq:Bdef}, 
\bea
 e^{itA}P_c &=& e^{itB}P_c -i \int_0^t e^{i(t-s)B}\, V e^{isA}P_c \,ds 
           = e^{itB}P_c - i \int_0^t e^{i(t-s)B}\, V P_c e^{isA} \,ds 
\label{eq:duhp1} \\
 P_c e^{isB} &=& P_c e^{isA} + i \int_0^s P_c e^{i(s-\tau)A}\, V e^{i\tau 
B}\, d\tau 
 	    = P_c e^{isA} + i \int_0^s  e^{i(s-\tau)A}\, P_c V e^{i\tau B}\, d\tau 
\label{eq:duhp2}
\eea
Inserting \eqref{eq:duhp2} into \eqref{eq:duhp1} yields 
\begin{equation}
e^{itA}P_c = e^{itB}P_c - i \int_0^t e^{i(t-s)B}\, V P_c e^{isB} \,ds 
+ \int_0^t \int_0^s  e^{i(t-s)B}\, V e^{i(s-\tau)A}\, P_c V e^{i\tau B}\, 
d\tau ds. \nn
\end{equation} 
Therefore, applying Proposition~\ref{prop:localL2} to the middle term 
in the double integral, yields
\bea
 \|e^{itA} P_c \vpsi_0\|_{L^2+L^\infty} 
 &\les& \la t\ra^{-\frac32} \|P_c\|_{L^1\cap L^2\to L^1\cap 
L^2}\;\|\vpsi_0\|_{L^1\cap L^2} \nn \\
&&   + \int_0^t \la t-s \ra^{-\frac32}  \la s\ra^{-\frac32}\,ds \;\|V 
P_c\|_{L^2+L^\infty\to L^1\cap L^2}\;\|\vpsi_0\|_{L^1\cap L^2} \nn \\
&&  + \int_0^t \int_0^s  \la t-s \ra^{-\frac32} \la 
s-\tau\ra^{-\frac32} \la \tau \ra^{-\frac32}\,d\tau
ds\;\|E_\eps^{-1}V\|^2_{L^2\to L^1\cap L^2}\; \|\vpsi_0\|_{L^1\cap L^2} \nn 
\\
&\les&  \la t\ra^{-\frac32} \|\vpsi_0\|_{L^1\cap L^2}. \label{eq:uff}
\eea
To pass to \eqref{eq:uff} one uses that 
\[ \|V P_c\|_{L^2+L^\infty\to L^1\cap L^2}<\infty, \quad 
 \|P_c\|_{L^1\cap L^2\to L^1\cap L^2} <\infty, \quad
\|VE_\eps^{-1}\|_{L^2\to L^1\cap L^2}<\infty,
\] 
which follow from~\eqref{eq:Pbexpl} and the assumption
that $V,\phi_j, \psi_i$ are exponentially decaying.

\noindent To remove $L^2$ we use the same device from~\cite{JSS}
that already appeared in Section~\ref{sec:infty} for 
the same purpose. However, some care is needed for
the case of systems since the analogue of~\eqref{eq:canc}, viz.
\[ \sup_{1\le p\le\infty}\|e^{itB}Ve^{-itB}\|_{p\to p}<\infty \]
{\em does not} hold now. Indeed, by~\eqref{eq:Bdef}, and with $H=-\Lapl+\mu$, 
\begin{equation}
\label{eq:mat_prod}
 e^{itB}Ve^{-itB} = 
\bm e^{itH}Ue^{-itH} & -e^{itH}We^{itH} \\
    e^{-itH}We^{-itH} & -e^{-itH} Ue^{itH}
\endm, 
\end{equation}
and the off-diagonal terms are not bounded on $L^p$ if $p\ne2$, $W\ne0$.
As in~\eqref{eq:duha},
\begin{align}
& e^{itA}\vpsi_0 = e^{itB}\vpsi_{0} - i \int_{0}^{t} 
e^{i(t-s)B} V e^{isA}\vpsi_0\, ds \nn\\
&= e^{it B}\vpsi_{0} - i \int_{0}^{t} 
e^{i(t-s)B} V e^{isB}\vpsi_{0}\, ds
- \int_{0}^{t}\int_{0}^{s}e^{i(t-s)B} V
e^{i(s-\tau)B} V e^{i\tau A}\vpsi_0\, d\tau\, ds.
\label{eq:duha_sys}
\end{align}
In what follows it suffices to assume $U=0$, since the 
diagonal terms in~\eqref{eq:mat_prod} have the same cancellation
as in the scalar case, see~\eqref{eq:canc}. 
Consider one of the off-diagonal entries of 
the second term in~\eqref{eq:duha_sys} for times $t\ge\half$. 
Thus,
\bea
&& \Big\|\int_0^t e^{i(t-s)H} W e^{-isH}\vpsi_0 \, ds\Big\|_\infty \le
\Big\|\int_0^{\frac14} e^{i(t-s)H} W e^{-isH}\vpsi_0 \,ds \Big\|_\infty 
+ \Big\|\int_{\frac14}^{t-\frac14} e^{i(t-s)H} W e^{-isH}\vpsi_0 \,ds 
\Big\|_\infty  \nn \\
&& \mbox{\hspace{2.0in}} 
+ \Big\|\int_{t-\frac14}^t e^{i(t-s)H} W e^{-isH}\vpsi_0 \,ds 
\Big\|_\infty  \nn \\
&& \les \int_0^{\frac14} |t-2s|^{-\frac32}\, \|e^{isH}W 
e^{-isH}\|_{1\to1}\|\vpsi_0 \|_1\,ds 
+ \int_{\frac14}^{t-\frac14} |t-s|^{-\frac32}\,\|W\|_1\,
|s|^{-\frac32}\|\vpsi_0\|_1 \,ds \nn \\
&& \mbox{\hspace{2.0in}} 
+ \int_{t-\frac14}^t \|e^{i(t-s)H} W e^{-i(t-s)H}\|_{\infty\to\infty}
\|e^{i(t-2s)H}\vpsi_0\|_\infty \,ds 
\nn \\
&&\qquad \les (\|W\|_1+\|\hat W\|_1)\,|t|^{-\frac32} \|\vpsi_0\|_1, 
\label{eq:erst}
\eea
where we have employed the cancellation property~\eqref{eq:canc} twice.
For small times $0<t$ one obtains similarly
\bea
&& \Big\|\int_0^t e^{i(t-s)H} W e^{-isH}\vpsi_0 \, ds\Big\|_\infty \le
\Big\|\int_0^{\half t-t^{2}}  e^{i(t-s)H} 
W e^{-isH}\vpsi_0 \,ds \Big\|_\infty  \nn \\
&& \qquad\qquad + \Big\|\int_{|\frac{t}{2}-s|<t^{2}} \!\!\! e^{i(t-s)H} 
W e^{-isH}\vpsi_0 \,ds \Big\|_\infty 
+ \Big\|\int_{\half t+t^{2}}^t e^{i(t-s)H} W e^{-isH}
\vpsi_0 \,ds \Big\|_\infty  \nn \\
&& \les \int_0^{\half t-t^{2}} |t-2s|^{-\frac32}\, 
\|e^{isH}W e^{-isH}\|_{1\to1}\|\vpsi_0 \|_1\,ds 
+ \int_{|\frac{t}{2}-s|<t^{2}} |t-s|^{-\frac32}\,
\|W\|_1\,|s|^{-\frac32}
\|\vpsi_0\|_1 \,ds \nn \\
&& \mbox{\hspace{1.7in}} 
+ \int_{\half t+t^{2}}^t \|e^{i(t-s)H} W 
e^{-i(t-s)H}\|_{\infty\to\infty}
\|e^{i(t-2s)H}\vpsi_0\|_\infty \,ds \nn \\
&&  \label{eq:zweit} \qquad \les (\|W\|_1+\|\hat W\|_1)\,
|t|^{-1} \|\vpsi_0\|_1.
\eea
We now deal with the third term in~\eqref{eq:duha_sys}, 
and we select
one representative off-diagonal term of the form
\begin{equation}
\label{eq:teil_2}
 \int_0^t\int_0^s e^{i(t-s)H} W e^{i(\tau-s)H} 
U\psi_1(\tau)\,d\tau ds,
\end{equation}
where $\psi_1(\tau)$ is the first component of the 
solution $\vpsi(\tau):= e^{i\tau A}\vpsi_0$. For $t\ge 4$,
introducing the new variable $\tau':= 2s-\tau$, we rewrite
\eqref{eq:teil_2} in the form
\bea
& &\int_0^t\int_s^{2s} e^{i(t-s)H} W e^{i(s-\tau')H} 
U\psi_1(2s-\tau')\,d\tau' ds = \nn\\
& &\bigg (\int_0^{t-1}\int_s^{2s} + \int_{t-1}^t\int_s^{2t-s} + 
\int_{t-1}^{t}\int_{2t-s}^{2s}  \bigg ) e^{i(t-s)H} W e^{i(s-\tau')H} 
U\psi_1(2s-\tau')\,d\tau' ds \label{eq:tail1}
\eea
The estimate for the first term in \eqref{eq:tail1} is 
straightforward
\begin{align}
\|\int_0^{t-1}\int_s^{2s} &e^{i(t-s)H} W e^{i(s-\tau')H} 
U\psi_1(2s-\tau')\,d\tau' ds\|_{\infty} \les\nn\\ 
&\int_0^{t-1}\int_s^{2s} \la t-s\ra^{-\frac 32} \|W\|_{1\cap 2}
\la s-\tau'\ra^{-\frac 32} \|U\|_{1\cap 2} \la 2s-\tau'\ra^{-\frac 32}
\|\vpsi_{0}\|_{1\cap 2}\,d\tau' ds \les\nn\\ &\la t\ra^{-\frac 32} 
\|\vpsi_{0}\|_{1\cap 2}\label{eq:tail2}
\end{align}
The second term in \eqref{eq:tail1} is dealt with in the following 
manner:
\begin{align}
\|\int_{t-1}^{t}\int_s^{2t-s}& e^{i(t-s)H} W e^{i(s-\tau')H} 
U\psi_1(2s-\tau')\,d\tau' ds\|_{\infty} \les\nn\\ 
&\int_{t-1}^{t}\int_s^{2t-s} |t-s|^{-\frac 32} \|W\|_{2} 
\|U\|_{2\cap\infty}
\la 2s-\tau'\ra^{-\frac 32}\|\vpsi_{0}\|_{1\cap 2}\,d\tau' ds \nn\\ &\les 
\la t\ra^{-\frac 32}
\|W\|_{2} \|U\|_{2\cap\infty}\int_{t-1}^{t}\int_s^{2t-s} |t-s|^{-\frac 
32}\|\vpsi_{0}\|_{1\cap 2}\,d\tau' ds \les\nn\\ &
\la t\ra^{-\frac 32}\int_{t-1}^{t} |t-s|^{-\frac 
12}\|\vpsi_{0}\|_{1\cap 2}\, ds 
\les t^{-\frac 32}\|\vpsi_{0}\|_{1\cap 2}\label{eq:tail3}
\end{align}
The last term in \eqref{eq:tail1} requires further splitting 
\bea
& &\int_{t-1}^{t}\int_{2t-s}^{2s} e^{i(t-s)H} W e^{i(s-\tau')H} 
U\psi_1(2s-\tau')\,d\tau' ds = \nn\\
& &\bigg (\int_{t-1}^{t}\int_{2t-s}^{t+\frac s2} + 
\int_{t-1}^{t}\int_{t+\frac s2}^{2s} \bigg ) 
e^{i(t-s)H} W e^{i(s-\tau')H} 
U\psi_1(2s-\tau')\,d\tau' ds\label{eq:tail4}
\eea
Observe that in the region $(s,\tau') \in [t-1,t]\times [2t-s, 
t+\frac s2]$ we have 
$2s-\tau'\ge \frac{t-3}2$. On the other hand, in the region 
$(s,\tau') \in [t-1,t]\times [t+\frac s2, 2s]$ we have 
$|t-\tau'|\ge \frac t2 - \frac 12$. 
Therefore,
\begin{align}
\|\int_{t-1}^{t}\int_{2t-s}^{t+\frac s2} &
e^{i(t-s)H} W e^{i(s-\tau')H} 
U\psi_1(2s-\tau')\,d\tau' ds \|_{\infty}=\nn\\
&\int_{t-1}^{t}\int_{2t-s}^{t+\frac s2}
\|e^{i(t-s)H} W  e^{-i(t-s)H}\|_{\infty\to\infty}
\|e^{i(t-\tau')H} 
U\psi_1(2s-\tau')\,d\tau' ds \|_{\infty}\,d\tau' ds \les\nn\\
&\|\hat W\|_{1} \|U\|_{1\cap 2}
\int_{t-1}^{t}\int_{2t-s}^{t+\frac s2} |t-\tau'|^{-\frac 32}
\la 2s -\tau'\ra^{-\frac 32}\|\vpsi_{0}\|_{1\cap 2}\, d\tau' ds \les\nn\\ 
&\la t\ra^{-\frac 32} \int_{t-1}^{t}\int_{2t-s}^{t+\frac s2} |t-\tau'|^{-\frac 32}
\|\vpsi_{0}\|_{1\cap 2}\, d\tau' ds 
\les\nn\\ &\la t\ra^{-\frac 32} \int_{t-1}^{t} \bigg ( |t-s|^{-\frac 
12} - |s/2|^{-\frac 12}\bigg )\|\vpsi_{0}\|_{1\cap 2}\, ds 
\les t^{-\frac 32} \|\vpsi_{0}\|_{1\cap 2}\label{eq:tail5}
\end{align}
Similarly,
\begin{align}
\|\int_{t-1}^{t}\int_{t+\frac s2}^{2s}&
e^{i(t-s)H} W e^{i(s-\tau')H} 
U\psi_1(2s-\tau')\,d\tau' ds \|_{\infty}\les\nn\\
&\|\hat W\|_{1} \|U\|_{1\cap 2} 
\int_{t-1}^{t}\int_{t+\frac s2}^{2s} |t-\tau'|^{-\frac 32} 
\la 2s-\tau'\ra^{-\frac 32} \|\vpsi_{0}\|_{1\cap 2}\, d\tau' ds \les\nn \\ 
&t^{-\frac 32} \int_{t-1}^{t}\int_{t+\frac s2}^{2s} 
\la 2s-\tau'\ra^{-\frac 32} \|\vpsi_{0}\|_{1\cap 2}\, d\tau' ds \les 
t^{-\frac 32} \|\vpsi_{0}\|_{1\cap 2}\label{eq:tail6}
\end{align}

For the values of $t<4$ we decompose \eqref{eq:teil_2}, with the 
variable $\tau'=2s-\tau$, as follows.
\bea
& &\int_0^t\int_s^{2s} e^{i(t-s)H} W e^{i(s-\tau')H} 
U\psi_1(2s-\tau')\,d\tau' ds = \nn\\
& &\bigg (\int_0^{\frac {3t}4}\int_s^{2s} + \int_{\frac {3t}4}^t\int_s^{2t-s} + 
\int_{\frac {3t}4}^{t}\int_{2t-s}^{2s}  \bigg ) e^{i(t-s)H} W e^{i(s-\tau')H} 
U\psi_1(2s-\tau')\,d\tau' ds \label{eq:tail7}
\eea
Note that $(2t-s) \le 2s$ provided that $s\ge \frac {2t}3$. In 
particular, this holds for $s\ge \frac {3t}4$. 
We have 
\begin{align}
\|\int_{0}^{\frac {3t}4}\int_s^{2s}& e^{i(t-s)H} W e^{i(s-\tau')H} 
U\psi_1(2s-\tau')\,d\tau' ds\|_{\infty} \les\nn\\ 
&\int_{0}^{\frac {3t}4}\int_s^{2s} |t-s|^{-\frac 32} \|W\|_{2} \|U\|_{2}
\la 2s-\tau'\ra^{-\frac 32} \|\vpsi_{0}\|_{1\cap 2}\,d\tau' ds \nn\\ &\les 
t^{-\frac 32}
\|W\|_{2} \|U\|_{2}\int_{0}^{\frac {3t}4}\int_s^{2s}
\la 2s-\tau'\ra^{-\frac 32} \|\vpsi_{0}\|_{1\cap 2} \,d\tau' ds \les\nn\\ 
&t^{-\frac 32} \|\vpsi_{0}\|_{1\cap 2}\label{eq:tail8}
\end{align}
Now observe that in the region $(s,\tau')\in [\frac {3t}4,t]\times [s, 
2t-s]$ the value $(2s-\tau') \ge \frac t4$. Therefore, 
for the second term in \eqref{eq:tail7} we obtain that
\begin{align}
\|\int_{\frac {3t}4}^{t}\int_s^{2t-s}& e^{i(t-s)H} W e^{i(s-\tau')H} 
U\psi_1(2s-\tau')\,d\tau' ds\|_{\infty} \les\nn\\  
&\int_{\frac {3t}4}^{t}\int_s^{2t-s} |t-s|^{-\frac 32} \|W\|_{2} \|U\|_{2\cap\infty}
\la 2s-\tau'\ra^{-\frac 32} \|\vpsi_{0}\|_{1\cap 2}\,d\tau' ds \nn\\ &\les 
\la t\ra^{-\frac 32}
\|W\|_{2} \|U\|_{2}\int_{\frac {3t}4}^{t}\int_s^{2t-s} |t-s|^{-\frac 
32} \|\vpsi_{0}\|_{1\cap 2}\,d\tau' ds \les\nn\\ &\la t\ra^{-\frac 32}\int_{\frac {3t}4}^{t} |t-s|^{-\frac 
12} \|\vpsi_{0}\|_{1\cap 2}\, ds \les  t^{-\frac 32}\|\vpsi_{0}\|_{1\cap 
2}\label{eq:tail9}
\end{align}
Similarly to \eqref{eq:tail4} we split the last term in 
\eqref{eq:tail7} 
\bea
& &\int_{\frac {3t}4}^{t}\int_{2t-s}^{2s} e^{i(t-s)H} W e^{i(s-\tau')H} 
U\psi_1(2s-\tau')\,d\tau' ds = \nn\\
& &\bigg (\int_{\frac {3t}4}^{t}\int_{2t-s}^{t+\frac s2} + 
\int_{\frac {3t}4}^{t}\int_{t+\frac s2}^{2s} \bigg ) 
e^{i(t-s)H} W e^{i(s-\tau')H} 
U\psi_1(2s-\tau')\,d\tau' ds\label{eq:tail10}
\eea
Observe that in the region $(s,\tau') \in [\frac {3t}4,t]\times [2t-s, 
t+\frac s2]$ we have 
$2s-\tau'\ge \frac t8$. On the other hand, in the region 
$(s,\tau') \in [\frac {3t}4,t]\times [t+\frac s2, 2s]$ we have 
$|t-\tau'|\ge  \frac {3t}8$. 
Therefore,
\begin{align}
\|\int_{\frac {3t}4}^{t}\int_{2t-s}^{t+\frac s2} &
e^{i(t-s)H} W e^{i(s-\tau')H} 
U\psi_1(2s-\tau')\,d\tau' ds \|_{\infty}=\nn\\
&\int_{\frac {3t}4}^{t}\int_{2t-s}^{t+\frac s2}
\|e^{i(t-s)H} W  e^{-i(t-s)H}\|_{\infty\to\infty}
\|e^{i(t-\tau')H} 
U\psi_1(2s-\tau')\|_{\infty}\,d\tau' ds \les\nn\\
&\|\hat W\|_{1} \|U\|_{1\cap 2}
\int_{\frac {3t}4}^{t}\int_{2t-s}^{t+\frac s2} |t-\tau'|^{-\frac 32}
\la 2s -\tau'\ra^{-\frac 32}\|\vpsi_{0}\|_{1\cap 2}\, d\tau' ds \les\nn\\ 
&\la t\ra^{-\frac 32} \int_{\frac {3t}4}^{t}
\int_{2t-s}^{t+\frac s2} |t-\tau'|^{-\frac 32}
\|\vpsi_{0}\|_{1\cap 2}\, d\tau' ds 
\les\nn \\ &\la t\ra^{-\frac 32} \int_{\frac {3t}4}^{t} \bigg ( |t-s|^{-\frac 
12} - |s/2|^{-\frac 12}\bigg )\|\vpsi_{0}\|_{1\cap 2}\, ds 
\les t^{-\frac 32} \|\vpsi_{0}\|_{1\cap 2}\nn
\end{align}
Similarly,
\begin{align}
\|\int_{\frac {3t}4}^{t}\int_{t+\frac s2}^{2s}&
e^{i(t-s)H} W e^{i(s-\tau')H} 
U\psi_1(2s-\tau')\,d\tau' ds \|_{\infty}\les\nn\\
&\|\hat W\|_{1} \|U\|_{1\cap 2} 
\int_{\frac {3t}4}^{t}\int_{t+\frac s2}^{2s} |t-\tau'|^{-\frac 32} 
\la 2s-\tau'\ra^{-\frac 32} \|\vpsi_{0}\|_{1\cap 2}\, d\tau' ds \les\nn \\ 
&t^{-\frac 32} \int_{\frac {3t}4}^{t}\int_{t+\frac s2}^{2s} 
\la 2s-\tau'\ra^{-\frac 32} \|\vpsi_{0}\|_{1\cap 2}\, d\tau' ds \les 
t^{-\frac 32} \|\vpsi_{0}\|_{1\cap 2}\nn
\end{align}

\noindent
The conclusion from the preceding is that 
\begin{equation}
\label{eq:tail11}
\|e^{itA}P_c \vpsi_0\|_{\infty}=
\|P_ce^{itA}\vpsi_0\|_{\infty}\les |t|^{-\frac32}\|\vpsi_0\|_{1\cap2}.
\end{equation}
as desired.
\end{proof}

\begin{remark}
Dualizing \eqref{eq:tail11} yields
\[ \|e^{itA^*}P_c^* \vpsi_0\|_{2+\infty} \les |t|^{-\frac32}\|\vpsi_0\|_{1}.\]
By construction, $\Ran(P_c^*)=(\ker P_c)^\perp = \Bigl(\sum_{j=1}^k L_j\Bigr)^\perp$,
$\ker(P_c^*)=(\Ran P_c)^\perp = \sum_{j=1}^k L_j^*$, and $(P_c^*)^2=P_c^*$.
It follows from these properties that $P_c^*$ is precisely the projection onto
$\Bigl(\sum_{j=1}^k L_j\Bigr)^\perp$ corresponding to the decomposition~\eqref{eq:split}
with the roles of~$L_j$ and~$L_j^*$ interchanged. Therefore, one has the corresponding
estimate for $e^{itA}$ as well, namely
$$ 
\|e^{itA}P_c \vpsi_0\|_{2+\infty} \les |t|^{-\frac32}\|\vpsi_0\|_{1}.
$$
\end{remark}

\begin{remark} 
We have made no use of the wave operators 
$ \Omega:=\lim_{t\to\infty} e^{-itA}e^{itB}$
in this section. 
It is easy to see from the dispersive bound for the free 
evolution~$e^{itB}$ and Cook's method
that this limit exists. There is also a corresponding asymptotic 
completeness statement, namely that~$\Omega$
 is an isomorphism from $\Hil$  onto $\Ran P_c$. 
A first step in this direction is to check 
the existence of the limit $\lim_{t\to\infty} e^{-itB}e^{itA}P_c$. This 
is nontrivial, 
and can be obtained by Cook's method using 
Theorem~\ref{thm:vec_dec}.  
Then there is an additional issue of showing that
$\Ran(\Omega)=\Ran(P_c)$. For this and much more on the wave 
operators see~\cite{Cuc}.
Let us just mention that 
Lemmas~\ref{lem:ergodic} and~\ref{lem:resol} are both immediate 
consequences 
of this asymptotic completeness. 
Indeed, Lemma~\ref{lem:resol} follows from
\[ (A-z)^{-1}P_c = \Omega (B-z)^{-1} (P_c\Omega)^{-1} P_c,\]
whereas Lemma~\ref{lem:ergodic} is an immediate consequence of the 
fact that
\[
\frac{1}{T}\int_0^T e^{-i\omega t}e^{tA}P_c \, dt = \Omega 
\frac{1}{T}\int_0^T e^{-i\omega t}e^{tB}
(P_c\Omega)^{-1}P_c \, dt
\]
via the usual ergodic theorem for the unitary group $e^{tB}$. But 
since the proof 
of asymptotic completeness requires Theorem~\ref{thm:vec_dec}, we 
needed to give 
a direct proof of these lemmas. 
\end{remark}

\subsection{The case of higher dimensions}

The proof that was just presented for $n=3$ also works 
in odd dimensions $n\ge3$. We do not address the situation
in even dimensions here since it is not clear at the moment 
whether these methods can be extended to that case. 

\noindent 
There are three places where the dimension becomes relevant
in the previous proof. The first instance is Lemma~\ref{lem:free_cont},
which establishes the analytic continuation of the weighted free resolvent
$E_\eps (iB-z)^{-1}E_\eps$ from $\Re z>0$ 
to the region $\Reg_{\Gamma_{\frac{\eps}{4}}}$, together with 
the bound $(1+|z|)^{-\half}$. Since
\[ 
E_\eps (B+iz)^{-1} E_\eps = \bm E_\eps(-\Lapl+1+iz)^{-1}E_\eps & 0 \\
0 & E_\eps(\Lapl-1+iz)^{-1}E_\eps \endm, 
\]
this reduces to the same statements for the {\em scalar}
operators $(-\Lapl+1+iz)^{-1}=(-\Lapl+\zeta^2)^{-1}$ and
$(\Lapl-1+iz)^{-1} = (-\Lapl-\zeta^2)^{-1}$, where 
we have used the maps $z=U_1(\zeta)=i-i\zeta^2$, and
$z=U_2(\zeta)=-i+i\zeta^2$, respectively. As far as
the region is concerned, observe that 
$\Reg_{\Gamma_{\frac{\eps}{4}}}\subset U_1(\Reg_1\cup\Reg_2)
\cap U_2(\Reg_1^*\cup\Reg_2^*).$ 
Rauch proves the desired bound on the analytic
continuation of the free resolvent by means of
the sharp Huyghens principle in dimension $n=3$. 
The very same argument also applies to any odd 
dimension $n\ge 3$. Indeed, one writes\footnote{The analytic 
continuation of the free resovent in the region $|\zeta|<\frac 
{\epsilon}4$, which avoids the introduction of the pole at $\zeta=0$, 
can be constructed via the formula
$$
\int_0^\infty e^{-\zeta t}\; 
E_\eps\frac{\sin(t\sqrt{-\triangle})}{\sqrt{-\triangle}} E_\eps \, dt
= E_\eps (-\triangle +\zeta^2)^{-1} E_\eps .
$$
The energy estimate 
$\|\partial_{t}\frac{\sin(t\sqrt{-\triangle})}{\sqrt{-\triangle}} f 
\|_{2}\le \|f\|_{2}$ together with the sharp Huyghens principle lead
to the bound
$\|\,E_\eps\frac {\sin(t\sqrt{-\triangle})}{\sqrt{-\triangle}}E_\eps\|_{2\to2} 
\les e^{-\eps t/2}$, which is sufficient for analytic continuation.
} 
\begin{equation}
\label{eq:lapl_wave}
\int_0^\infty e^{-\zeta t}\; E_\eps\cos(t\sqrt{-\triangle})E_\eps \, dt
= \zeta E_\eps (-\triangle +\zeta^2)^{-1} E_\eps,\qquad \Re\zeta>0.
\end{equation}
If $n$ is odd, then the sharp Huyghens principle implies that
$\|E_\eps\cos(t\sqrt{-\triangle})E_\eps\|_{2\to2} \les e^{-\eps t/2}$ 
for large $t$, so that~\eqref{eq:lapl_wave} allows one
to analytically continue the right-hand side to $\Re \zeta>-\frac{\eps}{2}$
together with the desired bound. See~\cite{R} for details.

\noindent The second place where the dimension becomes
relevant is in the proof of local decay, see the term involving $C(\tau)$
in~\eqref{eq:laurent}. 
In general (odd) dimensions, there is still
the Laurent expansion
\[ G_\eps(z) = \sum_{j=-2}^\infty B_j (z-i)^{\frac{j}{2}}, \]
and similarly about the point $z=-i$. 
By our assumption on absence of eigenvalue and resonance at $i$ 
one has $B_{-2}=B_{-1}=0$. In order to obtain the desired decay estimate
on the evolution $e^{itA}$ one needs to show that the
coefficients $B_j$ of \underline{odd} powers $j$ vanish for $1\le j<n-2$.
This follows from the corresponding property of the ``free'' 
function~$F_\eps(z)$. Indeed, by the previous paragraph,
\begin{equation}
\label{eq:free_laurent}
 F_\eps(z)=\sum_{j=0}^\infty C_j\,(z-i)^{\frac{j}{2}} \text{\ \ for\ \ }|z-i|<\frac{\eps^2}{16}.
\end{equation}
where the coefficients are compact operators in $L^2$. 
Indeed, in Lemma~\ref{lem:free_cont} (and its generalization to $n\ge3$ odd)
we observed that $F_\eps(i-i\zeta^2)$ is analytic for $|\zeta|<\frac{\eps}{4}$. 
Thus, the series in~\eqref{eq:free_laurent} is absolutely convergent
in the slit region $U_1(\Reg_1)$ for small $z$. Moreover, since the rate of decay
for the free evolution $e^{itB}$ is known to be $|t|^{-\frac{n}{2}}$,  
the argument from the proof of Proposition~\ref{prop:localL2}, applied to 
the free evolution, yields that all $C_j=0$ for $1\le j<n-2$ odd. 
Finally, by the resolvent identity one has, see \eqref{eq:Geps}
\bea
 G_\eps(z) &=& \Bigl(1+F_\eps(z) \, E_\eps^{-2}V\Bigr)^{-1}\;F_\eps(z) \nn \\
&=& \Bigl(1+F_\eps(i) \, E_\eps^{-2}V\Bigr)^{-1}
\Bigl[1+(F_\eps(z) - F_\eps(i))\, E_\eps^{-2}V\Bigl(1+F_\eps(i) 
\, E_\eps^{-2}V\Bigr)^{-1}\Bigr]^{-1}\;F_\eps(z). \label{eq:neuman}
\eea
This is justified, since by our assumption of $\pm i$ being neither eigenvalue
nor resonance of~$A$ the inverse  $\Bigl(1+F_\eps(i) \, E_\eps^{-2}V\Bigr)^{-1}$
exists as a bounded operator on $L^2$, cf.~\eqref{eq:Geps_2}. 
Furthermore, the inverse in square brackets in~\eqref{eq:neuman}
exists as an $L^2$ convergent Neuman series. 
Inserting~\eqref{eq:free_laurent} into~\eqref{eq:neuman} reveals
that all coefficients of the entire right-hand side in~\eqref{eq:neuman}
corresponding to powers $(z-i)^{\frac{j}{2}}$ with $1\le j<n-2$ odd,
have to vanish. Therefore, \eqref{eq:laurent} holds with 
$\|C(\tau)\|\les |\tau-i|^{\frac{n-2}{2}}$, and~\eqref{eq:rauch} 
thus holds with $\la t\ra^{-\frac{n}{2}}$. 

\noindent
The third place where the dimension becomes relevant 
is the argument allowing for $L^2$ removal in Theorem~\ref{thm:vec_dec}.
As in the scalar case, a tedious calculation involving further Duhamel
expansion would extend the 3-d argument to higher dimensions.

\section{Matrix charge transfer model}

\begin{defi}
\label{def:chargetransm} 
By a {\em matrix charge transfer model} we mean a system
\bea
&& \frac{1}{i} \partial_t \vpsi + \bm -\Lapl & 0 \\ 0 & \Lapl \endm\vpsi + \sum^\nu_{j =1}
V_j(\cdot - \vec{v_j} t) \vpsi = 0 \label{eq:transferm} \\
&& \vpsi |_{t=0} = \vpsi_0,  \nn
\eea
where $\vec v_j$ are distinct vectors in $\R^3$, and $V_j$ are matrix 
potentials of the form
\[ V_j(t,x) = \bm U_j(x) & -e^{i\theta_j(t,x)}\,W_j(x) \\ 
e^{-i\theta_j(t,x)}\,W_j(x) & -U_j(x) \endm, \]
where $\theta_j(t,x)=(|\vec v_j\,|^2+\alpha_j^2)t+2x\cdot \vec v_j + 
\gamma_j$, $\alpha_j,\gamma_j\in\R$,
$\alpha_j\not=0$. Furthermore, we require that
each
\[ H_j = \bm -\Lapl + \half\alpha_j^2 + U_j & -W_j \\ W_j & \Lapl - 
\half\alpha_j^2 - U_j \endm \]
be admissible in the sense of Definition~\ref{def:spec_ass} and  
that it satisfy the stability
condition~\eqref{eq:stab}. 
\end{defi}

In comparison to the simple definition of a charge transform model in 
the scalar case, 
Definition~\ref{def:chargetransm} might seem unnatural. 
However, it is natural by virtue of
being the only Galilei invariant definition. 
In our context we need to use the following vector-valued Galilei 
transform (which should explain why we are using~\eqref{eq:s_2} 
rather than~\eqref{eq:s_1}):  
\begin{equation}
\label{eq:matrGal} 
\calG_{\vec v, y}(t) \binom{\psi_1}{\psi_2} := \binom{\calg_{\vec 
v,y}(t)\psi_1}{\overline{\calg_{\vec v,y}(t)\overline{\psi_2}}},
\end{equation}
where $\calg_{\vec v,y}$ are as in Section~2. 
Clearly, the transformations $\calG_{\vec v,y}(t)$ are isometries on 
all $L^p$ 
spaces. Since in our case always $y=0$, we set $\calG_{\vec v}(t):= 
\calG_{\vec v,0}(t)$. 
By~\eqref{eq:invg}, $\calG_{\vec v}(t)^{-1}=\calG_{-\vec{v}}(t)$ in 
that 
case. 
There is now the following analogue of~\eqref{eq:invg}, 
\eqref{eq:shifted}. 
In contrast to the scalar case, the transformation law of 
Lemma~\ref{lem:trans_law} also
involves a modulation~$\calM(t)$. 
It should be compared with~Definition~\ref{def:chargetrans}. 

\begin{lemma}
\label{lem:trans_law} Let $\alpha\in\R$ and set
\[ A := \bm -\Lapl + \half\alpha^2 + U & -W \\ W & \Lapl - 
\half\alpha^2 - U \endm \]
with real-valued $U,W$. Moreover, let 
$\vec v\in\R^3$,  $\theta(t,x)=(|\vec v\,|^2+\alpha^2)t+2x\cdot \vec 
v + \gamma$, $\gamma\in\R$, and define
\[ H(t) := \bm -\Lapl+U(\cdot-\vec vt) & -e^{i\theta(t,\cdot-\vec 
vt)}W(\cdot-\vec vt) \\
 e^{-i\theta(t,\cdot-\vec vt)}W(\cdot-\vec vt) & \Lapl-U(\cdot-\vec 
vt) \endm.
\]
Let $\calS(t)$, $\calS(0)=Id$, denote the propagator of the system
\[ \frac{1}{i}\partial_t \calS(t) + H(t)\calS(t) =0.\]
Finally, let 
\begin{equation}
\label{eq:M_def}
\calM(t)=\calM_{\alpha,\gamma}(t)=\bm e^{-i\omega(t)/2} & 0 \\ 0 & 
e^{i\omega(t)/2} \endm 
\end{equation}
where $\omega(t)=\alpha^2 t+\gamma$.  Then 
\begin{equation}
\label{eq:trans_law}
 \calS(t) = \calG_{\vec v}(t)^{-1} \calM(t)^{-1} e^{-itA} 
\calM(0)\calG_{\vec v}(0).
\end{equation}
\end{lemma}
\begin{proof} One has
\begin{equation}
\frac{1}{i} \partial_t \calM(t)\calG_{\vec v}(t)\calS(t) = \bm 
-\half\dot{\omega} & 0\\ 0 & \half\dot{\omega} \endm 
\calM(t)\calG_{\vec v}(t)\calS(t) 
 + \calM(t)\frac{1}{i}\dot\calG_{\vec v}(t) \calS(t) - 
\calM(t)\calG_{\vec v}(t)H(t)\calS(t). \label{eq:big_dot}
\end{equation}
Let $\rho(t,x)=t|\vec v\,|^2+2x\cdot\vec v$. One now checks the 
following properties by differentiation:
\bea
\calM(t)\frac{1}{i}\dot\calG_{\vec v}(t) &=& \bm \half|\vec 
v\,|^2+\vec v\cdot\vec p & 0 \\
0 & -\half|\vec v\,|^2+\vec v\cdot\vec p \endm \calM(t)\calG_{\vec 
v}(t) \nn \\
\calM(t)\calG_{\vec v}(t)H(t) &=& \bm -\Lapl+U & 
-e^{i(\theta-\rho-\omega)}W \\
 e^{-i(\theta-\rho-\omega)}W & \Lapl-U \endm \calM(t)\calG_{\vec 
v}(t) \nn \\
&& \quad + \bm \half|\vec v\,|^2+\vec v\cdot\vec p & 0 \\
0 & -\half|\vec v\,|^2+\vec v\cdot\vec p \endm \calM(t)\calG_{\vec 
v}(t). \label{eq:diff_com}
\eea
In view of our definitions, $\theta-\rho-\omega=0$. Since 
$\dot{\omega}=\alpha^2$, the lemma
follows by inserting~\eqref{eq:diff_com} into~\eqref{eq:big_dot}.
\end{proof}

\noindent We now return to the matrix charge transfer 
problem~\eqref{eq:transferm}.
In order to state our main theorem, we need to impose an orthogonality
condition in the context of the charge transfer, analogous to the one 
used in the scalar case.
To do so, let $P_c(H_j)$ and~$P_b(H_j)$ be the projectors as 
in~Lemma~\ref{lem:split}. 
Abusing terminology somewhat, we refer to $\Ran(P_b(H_j))$ as the 
{\em bound states} of~$H_j$. Since $P_b$ and $P_c$ are no longer
orthogonal projections, we use ``scattering states'' instead
of ``asymptotically orthogonal to the bound states'' in the following
definition. 

\begin{defi} 
\label{def:asympm}  Let $U(t) \vpsi_0 = \vpsi(t, \cdot)$ be the
solution of~\eqref{eq:transferm}. We say that $\vpsi_0$ is a
{\em scattering state} 
relative to $H_j$ if 
$$
\|P_b(H_j,t)U(t) \vpsi_0 \|_{L^2} \to 0\text{ as }t\to +\infty.
$$
Here
\begin{equation}
\label{eq:Proj2m}
P_b(H_j,t) := \calG_{\vec v_j}(t)^{-1}\calM_j(t)^{-1} P_b(H_j)\, 
\calM_j(t)\calG_{\vec v_j}(t)
\end{equation}
with $\calM_j(t)=\calM_{\alpha_j,\gamma_j}(t)$ as in \eqref{eq:M_def}.
\end{defi}

The formula~\eqref{eq:Proj2m} is of course motivated 
by~\eqref{eq:trans_law}. 
Clearly,  $P_b(H_j,t)$ is the projection onto the bound states of 
$H_j$ that have been 
translated to the
position of the matrix potential $V_j(\cdot-t\vec{v}_j)$. 
Equivalently, one 
can think of it as translating the solution of~\eqref{eq:transferm} 
from that position to the 
origin, projecting onto the bound states of~$H_j$, and then 
translating back. 

\noindent We now formulate our decay estimate for matrix charge 
transfer models.

\begin{theorem}
\label{thm:mainm}
Consider the matrix charge transfer model as in 
Definition~\ref{def:chargetransm}.
Let $U(t)$ denote the propagator of the 
equation~\eqref{eq:transferm}. Then 
for any initial data $\vpsi_0 \in L^1\cap L^2$, which
is a scattering state relative to each $H_j$ 
in the sense of Definition~\ref{def:asympm}, 
one has the decay estimates
\begin{equation}
\| U(t) \vpsi_0 \|_{L^2 + L^\infty} \lesssim
\langle t\rangle^{-\frac32}\|\vpsi_0\|_{L^1\cap L^2}. \label{eq:mainm}
\end{equation}
If in addition the matrix potentials $V_{j}$ satisfy 
$\|\hat V_{j}\|_{L^{1}}<\infty$ for all $j=1,..,\nu$ then one also
has the bound 
\begin{equation}
\| U(t) \vpsi_0 \|_{L^\infty} \lesssim
|t|^{-\frac32}\|\vpsi_0\|_{L^1\cap L^2}. \label{eq:mainm1}
\end{equation}
\end{theorem}

\begin{remark}
Theorem \ref{thm:mainm} also holds in the case of higher dimensions 
$n\ge 3$. In particular, we have the corresponding estimate 
$$
\| U(t) \vpsi_0 \|_{L^\infty} \lesssim
|t|^{-\frac n2}\|\vpsi_0\|_{L^1\cap L^2}.
$$
As in the case of {\it {scalar}} charge transfer models the proof 
of Theorem \ref{thm:mainm} mainly relies on the decay estimates for
the corresponding matrix problems $H_{j}$ with a {\it{single}} 
time-independent potential $V_{j}(x)$
\begin{equation}
\label{eq:hiogh}
\|e^{itH_{j}} P_{c}(H_{j})\,f\|_{L^{2} + L^{\infty}}\les \la t\ra^{-\frac n2} 
\|f\|_{L^1\cap L^{2}}
\end{equation}
We have established such estimates for the admissible Hamiltonians 
$H_{j}$ in Theorem \ref{thm:vec_dec} for $n=3$. In section 7.2 we 
have also observed that our method for deriving \eqref{eq:hiogh} 
can be generalized to treat the 
case of an arbitrary odd dimension $n\ge 3$. While our method leaves
the even dimensional case open at the moment, one can, in principle, 
substitute Cuccagna's dispersive estimates \cite{Cuc} instead. 
\end{remark}

\begin{remark}
\label{rem:matr_pert}
The conclusions of Theorem~\ref{thm:mainm} are also valid
for {\em perturbed} matrix
charge transfer Hamiltonians. This refers to equations
of the type 
\bea
&& \frac{1}{i} \partial_t \vpsi + \bm -\Lapl & 0 \\ 0 & \Lapl 
\endm\vpsi + \sum^\nu_{j =1}
V_j(\cdot - \vec{v_j} t) \vpsi + V_0(t,\cdot)\vpsi= 0 \nn \\
&& \vpsi |_{t=0} = \vpsi_0,  \nn
\eea
where the charge transfer part is as in Definition~\ref{def:chargetransm}
and the perturbation satisfies
\[ \sup_{t} \|V_0(t,\cdot)\|_{1\cap\infty} < \eps.\]
See
the scalar case in Section~\ref{subsec:perturb}
for an exact formulation.
\end{remark}

As in the scalar case, \eqref{eq:mainm} is proved by means of  
a bootstrap argument
with the same assumption~\eqref{eq:boot}. 
Since the argument is basically identical with the scalar case, we do 
not write it out in full detail. As in the scalar case, the estimate
\eqref{eq:mainm1} is a consequence of the $L^{2}+L^{\infty}$ bound 
\eqref{eq:mainm}. This follows by combining the machinery developed
in Proposition \ref{pr:transit} with the cancellation property for 
a matrix problem with a {\it{single}} potential established in the 
proof of the estimate \eqref{eq:vec_dec_infty} in Theorem 
\ref{thm:vec_dec}.

\noindent The fact that systems have generalized eigenspaces
rather than just eigenspaces leads to some minor
changes from the scalar case. But this really only
affects the proof of Proposition~\ref{prop:bdstates}.
We show below that the statement still remains the same.

In the following we shall assume that the
number of potentials is $\nu = 2$ and that the velocities are 
$\vec{v_1} = 0,\vec{v_2} = (1,0, \dots 0) = \vec{e_1}$. 
This can be done without loss of generality.

\subsection{Bound states}

We now estimate the rate of convergence 
of the projections onto the bound states of 
solutions which evolve from scattering states. 

\begin{prop}
\label{prop:bdstatesm}
Let $\vpsi (t, x) = (U(t) \vpsi_0)(x)$ be a solution 
of~\eqref{eq:transferm} where $\vpsi_0$ is a
scattering state relative to $H_1$ and $H_2$ in 
the sense of Definition~\ref{def:asympm}. Then 
$$
\| P_b (H_1,t) U(t)\vpsi_0 \|_{L^2} + \| P_b (H_2,t) U(t) 
\vpsi_0\|_{L^2}
\lesssim e^{-\alpha t} \| \vpsi_0\|_{L^2}
$$
for some $\alpha>0$.
\end{prop}
\begin{proof}
By symmetry it suffices to prove the bound on the first part (one can 
switch the
roles of $V_1$ and~$V_2$ by means of a Galilei transform and a 
modulation). In view of Lemma~\ref{lem:split} one decomposes
\begin{equation}
\vphi(t):=\calM_{\alpha_1,\gamma_1}(t)U(t) \vpsi_0 = \sum^M_{j=0}  
\vf_j(t) + \vphi_1 (t, \cdot)\label{eq:decompm}
\end{equation}
relative to $H_1$ so that $\vphi_1(t,\cdot)$ lies in the continuous
subspace of $H_1$, i.e.,  $P_c(H_1) \vphi_1 = \vphi_1$ and 
$P_b(H_1)\vphi_1 = 0$ for all times. 
Furthermore, $\vf_j(t)\in L_j$ where $L_j$ are the spaces from 
Lemma~\ref{lem:split}
for the operator~$H_1$. By assumption, 
\begin{equation}
\label{eq:inf_dec}
\sum^M_{j=0} |\vf_j(t)|^2 \to 0 \text{ as } t \to \infty.
\end{equation}
One checks that $\vphi$ satisfies the equation
\begin{equation}
\label{eq:mod_eq}
 \frac{1}{i}\partial_t\vphi+ H_1 \vphi + \tilde{V}_2(\cdot-\vec e_1 
t)\vphi = 0
\end{equation}
where $\tilde{V}_2 = \calM_{\alpha_1,\gamma_1}(t) V_2 
\calM_{\alpha_1,\gamma_1}(t)^{-1}$. 
Substituting \eqref{eq:decompm} into \eqref{eq:mod_eq} yields
\begin{align}
&\frac{1}{i}\partial_t \vphi + H_1 \vphi + \tilde{V}_2(\cdot 
-t\vec{e_1}) \vphi\nn \\
& = \sum^M_{j=0}\frac{1}{i} \dot{\vf}_j + \frac{1}{i}\dot{\vphi}_1(t) 
+ 
+\sum_{j=0}^M H_1 \vf_j(t) + H_1\vphi_1 +  \tilde{V}_2(\cdot - 
t\vec{e_1}) \vphi =0 
\label{eq:aim}
\end{align}
Since $P_c(H_1) \vphi_1 = \vphi_1$ and $\Ran(P_c(H_1))$ is closed by 
assumption, one has 
$$
H_1\vphi_1 =  P_c(H_1) H_1 \vphi_1,\qquad 
\partial_t \vphi_1 = P_c(H_1) \partial_t
\vphi_1.
$$
In particular 
$$P_b(H_1) \left(\frac{1}{i}\partial_t \vphi_1 + H_1 \vphi_1\right) 
=0.
$$
Therefore, applying $P_b(H_1)$ to~\eqref{eq:aim} yields
\[ 
\sum_{j=0}^M \Bigl[\frac{1}{i}\dot{\vf}_j(t)  + H_1 \vf_j(t)\Bigr] = 
\vg(t) := 
-P_b(H_1)\tilde{V}_2(\cdot - t\vec{e_1}) \vphi.
\]
In view of the explicit expression \eqref{eq:Pbexpl} for $P_b(H_1)$ 
one has $|\vg(t)|\les e^{-\eps t}$. 
Moreover, decomposing $\vg(t)=\sum_{j=0}^M \vg_j(t)$ with 
$\vg_j(t)\in L_i$, one
also has $|\vg_j(t)|\les e^{-\eps t}$ for all $0\le j\le M$. 
Hence,
\bea
 \frac{1}{i}\dot{\vf}_0 + H_1 \vf_0 &=& \vg_0 \label{eq:f0} \\
 \frac{1}{i}\dot{\vf}_j + \omega_j  \vf_j &=& \vg_j \text{\ \ for\ \ 
}1\le j\le M\nn
\eea
with exponentially decaying right-hand sides.
It follows from \eqref{eq:inf_dec} that $|\vf_j(t)|\les e^{-\eps t}$ 
for $1\le j\le M$. Moreover, 
applying~$H_1$ to~\eqref{eq:f0} one obtains
\[ (H_1 \vf_0)^{\cdot} = i H_1 \vg_0(t). \]
Since the right-hand side decays exponentially, and since 
$H_1\vf_0(t)\to0$ as~$t\to\infty$
by our assumption of starting from a scattering state, 
it follows that $H_1 \vf_0(t)$ decays exponentially,
and therefore also~$\vf_0(t)$. 
\end{proof}

\section {Generalized decay estimates for the charge transfer model}
\label{se:Sob}
Consider the time-dependent matrix charge transfer problem 
$$
i \pa_{t}\vpsi + H(\sigma,t) \vpsi = F
$$
where the matrix charge transfer Hamiltonian $H(\sigma,t)$ 
is of the form
$$
H(\sigma,t)= \bm \Lapl & 0 \\ 0 & -\Lapl 
\endm+ \sum^\nu_{j =1}
V_j(\cdot - \vec{v_j} t) 
$$
where $\vec v_j$ are distinct vectors in $\R^3$, and $V_j$ are matrix 
potentials of the form
\[ V_j(t,x) = \bm U_j(x) & -e^{i\theta_j(t,x)}\,W_j(x) \\ 
e^{-i\theta_j(t,x)}\,W_j(x) & -U_j(x) \endm, \]
where $\theta_j(t,x)=(|\vec v_j\,|^2+\alpha_j^2)t+2x\cdot \vec v_j + 
\gamma_j$, $\alpha_j,\gamma_j\in\R$,
$\alpha_j\not=0$. 
Our goal is to extend the dispersive estimate 
\be 
\label{eq:first}
\|\vpsi(t)\|_{L^{2}+L^{\infty}}\les (1+t)^{-\frac n2} 
\Big (\|\vpsi_{0}\|_{L^{1}\cap L^{2}}  + \||F\||| + B\Big )
\ee
with
\be
\label{eq:tripleF}
\||F\||:=  \sup_{t\ge 0} \int_{0}^{t} 
\|F(\tau)\|_{L^{1}}\,d\tau 
+ (1+t)^{\frac n2+1} 
\|F(t,\cdot)\|_{L^{2}}
\ee
to the corresponding estimates for the derivatives of 
$\vpsi(t)$.
The estimate \eqref{eq:first} holds only
for the solutions $\vpsi(t)$ which are scattering states, i.e.,
for $\vpsi$ obeying the a priori condition that
$$
\|P_{b}(H_{j},t) \vpsi(t)\|_{L^{2}} \le B(1+t)^{-\frac n2}
$$
for all $j=1,..,\nu$.
Our first lemma shows that the functions 
$$
\vpsi_{k}(t):= \nabla^{k}\vpsi(t),\qquad k\in Z_{+}
$$
are scattering states as well.
\begin{lemma}
The functions $\vpsi_{k}(t)$ obey the estimates 
\be
\|P_{b}(H_{j},t) \vpsi_{k}(t)\|_{L^{2}} \les C_{k}(1+t)^{-\frac n2}
\label{eq:secbound}
\ee
\label{le:boundsta}
\end{lemma}
\begin{proof}
Let $\veta(x)$ be an arbitrary $C^{\infty}$ exponentially localized
function. Then for any $y\in \R^{n}$
$$
\int_{\R^{n}} \vpsi_{k}(t,x)\cdot \veta (x-y)\, dx = 
(-1)^{k} \int_{\R^{n}} \vpsi(t,x) \nabla^{k} \veta\,(x-y)\,dx
\les \|\vpsi(t)\|_{L^{2}+L^{\infty}} \|\nabla^{k} \veta\,\|_{L^{1}\cap 
L^{2}} \les (1+t)^{-\frac n2} 
$$
Now recall that the projection 
$$
P_b(H_j,t) := \calG_{\vec v_j}(t)^{-1}\calM_j(t)^{-1} P_b(H_j)\, 
\calM_j(t)\calG_{\vec v_j}(t)
$$
with the $P_{b}(H_{j}$ is given exlicitly 
$$
P_{b}(H_{j}) f = \sum_{\alpha\beta} c_{\alpha\beta} u_{\alpha} 
(f,v_{\beta} )
$$
where $c_{\alpha\beta}$ are given constants and $u_{\alpha}, 
v_{\beta}$ are exponentially localized functions. The result 
now follows.
\end{proof}
\begin{prop}
The functions $\vpsi_{k}=\nabla^{k}\vpsi$ satisfy the 
$L^{2}+L^{\infty}$
dispersive estimate
\be 
\label{eq:second}
\|\nabla^{k}\vpsi(t)\|_{L^{2}+L^{\infty}}\les (1+t)^{-\frac n2} 
\sum_{\ell=0}^{k} \Big (\|\nabla^{\ell}\vpsi_{0}\|_{L^{1}\cap L^{2}}  + 
\||\nabla^{\ell} F\|| + B\Big )
\ee
\label{prop:twoinf}
\end{prop}
\begin{proof}
We have already shown that $\nabla^{k}\psi$ is a scattering
state. Moreover, differentiating the equation $k$ times we obtain 
$$
i \pa_{t} \nabla^{k}\vpsi + H(t,\sigma) \nabla^{k}\vpsi = 
F_{k}:=\sum_{\ell=0}^{k-1} G_{\ell}(t,x) \nabla^{\ell}\vpsi + 
\nabla^{k} F 
$$
where $G_{\ell}(t,x)$ are smooth exponentially localized 
potentials uniformly bounded in time.
Therefore $\nabla^{k}$ is a scattering state solving an 
inhomogeneous charge transfer problem. Using the estimate 
\eqref{eq:first} we then have
\be
\label{eq:Fk}
\|\nabla^{k}\vpsi(t)\|_{L^{2}+L^{\infty}}\les (1+t)^{-\frac n2} 
\Big (\|\nabla^{k}\vpsi_{0}\|_{L^{1}\cap L^{2}}  +
\||F_{k}(\tau)\|| + B\Big )
\ee
We use that for any $p\in [1,2]$
$$
\| G_{\ell}(t,x) \nabla^{\ell}\vpsi \|_{L^{p}} \les 
\|\nabla^{\ell} \vpsi\|_{L^{2}+L^{\infty}}
$$
Proceeding by induction on $k$ we conclude that for 
any $\ell <k$
\begin{align*}
\int_{0}^{t}
\| G_{\ell}(t,x) \nabla^{\ell}\vpsi \|_{L^{1}} \les &
\int_{0}^{t} (1+\tau)^{-\frac n2}\,d\tau 
\sum_{m=0}^{\ell} \Big (\|\nabla^{m} \vpsi_{0}\|_{L^{1}\cap L^{2}}  + 
 \||\nabla^{m} F\||\Big )\\ \les & 
\sum_{m=0}^{\ell} (\|\nabla^{m} \vpsi_{0}\|_{L^{1}\cap L^{2}}  + 
 \||\nabla^{m} F\||)
\end{align*}
and that 
$$
(1+t)^{\frac n2}
\| G_{\ell}(t,x) \nabla^{\ell}\vpsi (t)\|_{L^{2}} \les 
\sum_{m=0}^{\ell} (\|\nabla^{m} \vpsi_{0}\|_{L^{1}\cap L^{2}}  + 
 \||\nabla^{m} F\||)
$$
The result now follows from \eqref{eq:Fk} and 
the inequality
$$
\||F_{k}(\tau)\||\le \||\nabla^{k} F\|| + 
\sum_{\ell=0}^{k-1}\|| G_{\ell}(t,x) \nabla^{\ell}\vpsi \||
$$
\end{proof}
Consider the following Banach spaces $\Xl$ and $\Yl$ of functions of 
$(t,x)$:
\bea
\|\psi\|_{\Xl_{s}} &=&  \sup_{t\ge 0}\Big (\|\psi(t,\cdot)\|_{H^{s}} +  
(1+t)^{\frac n2}\sum_{k=0}^{s}
\|\nabla^{k}\psi(t,\cdot)\|_{L^{2}+L^\infty}\Big )\label{eq:Xls}\\
\|F\|_{\Yl_{s}} &=& \sup_{t\ge 0}\sum_{k=0}^{s}
\Big (\int_{0}^{t} \|\nabla^{k}F(\tau,\cdot)\|_{L^{1}}\,d\tau +
(1+t)^{\frac n2+1}\|\nabla^{k} F(t,\cdot)\|_{L^{2}}\Big ).  \label{eq:Yls}
\eea
These spaces are relevant for our companion paper~\cite{RSS}. 
We can summarize our estimates for the charge transfer model in the
following proposition.
\begin{prop}
Let $\vpsi$ be a solution of the matrix charge transfer problem
$$
i\pa_{t}\vpsi + H(t,\sigma)\vpsi = F
$$
satisfying the condition that for every $j=1,..,\nu$
\be
\label{eq:orth}
\|P_{b}(H_{j}(\sigma,t)) \vpsi\|_{L^{2}}\les B (1+t)^{-\frac n2}
\ee 
Then for any integer $s\ge 0$
\be
\label{eq:estv}
\|\vpsi\|_{\Xl_{s}}\les \sum_{k=0}^{s}\|\nabla^{k} 
\psi(0,\cdot)\|_{L^{1}\cap L^{2}} + \|F\|_{\Yl_{s}} + B
\ee
\label{prop:Finest}
\end{prop}

\bibliographystyle{amsplain}

\noindent
\textsc{Rodnianski: Department of Mathematics, Princeton University, 
Fine Hall, Princeton N.J. 08544, U.S.A.}\\
{\em email: }\textsf{\bf irod@math.princeton.edu}

\medskip\noindent
\textsc{Schlag: Division of Astronomy, Mathematics, and Physics, 
253-37 Caltech, Pasadena, C.A. 91125, U.S.A.}\\
{\em email: }\textsf{\bf schlag@its.caltech.edu}

\medskip\noindent
\textsc{Soffer: Mathematics Department, Rutgers University, New 
Brunswick, N.J. 08903, U.S.A.}\\
{\em email: }\textsf{\bf soffer@math.rutgers.edu}

\end{document}